\documentclass[10pt,a4paper,twoside,fleqn]{article}
\usepackage{amsfonts,amssymb,amsmath,amsthm,latexsym,mathrsfs}
\usepackage{a4wide,fancyhdr}
\usepackage{graphicx,epic,eepic}
\usepackage{titlesec}
\titleformat{\section}{\large\bfseries}{\thesection}{1em}{}

\numberwithin{equation}{section}

\newtheorem{theorem}{Theorem}[section]
\newtheorem{lemma}[theorem]{Lemma}
\newtheorem{proposition}[theorem]{Proposition}
\newtheorem{corollary}[theorem]{Corollary}

\theoremstyle{definition}

\newcommand{\mres}{\mathbin{\vrule height 1.6ex depth 0pt width
0.13ex\vrule height 0.13ex depth 0pt width 1.3ex}}

\def\XXint#1#2#3{{\setbox0=\hbox{$#1{#2#3}{\int}$ } 
\vcenter{\hbox{$#2#3$ }}\kern-.6\wd0}}

\begin{document}

\begin{center}
{\bf\large An isoperimetric inequality in the plane with a log-convex density}
\end{center}

\smallskip

\begin{center} 
I. McGillivray
\end{center}

\smallskip

\begin{center}
School of Mathematics, University of Bristol\\
University Walk, Bristol BS8 1TW, UK\\
maiemg@bristol.ac.uk
\end{center}

\smallskip

\begin{abstract}
\noindent Given a positive lower semi-continuous density  $f$ on $\mathbb{R}^2$ the weighted volume $V_f:=f\mathscr{L}^2$ is defined on the $\mathscr{L}^2$-measurable sets in $\mathbb{R}^2$. The $f$-weighted perimeter of a set of finite perimeter $E$ in $\mathbb{R}^2$ is written $P_f(E)$. We study minimisers for the weighted isoperimetric problem
\[
I_f(v):=\inf\Big\{ P_f(E):E\text{ is a set of finite perimeter in }\mathbb{R}^2\text{ and }V_f(E)=v\Big\}
\]
for $v>0$. Suppose $f$ takes the form $f:\mathbb{R}^2\rightarrow(0,+\infty);x\mapsto e^{h(|x|)}$ where $h:[0,+\infty)\rightarrow\mathbb{R}$ is a non-decreasing convex function. Let $v>0$ and $B$ a centred ball in $\mathbb{R}^2$ with $V_f(B)=v$. We show that $B$ is a minimiser for the above variational problem and obtain a uniqueness result. 
\end{abstract} 

\smallskip


\smallskip


\section{Introduction}

\smallskip

\noindent Let $f$ be a positive lower semi-continuous density on $\mathbb{R}^2$. The weighted volume $V_f:=f\mathscr{L}^2$ is defined on the $\mathscr{L}^2$-measurable sets in $\mathbb{R}^2$. Let $E$ be a set of finite perimeter in $\mathbb{R}^2$. The weighted perimeter of $E$ is defined by
\begin{equation}\label{weighted_perimeter}
P_f(E):=\int_{\mathbb{R}^2}f\,d|D\chi_E|\in[0,+\infty].
\end{equation}
We study minimisers for the weighted isoperimetric problem
\begin{equation}\label{isoperimetric_problem}
I_f(v):=\inf\Big\{ P_f(E):E\text{ is a set of finite perimeter in }\mathbb{R}^2\text{ and }V_f(E)=v\Big\}
\end{equation}
for $v>0$. To be more specific we suppose that $f$ takes the form
\begin{equation}\label{form_of_f}
f:\mathbb{R}^2\rightarrow(0,+\infty);x\mapsto e^{h(|x|)}
\end{equation}
where $h:[0,+\infty)\rightarrow\mathbb{R}$ is a non-decreasing convex function. Our first main result is the following. It contains the classical isoperimetric inequality (cf. \cite{deGiorgi1958}, \cite{Fusco2004}) as a special case; namely, when $h$ is constant on $[0,+\infty)$. 

\smallskip

\begin{theorem}\label{main_theorem}
Let $f$ be as in (\ref{form_of_f}) where $h:[0,+\infty)\rightarrow\mathbb{R}$ is a non-decreasing convex function. Let $v>0$ and $B$ a centred ball in $\mathbb{R}^2$ with $V_f(B)=v$. Then $B$ is a minimiser for (\ref{isoperimetric_problem}).
\end{theorem}

\smallskip

\noindent 
\noindent For $x\geq 0$ and $v\geq 0$ define the directional derivative of $h$ in direction $v$ by
\[
h^\prime_+(x,v):=\lim_{t\downarrow 0}\frac{h(x+tv)-h(x)}{t}\in\mathbb{R}
\]
and define $h^\prime_-(x,v)$ similarly for $x>0$ and $v\leq 0$. We introduce the notation
\[
\varrho_+:=h^\prime_+(\cdot,+1),
\varrho_-:=-h^\prime_+(\cdot,-1)
\text{ and }
\varrho:=(1/2)(\varrho_++\varrho_-)
\]
on $(0,+\infty)$. The function $h$ is locally of bounded variation and is differentiable a.e. with $h^\prime=\varrho$ a.e. on $(0,+\infty)$. Our second main result is a uniqueness theorem. 

\smallskip

\begin{theorem}\label{uniqueness_theorem}
Let $f$ be as in (\ref{form_of_f}) where $h:[0,+\infty)\rightarrow\mathbb{R}$ is a non-decreasing convex function. Suppose that $R:=\inf\{\varrho>0\}\in[0,+\infty)$ and set $v_0:=V(B(0,R))$. Let $v>0$ and $E$ a minimiser for (\ref{isoperimetric_problem}). The following hold:
\begin{itemize}
\item[(i)] if $v\leq v_0$ then $E$ is a.e. equivalent to a ball $B$ in $\overline{B}(0,R)$ with $V(B)=V(E)$;
\item[(ii)] if $v>v_0$ then $E$ is a.e. equivalent to a centred ball $B$ with $V(B)=V(E)$.
\end{itemize}
\end{theorem}

\smallskip

\noindent Theorem \ref{main_theorem} is a generalisation of Conjecture 3.12 in \cite{Rosaleset2006} (due to K. Brakke) in the sense that less regularity is required of the density $f$: in the latter, $h$ is supposed to be smooth on $(0,+\infty)$ as well as convex and non-decreasing. This conjecture springs in part from the observation that the weighted perimeter of a local volume-preserving perturbation of a centred ball is non-decreasing (\cite{Rosaleset2006} Theorem 3.10).  In addition, the conjecture holds for log-convex Gaussian densities of the form $h:[0,+\infty)\rightarrow\mathbb{R};t\mapsto e^{ct^2}$ with $c>0$ (\cite{Borell1986}, \cite{Rosaleset2006} Theorem 5.2). In subsequent work partial forms of the conjecture were proved in the literature. In  \cite{KolesnikovZhdanov2011} it is shown to hold for large $v$ provided that $h$ is uniformly convex in the sense that $h^{\prime\prime}\geq 1$ on $(0,+\infty)$ (see \cite{KolesnikovZhdanov2011} Corollary 6.8). A complementary result is contained in \cite{FigalliMaggi2013} Theorem 1.1 which establishes the conjecture for small $v$ on condition that $h^{\prime\prime}$ is locally uniformly bounded away from zero on $[0,+\infty)$. The above-mentioned conjecture is proved in large part in \cite{Chambers2013} (see Theorem 1.1) in dimension $n\geq 2$ (see also \cite{Boyeretal2016}). There it is assumed that the function $h$ is of class $C^3$ on $(0,+\infty)$ and is convex and even (meaning that $h$ is the restriction of an even function on $\mathbb{R}$ to $[0,+\infty)$). A uniqueness result is also obtained (\cite{Chambers2013} Theorem 1.2). We obtain these results under weaker hypotheses in the $2$-dimensional case and our proofs proceed along different lines. 

\smallskip

\noindent We give a brief outline of the article. In Section 2 we discuss some preliminary material. In Section 3 we show that (\ref{isoperimetric_problem}) admits an open minimiser $E$ with $C^1$ boundary $M$ (Theorem \ref{C1_property_of_reduced_boundary}). The argument draws upon the regularity theory for almost minimal sets (cf. \cite{Tamanini1984}) and includes an adaptation of \cite{Morgan2003} Proposition 3.1. In Section 4 it is shown that the boundary $M$ is of class $C^{1,1}$ (and has weakly bounded curvature). This result is contained in \cite{Morgan2003} Corollary 3.7 (see also \cite{CintiPratelli2015}) but we include a proof for completeness. This Section also includes the result that $E$ may be supposed to possess spherical cap symmetry (Theorem \ref{M_is_spherical_cap_symmetric}). Section 5 contains further results on spherical cap symmetric sets useful in the sequel. The main result of Section 6 is Theorem \ref{constant_weighted_mean_curvature} which shows that the generalised (mean) curvature is conserved along $M$ in a weak sense. In Section 7 it is shown that there exist convex minimisers of (\ref{isoperimetric_problem}). Sections 8 and 9 comprise an analytic interlude and are devoted to the study of solutions of the first-order differential equation that appears in Theorem \ref{ode_for_y} subject to Dirichlet boundary conditions. Section 9 for example contains a comparison theorem for solutions to a Ricatti equation (Theorem \ref{distribution_function_inequality_for_Ricatti_equation} and Corollary \ref{integral_of_w}). These are new as far as the author is aware. Section 10 concludes the proof of our main theorems.
\section{Some preliminaries}

\smallskip

\noindent{\em Geometric measure theory.} We use $|\cdot|$ to signify the Lebesgue measure on $\mathbb{R}^2$ (or occasionally $\mathscr{L}^2$).  Let $E$ be a $\mathscr{L}^2$-measurable set in $\mathbb{R}^2$. The set of points in $E$ with density $t\in[0,1]$ is given by
\[
E^t:=\left\{x\in\mathbb{R}^2:\,\lim_{\rho\downarrow 0}\frac{|E\cap B(x,\rho)|}{|B(x,\rho)|}=t\right\}.
\]
As usual $B(x,\rho)$ denotes the open ball in $\mathbb{R}^2$ with centre $x\in\mathbb{R}^2$ and radius $\varrho>0$. The set $E^1$ is the measure-theoretic interior of $E$ while $E^0$ is the measure-theoretic exterior of $E$. The essential boundary of $E$ is the set $\partial^\star E:=\mathbb{R}^2\setminus(E^0\cup E^1)$. 

\smallskip

\noindent Recall that an integrable function $u$ on $\mathbb{R}^2$ is said to have bounded variation if the distributional derivative of $u$ is representable by a finite Radon measure $Du$ (cf. \cite{Ambrosio2000} Definition 3.1 for example) with total variation $|Du|$; in this case, we write $u\in\mathrm{BV}(\mathbb{R}^2)$. The set $E$ has finite perimeter if $\chi_E$ belongs to $\mathrm{BV}_{\mathrm{loc}}(\mathbb{R}^2)$. The reduced boundary $\mathscr{F}E$ of $E$ is defined by
\[
\mathscr{F}E:=\Big\{x\in\mathrm{supp}|D\chi_E|:\nu^E(x):=\lim_{\rho\downarrow 0}\frac{D\chi_E(B(x,\rho))}{|D\chi_E|(B(x,\rho))}\text{ exists in }\mathbb{R}^2\text{ and }|\nu^E(x)|=1\Big\}
\]
(cf. \cite{Ambrosio2000} Definition 3.54) and is a Borel set (cf. \cite{Ambrosio2000} Theorem 2.22 for example). We use $\mathscr{H}^k$ ($k\in[0,+\infty)$) to stand for $k$-dimensional Hausdorff measure. If $E$ is a set of finite perimeter in $\mathbb{R}^2$ then
\begin{equation}\label{relation_between_boundaries}
\mathscr{F}E\subset E^{1/2}\subset\partial^* E
\text{ and }
\mathscr{H}^1(\partial^* E\setminus\mathscr{F}E)=0
\end{equation}
by \cite{Ambrosio2000} Theorem 3.61.

\smallskip

\noindent Let $f$ be a positive locally Lipschitz density on $\mathbb{R}^2$. Let $E$ be a set of finite perimeter and $U$ a bounded open set in $\mathbb{R}^2$. The weighted perimeter of $E$ relative to $U$ is defined by
\[
P_f(E,U):=\sup\Big\{\int_U\mathrm{div}(fX)\,dx:X\in C^\infty_c(U,\mathbb{R}^2),\|X\|_\infty\leq 1\Big\}.
\]
By the Gauss-Green formula (\cite{Ambrosio2000} Theorem 3.36 for example) and a convolution argument,
\begin{align}
P_f(E,U)&=\sup\Big\{\int_{\mathbb{R}^2}f\langle\nu^E,X\rangle\,d|D\chi_E|:
X\in C^\infty_c(\mathbb{R}^2,\mathbb{R}^2),\mathrm{supp}[X]\subset U,\|X\|_\infty\leq 1\Big\}\nonumber\\
&=\sup\Big\{\int_{\mathbb{R}^2}f\langle\nu^E,X\rangle\,d|D\chi_E|:
X\in C_c(\mathbb{R}^2,\mathbb{R}^2),\mathrm{supp}[X]\subset U,\|X\|_\infty\leq 1\Big\}\nonumber\\
&=\int_U f\,d|D\chi_E|\label{weighted_perimeter_of_E_relative_to_U}
\end{align}
where we have also used \cite{Ambrosio2000} Propositions 1.47 and 1.23.  

\smallskip

\begin{lemma}\label{diffeomorhism_invariance_of_boundary}
Let $\varphi$ be a $C^1$ diffeomeorphism of $\mathbb{R}^2$ which coincides with the identity map on the complement of a compact set and $E\subset\mathbb{R}^2$ with $\chi_E\in\mathrm{BV}(\mathbb{R}^2)$. Then 
\begin{itemize}
\item[(i)] $\chi_{\varphi(E)}\in\mathrm{BV}(\mathbb{R}^2)$;
\item[(ii)] $\partial^\star\varphi(E)=\varphi(\partial^\star E)$;
\item[(iii)] $\mathscr{H}^1(\mathscr{F}\varphi(E)\Delta\varphi(\mathscr{F}E))=0$.
\end{itemize}
\end{lemma}

\smallskip

\noindent{\em Proof.}
Part {\em (i)} follows from \cite{Ambrosio2000} Theorem 3.16 as $\varphi$ is a proper Lipschitz function. Given $x\in E^0$ we claim that $y:=\varphi(x)\in\varphi(E)^0$. 
Let $M$ stand for the Lipschitz constant of $\varphi$ and $L$ stand for the Lipschitz constant of $\varphi^{-1}$. Note that
$
B(y,r)\subset\varphi(B(x,Lr))
$
for each $r>0$. As $\varphi$ is a bijection and using \cite{Ambrosio2000} Proposition 2.49,
\begin{align}
|\varphi(E)\cap B(y,r)|&\leq|\varphi(E)\cap\varphi(B(x,Lr)|
=|\varphi(E\cap B(x,Lr))|\leq M^2|E\cap B(x,Lr)|.\nonumber
\end{align}
This means that
\[
\frac{|\varphi(E)\cap B(y,r)|}{|B(y,r)|}
\leq(LM)^2
\frac{|E\cap B(x,Lr)|}{|B(x,Lr)|}
\]
for $r>0$ and this proves the claim. This entails that $\varphi(E^0)\subset[\varphi(E)]^0$. The reverse inclusion can be seen using the fact that $\varphi$ is a bijection. In summary $\varphi(E^0)=[\varphi(E)]^0$. The corresponding identity for $E^1$ can be seen in a similar way. These identities entail {\em (ii)}. From (\ref{relation_between_boundaries}) and {\em (ii)} we may write $\mathscr{F}\varphi(E)\cup N_1=\varphi(\mathscr{F}E)\cup\varphi(N_2)$ for $\mathscr{H}^1$-null sets $N_1$, $N_2$ in $\mathbb{R}^2$. Item {\em (iii)} follows.
\qed

\smallskip

\noindent{\em Curves with weakly bounded curvature.} Suppose the open set $E$ in $\mathbb{R}^2$ has $C^1$ boundary $M$. Denote by $n:M\rightarrow\mathbb{S}^1$ the inner unit normal vector field. Given $p\in M$ we choose a tangent vector $t(p)\in\mathbb{S}^1$ in such a way that the pair $\{t(p),n(p)\}$ forms a positively oriented basis for $\mathbb{R}^2$. There exists a local parametrisation $\gamma_1:I\rightarrow M$ where $I=(-\delta,\delta)$ for some $\delta>0$ of class $C^1$ with $\gamma_1(0)=p$. We always assume that $\gamma_1$ is parametrised by arc-length and that $\dot{\gamma_1}(0)=t(p)$ where the dot signifies differentiation with respect to arc-length. Let $X$ be a vector field defined in some neighbourhood of $p$ in $M$. Then
\begin{align}
(D_{t}X)(p)&:=\frac{d}{ds}\Big\rvert_{s=0}(X\circ\gamma_1)(s)\label{directional_derivative}
\end{align}
if this limit exists and the divergence $\mathrm{div}^M X$ of $X$ along $M$ at $p$ is defined by
\begin{align}
\mathrm{div}^M X&:=\langle D_{t}X,t\rangle\label{divergence}
\end{align}
evaluated at $p$. Suppose that $X$ is a vector field in $C^1(U,\mathbb{R}^2)$ where $U$ is an open neighbourhood of $p$ in $\mathbb{R}^2$. Then  
\begin{align}
\mathrm{div}\,X=\mathrm{div}^M X+\langle D_nX,n\rangle\label{decomposition_of_divergence}
\end{align}
at $p$. If $p\in M\setminus\{0\}$ let $\sigma(p)$ stand for the angle measured anti-clockwise from the position vector $p$ to the tangent vector $t(p)$; $\sigma(p)$ is uniquely determined up to integer multiples of $2\pi$.

\smallskip 

\noindent Let $E$ be an open set in $\mathbb{R}^2$ with $C^{1,1}$ boundary $M$. Let $x\in M$ and $\gamma_1:I\rightarrow M$ a local parametrisation of $M$ in a neighbourhood of $x$. There exists a constant $c>0$ such that
\[
|\dot{\gamma}_1(s_2)-\dot{\gamma}_1(s_1)|\leq c|s_2-s_1|
\]
for $s_1,s_2\in I$; a constraint on average curvature (cf. \cite{Dubins1957}, \cite{HowardTreibergs1995}). That is, $\dot{\gamma}_1$ is Lipschitz on $I$. So $\dot{\gamma}_1$ is absolutely continuous and differentiable a.e. on $I$ with
\begin{align}
\dot{\gamma}_1(s_2)-\dot{\gamma}_1(s_1)&=\int_{s_1}^{s_2}\ddot{\gamma}_1\,ds\label{difference_of_dot_gamma}
\end{align}
for any $s_1,s_2\in I$ with $s_1<s_2$. Moreover, $|\ddot{\gamma}_1|\leq c$ a.e. on $I$ (cf. \cite{Ambrosio2000} Corollary 2.23). As $\langle\dot{\gamma}_1,\dot{\gamma}_1\rangle=1$ on $I$ we see that $\langle\dot{\gamma}_1,\ddot{\gamma}_1\rangle=0$ a.e. on $I$. The (geodesic) curvature $k_1$ is then defined a.e. on $I$ via the relation
 \begin{equation}
\ddot{\gamma}_1=k_1n_1
\end{equation}
as in \cite{HowardTreibergs1995}. The curvature $k$ of $M$ is defined $\mathscr{H}^1$-a.e. on $M$ by
 \begin{equation}
 k(x):=k_1(s)
 \end{equation}
 whenever $x=\gamma_1(s)$ for some $s\in I$ and $k_1(s)$ exists. We sometimes write $H(\cdot,E)=k$.

\smallskip 

\noindent Let $E$ be an open set in $\mathbb{R}^2$ with $C^{1}$ boundary $M$. Let $x\in M$ and $\gamma_1:I\rightarrow M$ a local parametrisation of $M$ in a neighbourhood of $x$.  In case $\gamma_1\neq 0$ let $\theta_1$ stand for the angle measured anti-clockwise from $e_1$ to the position vector $\gamma_1$ and $\sigma_1$ stand for the angle measured anti-clockwise from the position vector $\gamma_1$ to the tangent vector $t_1=\dot{\gamma}_1$. Put $r_1:=|\gamma_1|$ on $I$. Then $r_1,\theta_1\in C^1(I)$ and
\begin{align}
\dot{r}_1&=\cos\sigma_1;\label{cos_of_sigma}\\
r_1\dot{\theta}_1&=\sin\sigma_1;\label{sin_of_sigma}
\end{align}
on $I$ provided that $\gamma_1\neq 0$. Now suppose that $M$ is of class $C^{1,1}$. Let $\alpha_1$ stand for the angle measured anti-clockwise from the fixed vector $e_1$ to the tangent vector $t_1$ (uniquely determined up to integer multiples of $2\pi$). Then $t_1=(\cos\alpha_1,\sin\alpha_1)$ on $I$ so $\alpha_1$ is absolutely continuous on $I$. In particular, $\alpha_1$ is differentiable a.e. on $I$ with $\dot{\alpha}_1=k_1$ a.e. on $I$. This means that $\alpha_1\in C^{0,1}(I)$. In virtue of the identities $r_1\cos\sigma_1=\langle\gamma_1,t_1\rangle$ and $r_1\sin\sigma_1=-\langle\gamma_1,n_1\rangle$ we see that $\sigma_1$ is absolutely continuous on $I$ and $\sigma_1\in C^{0,1}(I)$. By choosing an appropriate branch we may assume that
\begin{align}
\alpha_1&=\theta_1+\sigma_1\label{alpha_1_as_sum}
\end{align}
on $I$. We may  choose $\sigma$ in such a way that $\sigma\circ\gamma_1=\sigma_1$ on $I$.

\smallskip

\noindent{\em Flows.} Recall that a diffeomorphism $\varphi:\mathbb{R}^2\rightarrow\mathbb{R}^2$ is said to be proper if $\varphi^{-1}(K)$ is compact whenever $K\subset\mathbb{R}^2$ is compact. Given $X\in C^\infty_c(\mathbb{R}^2,\mathbb{R}^2)$ there exists a $1$-parameter group of proper $C^\infty$ diffeomorphisms $\varphi:\mathbb{R}\times\mathbb{R}^2\rightarrow\mathbb{R}^2$ as in \cite{Lee2009} Lemma 2.99 that satisfy
\begin{equation}\label{diffeomorphism_group}
\left.
\begin{array}{rl}
\partial_t\varphi(t,x)&= X(\varphi(t,x))\text{ for each }(t,x)\in\mathbb{R}\times\mathbb{R}^2;\\
\varphi(0,x)            &=  x  \text{ for each } x\in\mathbb{R}^2.\\
\end{array}
\right.
\end{equation}
We often use $\varphi_t$ to refer to the diffeomorphism $\varphi(t,\cdot):\mathbb{R}^2\rightarrow\mathbb{R}^2$.

\smallskip

\begin{lemma}\label{Taylor_expansion_of_flow}
Let $X\in C^\infty_c(\mathbb{R}^2,\mathbb{R}^2)$ and $\varphi$ be the corresponding flow as above. Then 
\begin{itemize}
\item[(i)] there exists $R\in C^\infty(\mathbb{R}\times\mathbb{R}^2,\mathbb{R}^2)$ and $K>0$ such that
\[
\varphi(t,x)=
\left\{
\begin{array}{ll}
x+tX(x)+R(t,x) & \text{ for }x\in\mathrm{supp}[X];\\
x & \text{ for }x\not\in\mathrm{supp}[X];\\
\end{array}
\right.
\]
where $|R(t,x)|\leq K t^2$  for $(t,x)\in\mathbb{R}\times\mathbb{R}^2$;
\item[(ii)] there exists $R^{(1)}\in C^\infty(\mathbb{R}\times\mathbb{R}^2,M_2(\mathbb{R}))$  and $K_1>0$ such that
\[
d\varphi(t,x)=
\left\{
\begin{array}{ll}
I+tdX(x)+R^{(1)}(t,x) & \text{ for }x\in\mathrm{supp}[X];\\
I & \text{ for }x\not\in\mathrm{supp}[X];\\
\end{array}
\right.
\]
where $|R^{(1)}(t,x)|\leq K_1 t^2$  for $(t,x)\in\mathbb{R}\times\mathbb{R}^2$;
\item[(iii)] there exists $R^{(2)}\in C^\infty(\mathbb{R}\times\mathbb{R}^2,\mathbb{R})$ and $K_2>0$ such that
\[
J_2d\varphi(t,x)=
\left\{
\begin{array}{ll}
1+t\,\mathrm{div}\,X(x)+R^{(2)}(t,x) & \text{ for }x\in\mathrm{supp}[X];\\
1 & \text{ for }x\not\in\mathrm{supp}[X];\\
\end{array}
\right.
\]
where $|R^{(2)}(t,x)|\leq K_2 t^2$  for $(t,x)\in\mathbb{R}\times\mathbb{R}^2$.
\end{itemize}
Let $x\in\mathbb{R}^2$, $v$ a unit vector in $\mathbb{R}^2$ and $M$ the line though $x$ perpendicular to $v$. Then
\begin{itemize}
\item[(iv)] there exists $R^{(3)}\in C^\infty(\mathbb{R}\times\mathbb{R}^2,\mathbb{R})$ and $K_3>0$ such that
\[
J_1d^{M}\varphi(t,x)=
\left\{
\begin{array}{ll}
1+t(\mathrm{div}^{M}\,X)(x)+R^{(3)}(t,x) & \text{ for }x\in\mathrm{supp}[X];\\
1 & \text{ for }x\not\in\mathrm{supp}[X];\\
\end{array}
\right.
\]
where $|R^{(3)}(t,x)|\leq K_3 t^2$  for $(t,x)\in\mathbb{R}\times\mathbb{R}^2$.
\end{itemize}
\end{lemma}

\smallskip

\noindent{\em Proof.}
{\em (i)} First notice that $\varphi\in C^\infty(\mathbb{R}\times\mathbb{R}^2)$ by \cite{Hale1969} Theorem 3.3 and Exercise 3.4. The statement for $x\not\in\mathrm{supp}[X]$ follows by uniqueness (cf. \cite{Hale1969} Theorem 3.1); the assertion for $x\in\mathrm{supp}[X]$ follows from Taylor's theorem. {\em (ii)} follows likewise: note, for example, that
\[
[\partial_{tt}d\varphi]_{\alpha\beta}\vert_{t=0}=X^\alpha_{,\beta\delta}X^\delta+X^\alpha_{,\gamma}X^\gamma_{,\beta}
\]
where the subscript $_,$ signifies partial differentiation. {\em (iii)} follows from {\em (ii)} and the definition of the $2$-dimensional Jacobian (cf. \cite{Ambrosio2000} Definition 2.68). {\em (iv)} Using \cite{Ambrosio2000} Definition 2.68 together with the Cauchy-Binet formula \cite{Ambrosio2000} Proposition 2.69, $J_1d^M\varphi(t,x)=|d\varphi(t,x)v|$ for $t\in\mathbb{R}$ and the result follows from {\em (ii)}.
\qed

\smallskip

\noindent Let $I$ be an open interval in $\mathbb{R}$ containing $0$. Let $Z:I\times\mathbb{R}^2\rightarrow\mathbb{R}^2;(t,x)\mapsto Z(t,x)$
be a continuous time-dependent vector field on $\mathbb{R}^2$ with the properties
\begin{itemize}
\item[(Z.1)] $Z(t,\cdot)\in C^1_c(\mathbb{R}^2,\mathbb{R}^2)$ for each $t\in I$;
\item[(Z.2)] $\mathrm{supp}[Z(t,\cdot)]\subset K$ for each $t\in I$ for some  compact set $K\subset\mathbb{R}^2$.
\end{itemize} 
By \cite{Hale1969} Theorems I.1.1, I.2.1, I.3.1, I.3.3 there exists a unique flow $\varphi:I\times\mathbb{R}^2\rightarrow\mathbb{R}^2$ such that
\begin{itemize}
\item[(F.1)] $\varphi:I\times\mathbb{R}^2\rightarrow\mathbb{R}^2$ is of class $C^1$;
\item[(F.2)] $\varphi(0,x)=x$ for each $x\in\mathbb{R}^2$;
\item[(F.3)] $\partial_t\varphi(t,x)=Z(t,\varphi(x,t))$ for each $(t,x)\in I\times\mathbb{R}^2$;
\item[(F.4)] $\varphi_t:=\varphi(t,\cdot):\mathbb{R}^2\rightarrow\mathbb{R}^2$ is a proper diffeomorphism for each $t\in I$.
\end{itemize}

\smallskip

\begin{lemma}\label{expansion_of_time_dependent_flow}
Let $Z$ be a time-dependent vector field with the properties $(Z.1)$-$(Z.2)$ and $\varphi$ be the corresponding flow. Then
\begin{itemize}
\item[(i)] for $(t,x)\in I\times\mathbb{R}^2$, 
\[
d\varphi(t,x)=
\left\{
\begin{array}{ll}
I+tdZ_0(x)+tR(t,x) & \text{ for }x\in K;\\
I & \text{ for }x\not\in K;\\
\end{array}
\right.
\]
where $\sup_K|R(t,\cdot)|\rightarrow 0$ as $t\rightarrow 0$.
\end{itemize}
Let $x\in\mathbb{R}^2$, $v$ a unit vector in $\mathbb{R}^2$ and $M$ the line though $x$ perpendicular to $v$. Then 
\begin{itemize}
\item[(ii)] for $(t,x)\in I\times\mathbb{R}^2$, 
\[
J_1d^M\varphi(t,x)=
\left\{
\begin{array}{ll}
1+t(\mathrm{div}^{M}\,Z_0)(x)+tR^{(1)}(t,x) & \text{ for }x\in K;\\
1 & \text{ for }x\not\in K.\\
\end{array}
\right.
\]
where $\sup_K|R^{(1)}(t,\cdot)|\rightarrow 0$ as $t\rightarrow 0$.
\end{itemize}
\end{lemma}

\smallskip

\noindent{\em Proof.} {\em (i)} We first remark that the flow $\varphi:I\times\mathbb{R}^2\rightarrow\mathbb{R}^2$ associated to $Z$ is continuously differentiable in $t,x$ in virtue of (Z.1) by \cite{Hale1969} Theorem I.3.3. Put $y(t,x):=d\varphi(t,x)$ for $(t,x)\in I\times\mathbb{R}^2$. By \cite{Hale1969} Theorem I.3.3,
\[
\dot{y}(t,x)=dZ(t,\varphi(t,x))y(t,x)
\]
for each $(t,x)\in I\times\mathbb{R}^2$ and $y(0,x)=I$ for each $x\in\mathbb{R}^2$ where $I$ stands for the $2\times 2$-identity matrix. For $x\in\ K$ and $t\in I$,
\begin{align}
d\varphi(t,x)&=I+d\varphi(t,x)-d\varphi(0,x)\nonumber\\
&=I+t\dot{y}(0,x)+t\Big\{\frac{d\varphi(t,x)-d\varphi(0,x)}{t}-\dot{y}(0,x)\Big\}\nonumber\\
&=I+t dZ(0,x)+t\Big\{\frac{y(t,x)-y(0,x)}{t}-\dot{y}(0,x)\Big\}\nonumber\\
&=I+t dZ_0(x)+t\Big\{\frac{y(t,x)-y(0,x)}{t}-\dot{y}(0,x)\Big\}.\nonumber
\end{align}
Applying the mean-value theorem component-wise and using uniform continuity of the matrix $\dot{y}$ in its arguments we see that 
\[
\frac{y(t,\cdot)-y(0,\cdot)}{t}-\dot{y}(0,\cdot)\rightarrow 0
\]
uniformly on $K$ as $t\rightarrow 0$. This leads to {\em (i)}. Part {\em (ii)} follows as in Lemma \ref{Taylor_expansion_of_flow}.
\qed

\smallskip

\noindent Let $E$ be a set of finite perimeter in $\mathbb{R}^2$ with $V_f(E)<+\infty$. The first variation of weighted volume resp. perimeter along $X\in C^\infty_c(\mathbb{R}^2,\mathbb{R}^2)$ is defined by
\begin{align}
\delta V_f(X)&:=\frac{d}{dt}\Big\vert_{t=0}V_f(\varphi_t(E)),\label{first_variation_of_volume}\\
\delta P_f^+(X)&:=\lim_{t\downarrow 0}\frac{P_f(\varphi_t(E))-P_f(E)}{t},\label{first_variation_of_perimeter}
\end{align}
whenever the limit exists. By Lemma \ref{diffeomorhism_invariance_of_boundary} the $f$-perimeter in (\ref{first_variation_of_perimeter}) is well-defined.

\smallskip

\noindent{\em Convex functions.} Suppose that $h:[0,+\infty)\rightarrow\mathbb{R}$ is a convex function. For $x\geq 0$ and $v\geq 0$ define
\[
h^\prime_+(x,v):=\lim_{t\downarrow 0}\frac{h(x+tv)-h(x)}{t}\in\mathbb{R}
\]
and define $h^\prime_-(x,v)$ similarly for $x>0$ and $v\leq 0$. For future use we introduce the notation
\[
\varrho_+:=h^\prime(\cdot,+1),
\varrho_-:=-h^\prime(\cdot,-1)
\text{ and }
\varrho:=(1/2)(\varrho_++\varrho_-)
\]
on $(0,+\infty)$. It holds that $h$ is differentiable a.e. and $h^\prime=\varrho$ a.e. on $(0,+\infty)$. Define $[\varrho]:=\varrho_+-\varrho_-$. Then $[\varrho]\geq 0$ and vanishes a.e.  on $(0,+\infty)$.

\smallskip

\begin{lemma}\label{differentiability_of_density}
Suppose that the function $f$ takes the form (\ref{form_of_f}) where $h:[0,+\infty)\rightarrow\mathbb{R}$ is a convex function. Then
\begin{itemize}
\item[(i)] the directional derivative $f^\prime_+(x,v)$ exists in $\mathbb{R}$ for each $x\in\mathbb{R}^2$ and $v\in\mathbb{R}^2$;
\item[(ii)] for $v\in\mathbb{R}^2$,
\[
f^\prime_+(x,v)=
\left\{
\begin{array}{ll}
f(x)h^\prime_+(|x|,\mathrm{sgn}\langle x,v\rangle)\frac{|\langle x,v\rangle|}{|x|} & \text{ for }x\in\mathbb{R}^2\setminus\{0\};\\
f(0)h^\prime_+(0,+1)|v| & \text{ for }x=0;\\
\end{array}
\right.
\]
\item[(iii)] if $M$ is a $C^1$ hypersurface in $\mathbb{R}^2$ such that $\cos\sigma\neq 0$ on $M$ then $f$ is differentiable $\mathscr{H}^1$-a.e. on $M$ and
\[
(\nabla f)(x)=f(x)\varrho(|x|)\frac{\langle x,\cdot\rangle}{|x|}
\]
for $\mathscr{H}^1$-a.e. $x\in M$.
\end{itemize}
\end{lemma}

\smallskip

\noindent{\em Proof.} The assertion in {\em (i)} follows from the monotonicity of chords property while {\em (ii)} is straightforward. {\em (iii)} Let $x\in M$ and $\gamma_1:I\rightarrow M$ be a $C^1$-parametrisation of $M$ near $x$ as above. Now $r_1\in C^1(I)$ and $\dot{r_1}(0)=\cos\sigma(x)\neq 0$ so we may assume that $r_1:I\rightarrow r_1(I)\subset(0,+\infty)$ is a $C^1$ diffeomorphism. The differentiability set $D(h)$ of $h$ has full Lebesgue measure in $[0,+\infty)$. It follows  by \cite{Ambrosio2000} Proposition 2.49 that $r_1^{-1}(D(h))$ has full measure in $I$. This entails that $f$ is differentiable $\mathscr{H}^1$-a.e. on $\gamma_1(I)\subset M$.
\qed

\section{Existence and $C^1$ regularity}

\smallskip

\noindent We start with an existence theorem. 

\begin{theorem}\label{first_properties_of_isoperimetric_set}
Assume that $f$ is a positive radial lower-semicontinuous non-decreasing density on $\mathbb{R}^2$ which diverges to infinity. Then for each $v>0$,
\begin{itemize}
\item[(i)] (\ref{isoperimetric_problem}) admits a minimiser;
\item[(ii)] any minimiser of (\ref{isoperimetric_problem}) is essentially bounded.
\end{itemize}
\end{theorem}

\smallskip

\noindent{\em Proof.}
See \cite{MorganPratelli2011} Theorems 3.3 and 5.9.
\qed

\smallskip

\noindent But the bulk of this section will be devoted to a discussion of $C^1$ regularity.

\smallskip

\begin{proposition}\label{variation_of_volume}
Let $f$ be a positive locally Lipschitz density on $\mathbb{R}^2$. Let $E\subset\mathbb{R}^2$ be a bounded set with finite perimeter. Let $X\in C^\infty_c(\mathbb{R}^2,\mathbb{R}^2)$. Then
\[
\delta V_f(X)
=\int_{E}\mathrm{div}(fX)\,dx
=-\int_{\mathscr{F}E} f\,\langle\nu^E, X\rangle\,d\mathscr{H}^1.
\]
\end{proposition}

\smallskip

\noindent{\em Proof.}
Let $t\in\mathbb{R}$. By the area formula (\cite{Ambrosio2000} Theorem 2.71 and (2.74)),
\begin{equation}\label{V_f_of_Phi_of_E}
V_f(\varphi_t(E)) = \int_{\varphi_t(E)}f\,dx
= \int_{E}(f\circ\varphi_t)\,J_2 d(\varphi_t)_x\,dx
\end{equation}
and
\begin{align}
V_f(\varphi_t(E))-V_f(E)&=\int_E(f\circ\varphi_t)J_2 d\varphi_t - f\,dx\nonumber\\
&=\int_E(f\circ\varphi_t)(J_2 d\varphi_t-1)\,dx+\int_Ef\circ\varphi_t-f\,dx.\nonumber
\end{align}
The density $f$ is locally Lipschitz and in particular differentiable a.e. on $\mathbb{R}^2$ (see \cite{Ambrosio2000} 2.3 for example). By the dominated convergence theorem and Lemma \ref{Taylor_expansion_of_flow},
\begin{align}
\delta V_f(X)
&=\int_{E}\Big\{f\mathrm{div}(X)+\langle\nabla f,X\rangle\Big\}\,dx
=\int_{E}\mathrm{div}(fX)\,dx
=-\int_{\mathscr{F}E} f\,\langle\nu^E,X\rangle\,d\mathscr{H}^1\nonumber
\end{align}
by the generalised Gauss-Green formula \cite{Ambrosio2000} Theorem 3.36.
\qed

\smallskip

\begin{proposition}\label{variation_of_perimeter}
Let $f$ be a positive locally Lipschitz density on $\mathbb{R}^2$.
Let $E\subset\mathbb{R}^2$ be a bounded set with finite perimeter. Let $X\in C^\infty_c(\mathbb{R}^2,\mathbb{R}^2)$. Then there exist constants $C>0$ and $\delta>0$ such that
\[
|P_f(\varphi_t(E))-P_f(E)|\leq C|t|
\]
for $|t|<\delta$.
\end{proposition}

\smallskip

\noindent{\em Proof.}
Let $t\in\mathbb{R}$. By Lemma \ref{diffeomorhism_invariance_of_boundary} and \cite{Ambrosio2000} Theorem 3.59,
\[
P_f(\varphi_t(E))
=\int_{\mathbb{R}^2}f\,d|D\chi_{\varphi_t(E)}|
=\int_{\mathscr{F}\varphi_t(E)}f\,d\mathscr{H}^1
=\int_{\varphi_t(\mathscr{F}E)}f\,d\mathscr{H}^1.
\]
As $\mathscr{F}E$ is countably $1$-rectifiable (\cite{Ambrosio2000} Theorem 3.59) we may use the generalised area formula \cite{Ambrosio2000} Theorem 2.91 to write
\[
P_f(\varphi_t(E))=\int_{\mathscr{F}E}(f\circ\varphi_t)J_1d^{\mathscr{F}E}(\varphi_t)_x\,d\mathscr{H}^1.
\]
For each $x\in\mathscr{F}E$ and any $t\in\mathbb{R}$,
\[
|(f\circ\varphi_t)(x)-f(x)|\leq K|\varphi(t,x)-x|\leq K\|X\|_\infty|t|
\]
where $K$ is the Lipschitz constant of $f$ on $\mathrm{supp}[X]$. The result follows upon writing
\begin{align}
P_f(\varphi_t(E))-P_f(E)&=\int_{\mathscr{F}E}(f\circ\varphi_t)(J_1d^{\mathscr{F}E}(\varphi_t)_x-1)+[f\circ\varphi_t-f]
\,d\mathscr{H}^1\label{difference_of_perimeters}
\end{align}
and using Lemma \ref{Taylor_expansion_of_flow}.
\qed

\smallskip

\begin{lemma}\label{vector_field_lemma}
Let $f$ be a positive locally Lipschitz density on $\mathbb{R}^2$. Let $E\subset\mathbb{R}^2$ be a bounded set with finite perimeter and $p\in\mathscr{F}E$. For any $r>0$ there exists $X\in C^\infty_c(\mathbb{R}^2,\mathbb{R}^2)$ with $\mathrm{supp}[X]\subset B(p,r)$ such that $\delta V_f(X)=1$.
\end{lemma}

\smallskip

\noindent{\em Proof.}
By (\ref{weighted_perimeter_of_E_relative_to_U}) and \cite{Ambrosio2000} Theorem 3.59 and (3.57) in particular,
\[
P_f(E,B(p,r))=\int_{B(p,r)\cap\mathscr{F}E}f\,d\mathscr{H}^1>0
\]
for any $r>0$. By the variational characterisation of the $f$-perimeter relative to $B(p,r)$ we can find $Y\in C^\infty_c(\mathbb{R}^2,\mathbb{R}^2)$ with $\mathrm{supp}[Y]\subset B(p,r)$ such that
\[
0<\int_{E\cap B(p,r)}\mathrm{div}(fY)\,dx=-\int_{\mathscr{F}E\cap B(p,r)}f\langle\nu^E,Y\rangle\,d\mathscr{H}^1=:c
\]
where we make use of the generalised Gauss-Green formula (cf. \cite{Ambrosio2000} Theorem 3.36).
Put $X:=(1/c)Y$. Then $X\in C^\infty_c(\mathbb{R}^2,\mathbb{R}^2)$ with $\mathrm{supp}[X]\subset B(p,r)$  and $\delta V_f(X)=1$ according to Proposition \ref{variation_of_volume}. 
\qed

\smallskip

\begin{proposition}\label{perimeter_of_set_combinations}
Let $f$ be a positive lower semi-continuous density on $\mathbb{R}^2$. Let $U$ be a bounded open set in $\mathbb{R}^2$ with Lipschitz boundary. Let $E, F_1, F_2$ be bounded sets in $\mathbb{R}^2$ with finite perimeter. Assume that $E\Delta F_1\subset\subset U$ and $E\Delta F_2\subset\subset\mathbb{R}^2\setminus\overline{U}$. Define 
\[
F:=\Big[F_1\cap U\Big]\cup\Big[F_2\setminus U\Big].
\]
Then $F$ is a set of finite perimeter in $\mathbb{R}^2$ and 
\[
P_f(E)+P_f(F)= P_f(F_1)+P_f(F_2).
\]
\end{proposition}

\smallskip

\noindent{\em Proof.}
The function $\chi_E\vert_U\in\mathrm{BV}(U)$ and $D(\chi_E\vert_U)=(D\chi_E)\vert_U$. We write $\chi_E^U$ for the boundary trace of $\chi_E\vert_U$ (see \cite{Ambrosio2000} Theorem 3.87); then $\chi_E^U\in L^1(\partial U,\mathscr{H}^1\mres\partial U)$ (cf. \cite{Ambrosio2000} Theorem 3.88). We use similar notation elsewhere. By \cite{Ambrosio2000} Corollary 3.89,
\begin{align}
D\chi_E&=D\chi_E\mres U+(\chi_E^U-\chi_E^{\mathbb{R}^2\setminus\overline{U}})\nu^U\mathscr{H}^1\mres{\partial U}+D\chi_E\mres(\mathbb{R}^2\setminus\overline{U});\nonumber\\
D\chi_F&=D\chi_{F_1}\mres U+(\chi_{F_1}^U-\chi_{F_2}^{\mathbb{R}^2\setminus\overline{U}})\nu^U\mathscr{H}^1\mres{\partial U}+D\chi_{F_2}\mres(\mathbb{R}^2\setminus\overline{U});\nonumber\\
D\chi_{F_1}&=D\chi_{F_1}\mres U+(\chi_{F_1}^U-\chi_{E}^{\mathbb{R}^2\setminus\overline{U}})\nu^U\mathscr{H}^1\mres{\partial U}+D\chi_{E}\mres(\mathbb{R}^2\setminus\overline{U});\nonumber\\
D\chi_{F_2}&=D\chi_{E}\mres U+(\chi_{E}^U-\chi_{F_2}^{\mathbb{R}^2\setminus\overline{U}})\nu^U\mathscr{H}^1\mres{\partial U}+D\chi_{F_2}\mres(\mathbb{R}^2\setminus\overline{U}).\nonumber
\end{align}
From the definition of the total variation measure (\cite{Ambrosio2000} Definition 1.4),
\begin{align}
|D\chi_E|&=|D\chi_E|\mres U+|\chi_E^U-\chi_E^{\mathbb{R}^2\setminus\overline{U}}|\mathscr{H}^1\mres{\partial U}+|D\chi_E|\mres(\mathbb{R}^2\setminus\overline{U});\nonumber\\
|D\chi_F|&=|D\chi_{F_1}|\mres U+|\chi_E^U-\chi_E^{\mathbb{R}^2\setminus\overline{U}}|\mathscr{H}^1\mres{\partial U}+|D\chi_{F_2}|\mres(\mathbb{R}^2\setminus\overline{U});\nonumber\\
|D\chi_{F_1}|&=|D\chi_{F_1}|\mres U+|\chi_E^U-\chi_{E}^{\mathbb{R}^2\setminus\overline{U}}|\mathscr{H}^1\mres{\partial U}+|D\chi_{E}|\mres(\mathbb{R}^2\setminus\overline{U});\nonumber\\
|D\chi_{F_2}|&=|D\chi_{E}|\mres U+|\chi_{E}^U-\chi_E^{\mathbb{R}^2\setminus\overline{U}}|\mathscr{H}^1\mres{\partial U}+|D\chi_{F_2}|\mres(\mathbb{R}^2\setminus\overline{U});\nonumber
\end{align}
where we also use the fact that $\chi_{F_1}^U=\chi_E^U$ as $E\Delta F_1\subset\subset U$ and similarly for $F_2$. The result now follows.
\qed

\smallskip

\begin{proposition}\label{version_of_Morgan}
Assume that $f$ is a positive locally Lipschitz density on $\mathbb{R}^2$. Let $v>0$ and suppose that the set $E$ is a bounded minimiser of (\ref{isoperimetric_problem}). Let $U$ be a bounded open set in $\mathbb{R}^2$. There exist constants $C>0$ and $\delta>0$ with the following property. For any $x\in U$ and $0<r<\delta$,
\begin{equation}\label{Morgan_lemma_assertion_1}
P_f(E)-P_f(F)\leq C\big|V_f(E)-V_f(F)\big|
\end{equation}
where $F$ is any set with finite perimeter in $\mathbb{R}^2$ such that $E\Delta F\subset\subset B(x,r)$.\end{proposition}

\smallskip

\noindent{\em Proof.}
The proof follows that of \cite{Morgan2003} Proposition 3.1. We assume to the contrary that
 \[
(\forall\,C>0)(\forall\,\delta>0)(\exists\,x\in  U)(\exists\,r\in(0,\delta))(\exists\,F\subset\mathbb{R}^2)
\]
\begin{equation}\label{hypothesis_for_contradiction}
\Big[
F\Delta E\subset\subset B(x,r)
\wedge
\Delta P_f > C|\Delta V_f|
\Big]
\end{equation}
in the language of quantifiers where we have taken some liberties with notation.

\smallskip

\noindent Choose $p_1,p_2\in\mathscr{F}E$ with $p_1\neq p_2$. Choose $r_0>0$ such that the open balls $B(p_1,r_0)$ and $B(p_2,r_0)$ are disjoint. Choose vector fields 
$X_j\in C^\infty_c(\mathbb{R}^2,\mathbb{R}^2)$ with $\mathrm{supp}[X_j]\subset B(p_j,r_0)$ such that
\begin{equation}\label{choice_of_vector_field}
\delta V_f(X_j)=1
\text{ and }
|P_f(\varphi^{(j)}_t(E))-P_f(E)|\leq a_j|t|\text{ for }|t|<\delta_j\text{ and }j=1,2
\end{equation}
as in Lemma \ref{vector_field_lemma} and Proposition \ref{variation_of_perimeter}. Put $a:=\max\{a_1,a_2\}$. By (\ref{choice_of_vector_field}),
\[
V_f(\varphi^{(j)}_t(E))-V_f(E)=t + o(t)\text{ as }t\rightarrow 0\text{ for }j=1,2.
\]
So there exist $\varepsilon>0$ and $1>\eta>0$ such that
\begin{align}
&t-\eta|t| <  V_f(\varphi^{(j)}_t(E)) - V_f(E)<  t + \eta|t|;\label{f_volume_estimate}\\
&|P_f(\varphi_t^{(j)}(E)) - P_f(E) |<(a+1)|t|;\nonumber
\end{align}
for $|t|<\varepsilon$ and $j=1,2$. In particular,
\begin{align}
&|V_f(\varphi^{(j)}_t(E)) - V_f(E)|>(1-\eta)|t|;\nonumber\\
&|P_f(\varphi_t^{(j)}(E)) - P_f(E)| <\frac{1+a}{1-\eta}|V_f(\varphi^{(j)}_t(E)) - V_f(E)|
\text{ for }|t|<\varepsilon;\label{perimeter_of_deformed_sets}
\end{align}
for $|t|<\varepsilon$ and $j=1,2$.

\smallskip

\noindent In (\ref{hypothesis_for_contradiction}) choose $C=(1+a)/(1-\eta)$ and $\delta>0$ such that
\begin{itemize}
\item[(a)] $0<2\delta<\mathrm{dist}(B(p_1,r_0),B(p_2,r_0))$, 
\item[(b)] $\sup\{V_f(B(x,\delta)):\,x\in U\}<(1-\eta)\,\varepsilon$.
\end{itemize}
Choose $x,r$ and $F_1$ as in (\ref{hypothesis_for_contradiction}). In light of (a) we may assume that $B(x,r)\cap B(p_1,r_0)=\emptyset$. By (b),
\begin{equation}\label{f_volume_for_F_1}
|V_f(F_1) - V_f(E)|\leq V_f(B(x,r))\leq V_f(B(x,\delta))<(1-\eta)\,\varepsilon.
\end{equation}
From (\ref{f_volume_estimate}) and (\ref{f_volume_for_F_1}) we can find $t\in(-\varepsilon,\varepsilon)$ such that with $F_2:=\varphi_t^{(1)}(E)$,
\begin{equation}\label{equality_of_volumes}
V_f(F_2) - V_f(E)=-\Big\{V_f(F_1) - V_f(E)\Big\}
\end{equation}
by the intermediate value theorem. From (\ref{hypothesis_for_contradiction}),
\begin{equation}\label{Morgan_lemma_inequality_1}
P_f(F_1) < P_f(E)-C|V_f(F_1) - V_f(E)|
\end{equation}
while from (\ref{perimeter_of_deformed_sets}),
\begin{equation}\label{Morgan_lemma_inequality_2}
P_f(F_2) <  P_f(E) + C|V_f(F_2) - V_f(E)|.
\end{equation}
Let $F$ be the set
\[
F:=\Big[F_1\setminus B(p_1,r_0))\Big]\cup\Big[B(p_1,r_0)\cap F_2\Big].
\]
Note that $E\Delta F_2\subset\subset B(p_1,r_0)$. By Proposition \ref{perimeter_of_set_combinations}, $F$ is a bounded set of finite perimeter in $\mathbb{R}^2$ and
\[
P_f(E)+P_f(F)= P_f(F_1)+P_f(F_2).
\]
We then infer from (\ref{Morgan_lemma_inequality_1}), (\ref{Morgan_lemma_inequality_2}) and (\ref{equality_of_volumes}) that
\begin{align}
P_f(F)&= P_f(F_1)+P_f(F_2)-P_f(E)\nonumber\\
&<P_f(E)-C|V_f(F_1) - V_f(E)|+ P_f(E) + C|V_f(F_2) - V_f(E)|-P_f(E)\
=P_f(E).\nonumber
\end{align}
On the other hand, $V_f(F)=V_f(F_1)+V_f(F_2)-V_f(E)=V_f(E)$ by (\ref{equality_of_volumes}). 
We therefore obtain a contradiction to the $f$-isoperimetric property of $E$.
\qed

\smallskip

\noindent Let $E$ be a set of finite perimeter in $\mathbb{R}^2$ and $U$ a bounded open set in $\mathbb{R}^2$. The minimality excess is the function $\psi$ defined by
\begin{align}
\psi(E,U)&:=P(E,U)-\nu(E,U)
\end{align}
where
\begin{align}
\nu(E,U) &:=\inf\{P(F,U):F\text{ is a set of finite perimeter with }F\Delta E\subset\subset U\}\nonumber
\end{align}
as in \cite{Tamanini1984} (1.9). We recall that the boundary of $E$ is said to be almost minimal in $\mathbb{R}^2$ if for each bounded open set $U$ in $\mathbb{R}^2$ there exists $T>0$ and a positive constant $K$ such that for every $x\in U$ and $r\in(0,T)$,
\begin{equation}\label{almost_minimal_inequality}
\psi(E,B(x,r))\leq K r^2.
\end{equation}
This definition corresponds to \cite{Tamanini1984} Definition 1.5.

\smallskip

\begin{theorem}\label{boundary_is_almost_minimal}
Assume that $f$ is a positive locally Lipschitz density on $\mathbb{R}^2$. Let $v>0$ and assume that $E$ is a bounded minimiser of (\ref{isoperimetric_problem}). Then the boundary of $E$ is almost minimal in $\mathbb{R}^2$.
\end{theorem}

\smallskip

\noindent{\em Proof.} 
Let $U$ be a bounded open set in $\mathbb{R}^2$ and $C>0$ and $\delta>0$ as in Proposition \ref{version_of_Morgan}. The open $\delta$-neighbourhood of $U$ is denoted $I_\delta(U)$. Let $x\in U$ and $r\in(0,\delta)$. Put $V:=I_{2\delta}(U)$. For the sake of brevity write $m:=\inf_{B(x,r)}f$ and $M:=\sup_{B(x,r)}f$. Let $F$ be a set of finite perimeter in $\mathbb{R}^2$ such that $F\Delta E\subset\subset B(x,r)$. By Proposition \ref{version_of_Morgan},
\begin{align}
&P(E,B(x,r))-P(F,B(x,r))\nonumber\\
&\leq
\frac{1}{m}P_f(E,B(x,r))-\frac{1}{M}P_f(F,B(x,r))\nonumber\\
&=\frac{1}{m}
\Big(
P_f(E,B(x,r))-P_f(F,B(x,r))
\Big)+
\Big(\frac{1}{m}-\frac{1}{M}
\Big)
P_f(F,B(x,r))\nonumber\\
&\leq\frac{1}{m}
\Big(
P_f(E,B(x,r))-P_f(F,B(x,r))
\Big)+
\frac{M-m}{m^2}
P_f(F,B(x,r))\nonumber\\
&\leq\frac{C}{\inf_V f}|V_f(E)-V_f(F)|
+
(2Lr)\frac{\sup_V f}{(\inf_V f)^2}P(F,B(x,r))\nonumber\\
&\leq C\pi r^2\frac{\sup_V f}{\inf_V f}
+
(2Lr)\frac{\sup_V f}{(\inf_V f)^2}P(F,B(x,r))\nonumber
\end{align}
where $L$ stands for the Lipschitz constant of the restriction of $f$ to $V$. We then derive that
\[
\psi(E,B(x,r))\leq C\pi r^2\frac{\sup_V f}{\inf_V f}
+
(2Lr)\frac{\sup_V f}{(\inf_V f)^2}
\nu(E,B(x,r)).
\]
By \cite{Giusti1984} (5.14), $\nu(E,B(x,r))\leq\pi r$. The inequality in (\ref{almost_minimal_inequality}) now follows.
\qed

\smallskip

\begin{theorem}\label{C1_property_of_reduced_boundary}
Assume that $f$ is a positive locally Lipschitz density on $\mathbb{R}^2$. Let $v>0$ and suppose that $E$ is a bounded minimiser of (\ref{isoperimetric_problem}). Then there exists a set $\widetilde{E}\subset\mathbb{R}^2$ such that
\begin{itemize}
\item[(i)] $\widetilde{E}$ is a bounded minimiser of (\ref{isoperimetric_problem});
\item[(ii)] $\widetilde{E}$ is equivalent to $E$;
\item[(iii)] $\widetilde{E}$ is open and $\partial\widetilde{E}$ is a $C^1$ hypersurface in $\mathbb{R}^2$. 
\end{itemize}
\end{theorem}

\smallskip

\noindent{\em Proof.}
By \cite{Giusti1984} Proposition 3.1 there exists a Borel set $F$ equivalent to $E$ with the property that
\[
\partial F=\{x\in\mathbb{R}^2:0<|F\cap B(x,\rho)|<\pi\rho^2\text{ for each }\rho>0\}.
\] 
By Theorem \ref{boundary_is_almost_minimal} and \cite{Tamanini1984} Theorem 1.9, $\partial F$ is a $C^1$ hypersurface in $\mathbb{R}^2$ (taking note of differences in notation). The set
\[
\widetilde{E}:=\{x\in\mathbb{R}^2:|F\cap B(x,\rho)|=\pi\rho^2\text{ for some }\rho>0\}
\]
satisfies {\em (i)}-{\em (iii)}.
\qed

\section{Weakly bounded curvature and spherical cap symmetry}

\smallskip

\begin{theorem}\label{C11_property_of_boundary}
Assume that $f$ is a positive locally Lipschitz density on $\mathbb{R}^2$. Let $v>0$ and suppose that $E$ is a bounded minimiser of (\ref{isoperimetric_problem}). Then there exists a set $\widetilde{E}\subset\mathbb{R}^2$ such that
\begin{itemize}
\item[(i)] $\widetilde{E}$ is a bounded minimiser of (\ref{isoperimetric_problem});
\item[(ii)] $\widetilde{E}$ is equivalent to $E$;
\item[(iii)] $\widetilde{E}$ is open and $\partial\widetilde{E}$ is a $C^{1,1}$ hypersurface in $\mathbb{R}^2$. 
\end{itemize}
\end{theorem}

\smallskip

\noindent{\em Proof.}
We may assume that $E$ has the properties listed in Theorem \ref{C1_property_of_reduced_boundary}. Put $M:=\partial E$. Let $x\in M$ and $U$ a bounded open set containing $x$. Choose $C>0$ and $\delta>0$ as in Proposition \ref{version_of_Morgan}. Let $0<r<\delta$ and $X\in C^\infty_c(\mathbb{R}^2,\mathbb{R}^2)$ with $\mathrm{supp}[X]\subset B(x,r)$. Then
\[
P_f(E)-P_f(\varphi_t(E))\leq C|V_f(E)-V_f(\varphi_t(E))|
\]
for each $t\in\mathbb{R}$. From the identity (\ref{difference_of_perimeters}),
\begin{align}
-\int_M(f\circ\varphi_t)(J_1d^M(\varphi_t)_x-1)\,d\mathscr{H}^1&\leq C|V_f(E)-V_f(\varphi_t(E))|
+\int_M[f\circ\varphi_t-f]\,d\mathscr{H}^1\nonumber\\
&\leq C|V_f(E)-V_f(\varphi_t(E))|
+\sqrt{2}K\|X\|_\infty\mathscr{H}^1(M\cap\mathrm{supp}[X])t\nonumber
\end{align}
where $K$ stands for the Lipschitz constant of $f$ restricted to $U$. On dividing by $t$ and taking the limit $t\rightarrow 0$ we obtain
\begin{align}
-\int_M f\mathrm{div}^M X\,d\mathscr{H}^1&\leq C\Big|\int_M f\langle n,X\rangle\,d\mathscr{H}^1\Big|
+\sqrt{2}K\|X\|_\infty\mathscr{H}^1(M\cap\mathrm{supp}[X])\nonumber
\end{align}
upon using Lemma \ref{Taylor_expansion_of_flow} and Proposition \ref{variation_of_volume}. Replacing $X$ by $-X$ we derive that
\[
\Big|\int_M f\mathrm{div}^M X\,d\mathscr{H}^1\Big|\leq C_1\|X\|_\infty\mathscr{H}^1(M\cap\mathrm{supp}[X])
\]
where $C_1=C\|f\|_{L^\infty(U)}+\sqrt{2}K$. Let $\gamma_1:I\rightarrow M$ be a local $C^1$ parametrisation of $M$ near $x$. Suppose that $Y\in C^1_c(I,\mathbb{R}^2)$ with $\mathrm{supp}[Y]\subset I$ and that $\gamma_1(I)\subset M\cap B(x,r)$. Note that there exists $X\in C^\infty_c(\mathbb{R}^2,\mathbb{R}^2)$ with $\mathrm{supp}[X]\subset B(x,r)$ such that $X\circ\gamma_1=Y$ on $I$. The above estimate entails that
\[
\Big|\int_I(f\circ\gamma_1)\langle\dot{Y},t\rangle\,ds\Big|\leq 
C_1\Big|\mathrm{supp}[Y]\Big|\|Y\|_\infty.
\]
This means that the function $(f\circ\gamma_1)t$ belongs to $\mathrm{BV}(I)$ and this implies in turn that $t\in\mathrm{BV}(I)$. For $s_1,s_2\in I$ with $s_1<s_2$,
\begin{align}
|t(s_2)-t(s_1)|&=|Dt((s_1,s_2))|\leq |Dt|((s_1,s_2))\nonumber\\
&=\sup\Big\{\int_{(s_1,s_2)}\langle t,\dot{Y}\rangle\,ds:Y\in C^1_c((s_1,s_2))\text{ and }\|Y\|_\infty\leq 1\Big\}\nonumber\\
&\leq c\sup\Big\{\int_{(s_1,s_2)}(f\circ\gamma_1)\langle t,\dot{Y}\rangle\,ds:Y\in C^1_c((s_1,s_2))\text{ and }\|Y\|_\infty\leq 1\Big\}\nonumber\\
&\leq cC_1|s_2-s_1|\nonumber
\end{align}
where $1/c=\inf_{\overline{U}}f>0$. It follows that $M$ is of class $C^{1,1}$. \qed

\smallskip

\noindent We turn to the topic of spherical cap symmetrisation. Denote by $\mathbb{S}^1_\tau$ the centred circle in $\mathbb{R}^2$ with radius $\tau>0$. We sometimes write $\mathbb{S}^1$ for $\mathbb{S}^1_1$. Given $x\in\mathbb{R}^2$, $v\in\mathbb{S}^1$ and $\alpha\in(0,\pi]$ the open cone with vertex $x$, axis $v$ and opening angle $2\alpha$ is the set
\[
C(x,v,\alpha):=\Big\{y\in\mathbb{R}^2:\langle y-x,v\rangle>|y-x|\cos\alpha\Big\}.
\]

\smallskip

\noindent Let $E$ be an $\mathscr{L}^2$-measurable set in $\mathbb{R}^2$ and $\tau>0$. The $\tau$-section $E_\tau$ of $E$ is the set $E_\tau:=E\cap\mathbb{S}^1_\tau$.
Put
\begin{equation}\label{definition_of_L}
L(\tau)=L_E(\tau):=\mathscr{H}^1(E_\tau)\text{ for }\tau>0
\end{equation}
and $p(E) := \{\tau>0:L(\tau)>0\}$. The function $L$ is $\mathscr{L}^1$-measurable by \cite{Ambrosio2000} Theorem 2.93. Given $\tau>0$ and $0<\alpha\leq\pi$ the spherical cap $C(\tau,\alpha)$ is the set
\[
C(\tau,\alpha):=
\left\{
\begin{array}{ll}
\mathbb{S}^1_\tau\cap C(0,e_1,\alpha) & \text{ if }0<\alpha<\pi;\\
\mathbb{S}^1_\tau & \text{ if }\alpha=\pi;\\
\end{array}
\right.
\]
and has $\mathscr{H}^1$-measure $s(\tau,\alpha):=2\alpha\tau$. The spherical cap symmetral $E^{sc}$ of the set $E$ is defined by
\begin{equation}\label{definition_of_spherical_cap_symmetral}
E^{sc}:=\bigcup_{\tau\in p(E)}C(\tau,\alpha)
\end{equation}
where $\alpha\in(0,\pi]$ is determined by $s(\tau,\,\alpha)=L(\tau)$. Observe that $E^{sc}$ is a $\mathscr{L}^2$-measurable set in $\mathbb{R}^2$ and $V_f(E^{sc})=V_f(E)$. Note also that if $B$ is a centred open ball then $B^{sc}=B\setminus\{0\}$. We say that $E$ is spherical cap symmetric if $\mathscr{H}^1((E\Delta E^{sc})_\tau)=0$ for each $\tau>0$. This definition is broad but suits our purposes.

\smallskip

\noindent The result below is stated in \cite{MorganPratelli2011} Theorem 6.2 and a sketch proof given. A proof along the lines of \cite{Barchiesietal2013} Theorem 1.1 can be found in \cite{Perugini2016}. First, let $B$ be a Borel set in $(0,+\infty)$; then the annulus $A(B)$ over $B$ is the set $A(B):=\{x\in\mathbb{R}^2:|x|\in B\}$.

\smallskip

\begin{theorem}\label{spherical_cap_decreases_perimeter}
Let $E$ be a set of finite perimeter in $\mathbb{R}^2$. Then $E^{sc}$ is a set of finite perimeter and
\begin{equation}
P(E^{sc},A(B))\leq P(E,A(B))
\end{equation}
for any Borel set $B\subset(0,\infty)$ and the same inequality holds with $E^{sc}$ replaced by any set $F$ that is $\mathscr{L}^2$-equivalent to $E^{sc}$.
\end{theorem}

\smallskip

\begin{corollary}\label{f_perimeter_decreases_under_spherical_cap_symmetrisation}
Let $f$ be a positive lower semi-continuous radial function on $\mathbb{R}^2$. Let $E$ be a set of finite perimeter in $\mathbb{R}^2$. Then
$
P_f(E^{sc})\leq P_f(E)
$.
\end{corollary}

\smallskip

\noindent{\em Proof.}
Assume that $P_f(E)<+\infty$. We remark that $f$ is Borel measurable as $f$ is lower semi-continuous. Let $(f_h)$ be a sequence of simple Borel measurable radial functions on $\mathbb{R}^2$ such that $0\leq f_h\leq f$ and $f_h\uparrow f$ on $\mathbb{R}^2$ as $h\rightarrow\infty$. By Theorem \ref{spherical_cap_decreases_perimeter},
\[
P_{f_h}(E^{sc})=\int_{\mathbb{R}^2}f_h\,d|D\chi_{E^{sc}}|
\leq\int_{\mathbb{R}^2}f_h\,d|D\chi_{E}|= P_{f_h}(E)
\]
for each $h$. Taking the limit $h\rightarrow\infty$ the monotone convergence theorem gives $P_f(E^{sc})\leq P_f(E)$. 
\qed

\smallskip

\begin{lemma}\label{spherical_cap_symmetry_and_measure_theoretic_interior}
Let $E$ be an $\mathscr{L}^2$-measurable set in $\mathbb{R}^2$ such that $E\setminus\{0\}=E^{sc}$. Then there exists an $\mathscr{L}^2$-measurable set $F$ equivalent to $E$ such that
\begin{itemize}
\item[(i)] $\partial F=\{x\in\mathbb{R}^2:0<|F\cap B(x,\rho)|<|B(x,\rho)|\text{ for any }\rho>0\}$;
\item[(ii)] $F$ is spherical cap symmetric.
\end{itemize}
\end{lemma}

\smallskip

\noindent{\em Proof.}
Put 
\begin{align}
E_1&:=\{x\in\mathbb{R}^2:|E\cap B(x,\rho)|=|B(x,\rho)|\text{ for some }\rho>0\};\nonumber\\
E_0&:=\{x\in\mathbb{R}^2:|E\cap B(x,\rho)|=0\text{ for some }\rho>0\}.\nonumber
\end{align}
We claim that $E_1$ is spherical cap symmetric. For take $x\in E_1$ with $\tau=|x|>0$ and $|\theta(x)|\in(0,\pi]$. Now $|E\cap B(x,\rho)|=|B(x,\rho)|$ for some $\rho>0$. Let $y\in\mathbb{R}^2$ with $|y|=\tau$ and $|\theta(y)|<|\theta(x)|$. Choose a rotation $O\in\mathrm{SO}(2)$ such that $OB(x,\rho)=B(y,\rho)$. As $E\setminus\{0\}=E^{sc}$, $|E\cap B(y,\rho)|=|O(E\cap B(x,\rho))|=|E\cap B(x,\rho)|=|B(x,\rho)|=|B(y,\rho)|$. The claim follows. It follows in a similar way that $\mathbb{R}^2\setminus E_0$ is spherical cap symmetric. It can then be seen that the set $F:=(E_1\cup E)\setminus E_0$ inherits this property. As in \cite{Giusti1984} Proposition 3.1 the set $F$ is equivalent to $E$ and enjoys the property in {\em (i)}. 
\qed

\smallskip

\begin{theorem}\label{M_is_spherical_cap_symmetric}
Let $f$ be as in (\ref{form_of_f}) where $h:[0,+\infty)\rightarrow\mathbb{R}$ is a non-decreasing convex function. Given $v>0$ let $E$ be a bounded minimiser of (\ref{isoperimetric_problem}). Then there exists an $\mathscr{L}^2$-measurable set $\widetilde{E}$ with the properties
\begin{itemize}
\item[(i)] $\widetilde{E}$ is a minimiser of (\ref{isoperimetric_problem});
\item[(ii)] $L_{\widetilde{E}}=L$ a.e. on $(0,+\infty)$;
\item[(iii)] $\widetilde{E}$ is open, bounded and has $C^{1,1}$ boundary;
\item[(iv)] $\widetilde{E}\setminus\{0\}=\widetilde{E}^{sc}$.
\end{itemize}
\end{theorem}

\smallskip

\noindent{\em Proof.} Let $E$ be a bounded minimiser for (\ref{isoperimetric_problem}). Then $E_1:=E^{sc}$ is a bounded minimiser of (\ref{isoperimetric_problem}) by Corollary \ref{f_perimeter_decreases_under_spherical_cap_symmetrisation} and $L_E=L_{E_1}$ on $(0,+\infty)$. Now put $E_2:=F$ with $F$ as in Lemma \ref{spherical_cap_symmetry_and_measure_theoretic_interior}. Then $L_{E_2}=L$ a.e. on $(0,+\infty)$ as $E_2$ is equivalent to $E_1$, $E_2$ is a bounded minimiser of (\ref{isoperimetric_problem}) and $E_2$ is spherical cap symmetric. Moreover, $\partial E_2=\{x\in\mathbb{R}^2:0<|E_2\cap B(x,\rho)|<|B(x,\rho)|\text{ for any }\rho>0\}$. As in the proof of Theorem \ref{C1_property_of_reduced_boundary}, $\partial E_2$ is a $C^1$ hypersurface in $\mathbb{R}^2$. Put
\[
\widetilde{E}:=\{x\in\mathbb{R}^2:|E_2\cap B(x,\rho)|=|B(x,\rho)|\text{ for some }\rho>0\}.
\]
Then $\widetilde{E}$ is equivalent to $E_2$ so that {\em (ii)} holds, and is a bounded minimiser of (\ref{isoperimetric_problem}); $\widetilde{E}$ is open and $\partial\widetilde{E}=\partial E_2$ is $C^1$. In fact, $\partial\widetilde{E}$ is of class $C^{1,1}$ by Theorem \ref{C11_property_of_boundary}. As $E_2$ is spherical cap symmetric the same is true of $\widetilde{E}$. But $\widetilde{E}$ is open which entails that $\widetilde{E}\setminus\{0\}=\widetilde{E}^{sc}$. \qed

\section{More on spherical cap symmetry}

\smallskip

\noindent Let
\[
H:=\{x=(x_1,x_2)\in\mathbb{R}^2:x_2>0\}
\]
stand for the open upper half-plane in $\mathbb{R}^2$ and
\[
S:\mathbb{R}^2\rightarrow\mathbb{R}^2;x=(x_1,x_2)\mapsto(x_1,-x_2)
\]
for reflection in the $x_1$-axis. Let $O\in\mathrm{SO}(2)$ represent rotation anti-clockwise through $\pi/2$. 

\smallskip

\begin{lemma}\label{oddness_of_cos_sigma}
Let $E$ be an open set in $\mathbb{R}^2$ with $C^1$ boundary $M$ and assume that $E\setminus\{0\}=E^{sc}$. Let $x\in M\setminus\{0\}$. Then 
\begin{itemize}
\item[(i)] $Sx\in M\setminus\{0\}$;
\item[(ii)] $n(Sx)=Sn(x)$;
\item[(iii)] $\cos\sigma(Sx)=-\cos\sigma(x)$.
\end{itemize} 
\end{lemma}

\smallskip

\noindent{\em Proof.}
{\em (i)} The closure $\overline{E}$ of $E$ is spherical cap symmetric. The spherical cap symmetral $\overline{E}$ is invariant under $S$ from the representation (\ref{definition_of_spherical_cap_symmetral}). {\em (ii)} is a consequence of this last observation. {\em (iii)} Note that $t(Sx)=O^\star n(Sx)=O^\star Sn(x)$. Then
\begin{align}
\cos\sigma(Sx)&=\langle Sx,t(Sx)\rangle=\langle Sx,O^\star Sn(x)\rangle=\langle x,SO^\star Sn(x)\rangle\nonumber\\
&=\langle x,On(x)\rangle=-\langle x,O^\star n(x)\rangle=\cos\sigma(x)\nonumber
\end{align}
as $SO^\star S=O$ and $O=-O^\star$.
\qed

\smallskip

\noindent We introduce the projection $\pi:\mathbb{R}^2\rightarrow[0,+\infty);x\mapsto|x|$. 

\smallskip

\begin{lemma}\label{property_of_E_sc}
Let $E$ be an open set in $\mathbb{R}^2$ with boundary $M$ and assume that $E\setminus\{0\}=E^{sc}$. 
\begin{itemize}
\item[(i)] Suppose $0\neq x\in\mathbb{R}^2\setminus\overline{E}$ and $\theta(x)\in(0,\pi]$. Then there exists an open interval $I$ in $(0,+\infty)$ containing $\tau$ and $\alpha\in(0,\theta(x))$ such that
$A(I)\setminus\overline{S}(\alpha)\subset\mathbb{R}^2\setminus\overline{E}$.
\item[(ii)] Suppose $0\neq x\in E$ and $\theta(x)\in[0,\pi)$. Then there exists an open interval $I$ in $(0,+\infty)$ containing $\tau$ and $\alpha\in(\theta(x),\pi)$ such that $A(I)\cap S(\alpha)\subset E$.
\item[(iii)] For each $0<\tau\in \pi(M)$, $M_\tau$ is the union of two closed  spherical arcs in $\mathbb{S}^1_\tau$ symmetric about the $x_1$-axis.
\end{itemize}
\end{lemma}

\smallskip

\noindent{\em Proof.}
{\em (i)} We can find $\alpha\in(0,\theta(x))$ such that $\mathbb{S}^1_\tau\setminus S(\alpha)\subset\mathbb{R}^2\setminus\overline{E}$ as can be seen from definition (\ref{definition_of_spherical_cap_symmetral}). This latter set is compact so $\mathrm{dist}(\mathbb{S}^1_\tau\setminus S(\alpha),\overline{E})>0$. This means that the $\varepsilon$-neighbourhood of $\mathbb{S}^1_\tau\setminus S(\alpha)$ is contained in $\mathbb{R}^2\setminus\overline{E}$ for $\varepsilon>0$ small. The claim follows. {\em (ii)} Again from (\ref{definition_of_spherical_cap_symmetral}) we can find $\alpha\in(\theta(x),\pi)$ such that $\overline{\mathbb{S}^1_\tau\cap S(\alpha)}\subset E$ and the assertion follows as before. 

\smallskip

\noindent{\em (iii)} Suppose $x_1,x_2$ are distinct points in $M_\tau$ with $0\leq\theta(x_1)<\theta(x_2)\leq\pi$. Suppose $y$ lies in the interior of the spherical arc joining $x_1$ and $x_2$. If $y\in\mathbb{R}^2\setminus\overline{E}$ then $x_2\in\mathbb{R}^2\setminus\overline{E}$ by {\em (i)} and hence $x_2\not\in M$. If $y\in E$ we obtain the contradiction that $x_1\in E$ by {\em (ii)}. Therefore $y\in M$. We infer that the closed spherical arc joining $x_1$ and $x_2$ lies in $M_\tau$. The claim follows noting that $M_\tau$ is closed. 
\qed

\smallskip

\begin{lemma}\label{property_of_n}
Let $E$ be an open set in $\mathbb{R}^2$ with $C^1$ boundary $M$. Let $x\in M$. Then
\[
\liminf_{\overline{E}\ni y\rightarrow x}\Big\langle\frac{y-x}{|y-x|},n(x)\Big\rangle\geq 0.
\]
\end{lemma}

\smallskip

\noindent{\em Proof.} Assume for a contradiction that
\[
\liminf_{\overline{E}\ni y\rightarrow x}\Big\langle\frac{y-x}{|y-x|},n(x)\Big\rangle\in[-1,0).
\]
There exists $\eta\in(0,1)$ and a sequence $(y_h)$ in $E$ such that $y_h\rightarrow x$ as $h\rightarrow\infty$ and 
\begin{align}
\Big\langle\frac{y_h-x}{|y_h-x|},n(x)\Big\rangle&<-\eta\label{sequence_and_normal}
\end{align}
for each $h\in\mathbb{N}$. Choose $\alpha\in(0,\pi/2)$ such that $\cos\alpha=\eta$. As $M$ is $C^1$ there exists $r>0$ such that
\[
B(x,r)\cap C(x,-n(x),\alpha)\cap E=\emptyset.
\]
By choosing $h$ sufficiently large we can find $y_h\in B(x,r)$ with the additional property that $y_h\in C(x,-n(x),\alpha)$ by (\ref{sequence_and_normal}). We are thus led to a contradiction.
\qed

\smallskip

\begin{lemma}\label{properties_of_section_of_boundary_of_E}
Let $E$ be an open set in $\mathbb{R}^2$ with $C^1$ boundary $M$ and assume that $E\setminus\{0\}=E^{sc}$. For each $0<\tau\in \pi(M)$,
\begin{itemize}
\item[(i)] $|\cos\sigma|$ is constant on $M_\tau$;
\item[(ii)] $\cos\sigma=0$ on $M_\tau\cap\{x_2=0\}$;
\item[(iii)] $\langle Ox,n(x)\rangle\leq 0$ for $x\in M_\tau\cap H$ 
\item[(iv)] $\cos\sigma\leq 0$ on $M_\tau\cap H$;
\end{itemize}
and if $\cos\sigma\not\equiv0$ on $M_\tau$ then
\begin{itemize}
\item[(v)] $\tau\in p(E)$;
\item[(vi)]  $M_\tau$ consists of two disjoint singletons in $\mathbb{S}^1_\tau$ symmetric about the $x_1$-axis;
\item[(vii)] $L(\tau)\in(0,2\pi\tau)$;
\item[(viii)] $M_\tau=\{(\tau\cos(L(\tau)/2\tau),\pm\tau\sin(L(\tau)/2\tau)\}$.
\end{itemize}
\end{lemma}

\smallskip

\noindent{\em Proof.} 
{\em (i)} By Lemma \ref{property_of_E_sc}, $M_\tau$ is the union of two closed  spherical arcs in $\mathbb{S}^1_\tau$ symmetric about the $x_1$-axis. In case $M_\tau\cap\overline{H}$ consists of a singleton the assertion follows from Lemma \ref{oddness_of_cos_sigma}. Now suppose that $M_\tau\cap\overline{H}$ consists of a spherical arc in $\mathbb{S}^1_\tau$ with non-empty interior. It can be seen that $\cos\sigma$ vanishes on the interior of this arc as $0=r_1^\prime=\cos\sigma_1$ in a local parametrisation by (\ref{cos_of_sigma}). By continuity $\cos\sigma=0$ on $M_\tau$. {\em (ii)} follows from Lemma \ref{oddness_of_cos_sigma}. 
{\em (iii)} Let $x\in M_\tau\cap H$ so $\theta(x)\in(0,\pi)$. Then $S(\theta(x))\cap\mathbb{S}^1_\tau\subset\overline{E}$ as $\overline{E}$ is spherical cap symmetric. Then
\[
0\leq\lim_{S(\theta(x))\cap\mathbb{S}^1_\tau\ni y\rightarrow x}\Big\langle\frac{y-x}{|y-x|},n(x)\Big\rangle=-\langle Ox,n(x)\rangle
\]
by Lemma \ref{property_of_n}. {\em (iv)} The adjoint transformation $O^\star$ represents rotation clockwise through $\pi/2$. Let $x\in M_\tau\cap H$. By {\em (iii)},
\[
0\geq\langle Ox,n(x)\rangle=\langle x,O^\star n(x)\rangle=\langle x,t(x)\rangle=\tau\cos\sigma(x)
\]
and this leads to the result. {\em (v)} As $\cos\sigma\not\equiv0$ on $M_\tau$ we can find $x\in M_\tau\cap H$. We claim that $\mathbb{S}^1_\tau\cap S(\theta(x))\subset E$. For suppose that $y\in\mathbb{S}^1_\tau\cap S(\theta(x))$ but $y\not\in E$. We may suppose that $0\leq\theta(y)<\theta(x)<\pi$. If $y\in\mathbb{R}^2\setminus\overline{E}$ then $x\in\mathbb{R}^2\setminus\overline{E}$ by Lemma \ref{property_of_E_sc}. On the other hand, if $y\in M$ then the spherical arc in $H$ joining $y$ to $x$ is contained in $M$ again by Lemma \ref{property_of_E_sc}. This arc also has non-empty interior in $\mathbb{S}^1_\tau$. Now $\cos\sigma=0$ on its interior so $\cos(\sigma(x))=0$ by {\em (i)} contradicting the hypothesis. A similar argument deals with {\em (vi)} and this together with {\em (v)} in turn entails {\em (vii)} and {\em (viii)}.
\qed

\smallskip

\begin{lemma}\label{behaviour_at_0}
Let $E$ be an open set in $\mathbb{R}^2$ with $C^1$ boundary $M$ and assume that $E\setminus\{0\}=E^{sc}$. Suppose that $0\in M$. Then 
\begin{itemize}
\item[(i)] $(\sin\sigma)(0+)=0$;
\item[(ii)] $(\cos\sigma)(0+)=-1$.
\end{itemize}
\end{lemma}

\smallskip

\noindent{\em Proof.}
{\em (i)} Let $\gamma_1$ be a $C^1$ parametrisation of $M$ in a neighbourhood of $0$ with $\gamma_1(0)=0$ as above. Then $n(0)=n_1(0)=e_1$ and hence $t(0)=t_1(0)=-e_2$. By Taylor's Theorem $\gamma_1(s)=\gamma_1(0)+t_1(0)s+o(s)=-e_2s+o(s)$ for $s\in I$. This means that $r_1(s)=|\gamma_1(s)|=s+o(s)$ and 
\[
\cos\theta_1=\frac{\langle e_1,\gamma_1\rangle}{r_1}=\frac{\langle e_1,\gamma_1\rangle}{s}\frac{s}{r_1}\rightarrow 0
\]
as $s\rightarrow 0$ which entails that $(\cos\theta_1)(0-)=0$. Now $t_1$ is continuous on $I$ so $t_1=-e_2+o(1)$ and $\cos\alpha_1=\langle e_1,t_1\rangle=o(1)$. We infer that $(\cos\alpha_1)(0-)=0$. By (\ref{alpha_1_as_sum}), $\cos\alpha_1=\cos\sigma_1\cos\theta_1-\sin\sigma_1\sin\theta_1$ on $I$ and hence $(\sin\sigma_1)(0-)=0$. We deduce that $(\sin\sigma)(0+)=0$. Item {\em (ii)} follows from {\em (i)} and Lemma \ref{properties_of_section_of_boundary_of_E}.
\qed

\smallskip

\noindent The set
\begin{equation}\label{definition_of_Omega}
\Omega:=\pi\Big[(M\setminus\{0\})\cap\{\cos\sigma\neq 0\}\Big]
\end{equation}
plays an important r\^{o}le in the proof of Theorem \ref{main_theorem}. 

\smallskip

\begin{lemma}\label{Omega_is_open}
Let $E$ be an open set in $\mathbb{R}^2$ with $C^1$ boundary $M$ and assume that $E\setminus\{0\}=E^{sc}$.
Then $\Omega$ is an open set in $(0,+\infty)$.
\end{lemma}

\smallskip

\noindent{\em Proof.}
Suppose $0<\tau\in\Omega$. Choose $x\in M_\tau\cap\{\cos\sigma\neq 0\}$. Let $\gamma_1:I\rightarrow M$ be a local $C^1$ parametrisation of $M$ in a neighbourhood of $x$ such that $\gamma_1(0)=x$ as before. By shrinking $I$ if necessary we may assume that $r_1\neq 0$ and $\cos\sigma_1\neq 0$ on $I$. 
Then the set $\{r_1(s):s\in I\}\subset\Omega$ is connected and so an interval in $\mathbb{R}$ (see for example \cite{Simmons1963} Theorems 6.A and 6.B). By (\ref{cos_of_sigma}), $r_1^\prime(0)=\cos\sigma_1(0)=\cos\sigma(p)\neq 0$. This means that the set $\{r_1(s):s\in I\}$ contains an open interval about $\tau$.
\qed

\section{Generalised (mean) curvature}

\smallskip

\noindent Given a set $E$ of finite perimeter in $\mathbb{R}^2$ the first variation $\delta V_f(Z)$ resp. $\delta P_f^+(Z)$ of weighted volume and perimeter along a time-dependent vector field $Z$ are defined as in (\ref{first_variation_of_volume}) and (\ref{first_variation_of_perimeter}). 

\smallskip

\begin{proposition}\label{variation_of_perimeter_along_time_dependent_vector_field}
Let $f$ be as in (\ref{form_of_f}) where $h:[0,+\infty)\rightarrow\mathbb{R}$ is a non-decreasing convex function. Let $E$ be a bounded open set in $\mathbb{R}^2$ with $C^1$ boundary $M$. Let $Z$ be a time-dependent vector field. Then
\[
\delta P_f^+(Z)=\int_M f^\prime_+(\cdot,Z_0) + f\mathrm{div}^M Z_0\,d\mathscr{H}^1
\]
where $Z_0:=Z(0,\cdot)\in C^1_c(\mathbb{R}^2,\mathbb{R}^2)$. 
\end{proposition}

\smallskip

\noindent{\em Proof.} The identity (\ref{difference_of_perimeters}) holds for each $t\in I$ with $M$ in place of $\mathscr{F}E$. The assertion follows on appealing to Lemma \ref{expansion_of_time_dependent_flow} and Lemma \ref{differentiability_of_density} with the help of the dominated convergence theorem.
\qed

\smallskip

\noindent Given $X,Y\in C^\infty_c(\mathbb{R}^2,\mathbb{R}^2)$ let $\psi$ resp. $\chi$ stand for the $1$-parameter group of $C^\infty$ diffeomorphisms of $\mathbb{R}^2$ associated to the vector fields $X$ resp. $Y$ as in (\ref{diffeomorphism_group}). Let $I$ be an open interval in $\mathbb{R}$ containing the point $0$. Suppose that the function $\sigma:I\rightarrow\mathbb{R}$ is $C^1$. Define a flow via
\[
\varphi:I\times\mathbb{R}^2\rightarrow\mathbb{R}^2;
(t,x)\mapsto\chi(\sigma(t),\psi(t,x)).
\]

\smallskip

\begin{lemma}\label{time_dependent_vector_field_associated_with_varphi}
The time-dependent vector field $Z$ associated with the flow $\varphi$ is given by
\begin{align}\label{time_dependent_vector_field}
Z(t,x)=\sigma^\prime(t)Y(\chi(\sigma(t),\psi(t,x)))+d\chi(\sigma(t),\psi(t,x))X(\psi(t,x))
\end{align}
for $(t,x)\in I\times\mathbb{R}^2$ and satisfies (Z.1) and (Z.2).
\end{lemma}

\smallskip

\noindent{\em Proof.}
For $t\in I$ and $x\in\mathbb{R}^2$ we compute using (\ref{diffeomorphism_group}),
\begin{align}
\partial_t\varphi(t,x)
&=(\partial_t\chi)(\sigma(t),\psi(t,x))\sigma^\prime(t)+d\chi(\sigma(t),\psi(t,x))\partial_t\psi(t,x)\nonumber
\end{align}
and this gives (\ref{time_dependent_vector_field}). Put $K_1:=\mathrm{supp}[X]$, $K_2:=\mathrm{supp}[Y]$ and $K:=K_1\cup K_2$. Then $(Z.2)$ holds with this choice of $K$.
\qed

\smallskip

\noindent Let $E$ be a bounded open set in $\mathbb{R}^2$ with $C^1$ boundary $M$. Define $\Lambda:=(M\setminus\{0\})\cap\{\cos\sigma=0\}$ and 
\begin{align}
\Lambda_1&:=\{x\in M:\mathscr{H}^1(\Lambda\cap B(x,\rho))=\mathscr{H}^1(M\cap B(x,\rho))\text{ for some }\rho>0\}.\label{definition_of_Lambda_1}
\end{align}
For future reference put $\Lambda_1^{\pm}:=\Lambda_1\cap\{x\in M:\pm\langle x,n\rangle>0\}$.

\smallskip

\begin{lemma}\label{differentiability_of_f}
Let $f$ be as in (\ref{form_of_f}) where $h:[0,+\infty)\rightarrow\mathbb{R}$ is a non-decreasing convex function. Let $E$ be a bounded open set in $\mathbb{R}^2$ with $C^{1,1}$ boundary $M$ and suppose that $E\setminus\{0\}=E^{sc}$. Then
\begin{itemize}
\item[(i)] $\Lambda_1$ is a countable disjoint union of well-separated open circular arcs centred at $0$;
\item[(ii)] $\mathscr{H}^1(\overline{\Lambda_1}\setminus\Lambda_1)=0$; 
\item[(iii)] $f$ is differentiable $\mathscr{H}^1$-a.e. on $M\setminus\overline{\Lambda_1}$.
\end{itemize}
\end{lemma}

\smallskip

\noindent The term well-separated in {\em (i)} means the following: if $\Gamma$ is an open circular arc in $\Lambda_1$ with $\Gamma\cap(\Lambda_1\setminus\Gamma)=\emptyset$ then $d(\Gamma,\Lambda_1\setminus\Gamma)>0$.

\smallskip

\noindent{\em Proof.}
{\em (i)} Let $x\in\Lambda_1$ and $\gamma_1:I\rightarrow M$ a $C^{1,1}$ parametrisation of $M$ near $x$. By shrinking $I$ if necessary we may assume that $\gamma_1(I)\subset M\cap B(x,\rho)$ with $\rho$ as in (\ref{definition_of_Lambda_1}). So $\cos\sigma=0$ $\mathscr{H}^1$-a.e. on $\gamma_1(I)$ and hence $\cos\sigma_1=0$ a.e. on $I$. This means that $\cos\sigma_1=0$ on $I$ as $\sigma_1\in C^{0,1}(I)$ and that $r_1$ is constant on $I$ by (\ref{cos_of_sigma}). Using (\ref{sin_of_sigma}) it can be seen that $\gamma_1(I)$ is an open circular arc centred at $0$. By compactness of $M$ it follows that $\Lambda_1$ is a countable disjoint union of open circular arcs centred on $0$. The well-separated property flows from the fact that $M$ is $C^1$. {\em (ii)} follows as a consequence of this property. {\em (iii)} Let $x\in M\setminus\overline{\Lambda_1}$ and $\gamma_1:I\rightarrow M$ a $C^{1,1}$ parametrisation of $M$ near $x$ with properties as before. We assume that $x$ lies in the upper half-plane $H$. By shrinking $I$ if necessary we may assume that $\gamma_1(I)\subset(M\setminus\overline{\Lambda_1})\cap H$. Let $s_1,s_2, s_3\in I$ with $s_1<s_2<s_3$. Then $y:=\gamma_1(s_2)\in M\setminus\overline{\Lambda_1}$. So $\mathscr{H}^1(M\cap\{\cos\sigma\neq 0\}\cap B(y,\rho))>0$ for each $\rho>0$. This means that for small $\eta>0$ the set $\gamma_1((s_2-\eta,s_2+\eta))\cap\{\cos\sigma\neq 0\}$ has positive $\mathscr{H}^1$-measure. Consequently, $r_1(s_3)-r_1(s_1)=\int_{s_1}^{s_3}\cos\sigma_1\,ds<0$ bearing in mind Lemma \ref{properties_of_section_of_boundary_of_E}. This shows that $r_1$ is strictly decreasing on $I$. So $h$ is differentiable a.e. on $r_1(I)\subset(0,+\infty)$ in virtue of the fact that $h$ is convex and hence locally Lipschitz. This entails {\em (iii)}.
\qed

\smallskip

\begin{proposition}\label{on_generalised_mean_curvature}
Let $f$ be as in (\ref{form_of_f}) where $h:[0,+\infty)\rightarrow\mathbb{R}$ is a non-decreasing convex function. Given $v>0$ let $E$ be a minimiser of (\ref{isoperimetric_problem}). Assume that $E$ is a bounded open set in $\mathbb{R}^2$ with $C^1$ boundary $M$ and suppose that $E\setminus\{0\}=E^{sc}$. Suppose that $M\setminus\overline{\Lambda_1}\neq\emptyset$. Then there exists $\lambda\in\mathbb{R}$ such that for any $X\in C^1_c(\mathbb{R}^2,\mathbb{R}^2)$,
\[
0\leq\int_M\Big\{f^\prime_+(\cdot,X) + f\,\mathrm{div}^M X-\lambda f\langle n, X\rangle\Big\} \,d\mathscr{H}^1.
\]
\end{proposition}

\smallskip

\noindent{\em Proof.}
Let $X\in C^\infty_c(\mathbb{R}^2,\mathbb{R}^2)$. Let $x\in M$  and $r>0$ such that $M\cap B(x,r)\subset M\setminus\overline{\Lambda_1}$. Choose $Y\in C^\infty_c(\mathbb{R}^2,\mathbb{R}^2)$ with $\mathrm{supp}[Y]\subset B(x,r)$ as in Lemma \ref{vector_field_lemma}. Let $\psi$ resp. $\chi$ stand for the $1$-parameter group of $C^\infty$ diffeomorphisms of $\mathbb{R}^2$ associated to the vector fields $X$ resp. $Y$ as in (\ref{diffeomorphism_group}). For each $(s,t)\in\mathbb{R}^2$ the set $\chi_s(\psi_t(E))$ is an open set in $\mathbb{R}^2$ with $C^1$ boundary and $\partial(\chi_s\circ\psi_t)(E)=(\chi_s\circ\psi_t)(M)$ by Lemma \ref{diffeomorhism_invariance_of_boundary}. Define
\begin{align}
V(s,t) &:=V_f(\chi_t(\psi_s(E)))-V_f(E),\nonumber\\
P(s,t) &:=P_f(\chi_t(\psi_s(E))),\nonumber
\end{align}
for $(s,t)\in\mathbb{R}^2$. We write $F=(\chi_t\circ\psi_s)(E)$. Arguing as in Proposition \ref{variation_of_volume},
\begin{align}
\partial_t V(s,t)&=\lim_{h\rightarrow 0}(1/h)\{V_f(\chi_h(F))-V_f(F)\}
=\int_F\mathrm{div}(fY)\,dx\nonumber\\
&=\int_E(\mathrm{div}(fY)\circ\chi_t\circ\psi_s)\,J_2 d(\chi_t\circ\psi_s)_x\,dx\nonumber
\end{align}
with an application of the area formula (cf. \cite{Ambrosio2000} Theorem 2.71). This last varies continuously in $(s,t)$. The same holds for partial differentiation with respect to $s$. Indeed, put $\eta:=\chi_t\circ\psi_s$. Then noting that $J_2d(\eta\circ\psi_h)=(J_2d\eta)\circ\psi_h J_2d\psi_h$ and using the dominated convergence theorem,
\begin{align}
\partial_s V(s,t)&=\lim_{h\rightarrow 0}(1/h)
\Big\{
V_f(\eta(\psi_h(E)))-V_f(\eta(E))
\Big\}\nonumber\\
&=\lim_{h\rightarrow 0}(1/h)
\Big\{
\int_E(f\circ\eta\circ\psi_h)J_2d(\eta\circ\psi_h)_x\,dx-\int_E(f\circ\eta)J_2d\eta_x\,dx
\Big\}\nonumber\\
&=\lim_{h\rightarrow 0}(1/h)
\Big\{
\int_E[(f\circ\eta\circ\psi_h)-(f\circ\eta)]J_2d(\eta\circ\psi_h)_x\,dx\nonumber\\
&+\int_E(f\circ\eta)[(J_2d\eta\circ\psi_h-J_2d\eta]J_2d\psi_h\,dx
+\int_E(f\circ\eta)J_2d\eta[J_2d\psi_h-1]\,dx
\Big\}\nonumber\\
&=\int_E\langle\nabla(f\circ\eta),X\rangle J_2d\eta_x\,dx
+\int_E(f\circ\eta)\langle\nabla J_2d\eta,X\rangle\,dx
+\int_E(f\circ\eta)J_2d\eta\,\mathrm{div}\,X\,dx
\nonumber
\end{align}
where the explanation for the last term can be found in the proof of Proposition \ref{variation_of_volume}. In this regard we note that $d(d\chi_t)$ (for example) is continuous on $I\times\mathbb{R}^2$ (cf. \cite{Ambrosio2000} Theorem 3.3 and Exercise 3.2) and in particular $\nabla J_2d\chi_t$ is continuous on $I\times\mathbb{R}^2$. The expression above also varies continuously in $(s,t)$ as can be seen with the help of the dominated convergence theorem. This means that $V(\cdot,\cdot)$ is continuously differentiable on $\mathbb{R}^2$. Note that
\[
\partial_tV(0,0)=\int_E\mathrm{div}(fY)\,dx=1
\]
by choice of $Y$. By the implicit function theorem there exists $\eta>0$ and a $C^1$ function $\sigma:(-\eta,\eta)\rightarrow\mathbb{R}$ such that $\sigma(0)=0$ and $V(s,\sigma(s))=0$ for $s\in(-\eta,\eta)$; moreover,
\[
\sigma^\prime(0)=-\partial_s V(0,0)
=-\int_E\Big\{\langle\nabla f,X\rangle+f\,\mathrm{div}\,X\Big\}\,dx
=-\int_E\mathrm{div}(fX)\,dx
=\int_M f\,\langle n, X\rangle\,d\mathscr{H}^1
\]
by the Gauss-Green formula (cf. \cite{Ambrosio2000} Theorem 3.36). 

\smallskip

\noindent The mapping 
\[
\varphi:(-\eta,\eta)\times\mathbb{R}^2\rightarrow\mathbb{R}^2;t\mapsto\chi(\sigma(t),\psi(t,x))
\]
satisfies conditions (F.1)-(F.4) above with $I=(-\eta,\eta)$ where the associated time-dependent vector field $Z$ is given as in (\ref{time_dependent_vector_field}) and satisfies (Z.1) and (Z.2); moreover, $Z_0=Z(0,\cdot)=\sigma^\prime(0)Y+X$. Note that $Z_0=X$ on $M\setminus B(x,r)$.

\smallskip

\noindent The mapping $I\rightarrow\mathbb{R};t\mapsto P_f(\varphi_t(E))$ is right-differentiable at $t=0$ as can be seen from Proposition \ref{variation_of_perimeter_along_time_dependent_vector_field} and has non-negative right-derivative there. By Proposition \ref{variation_of_perimeter_along_time_dependent_vector_field} and Lemma \ref{differentiability_of_f},
\begin{align}
0&\leq\delta P_f^+(Z)
=\int_{M} f^\prime_+(\cdot, Z_0)+ f\,\mathrm{div}^M Z_0\,d\mathscr{H}^1\nonumber\\
&=\int_{M\setminus\overline{\Lambda_1}} f^\prime_+(\cdot, Z_0) + f\,\mathrm{div}^M Z_0\,d\mathscr{H}^1
+\int_{\overline{\Lambda_1}} f^\prime_+(\cdot, X) + f\,\mathrm{div}^M X\,d\mathscr{H}^1
\nonumber\\
&=\int_{M\setminus\overline{\Lambda_1}}\sigma^\prime(0)\langle\nabla f,Y\rangle+\langle\nabla f,X\rangle
+\sigma^\prime(0)\,f\,\mathrm{div}^M Y+f\,\mathrm{div}^M X\,d\mathscr{H}^1\nonumber\\
&+\int_{\overline{\Lambda_1}} f^\prime_+(\cdot, X) + f\,\mathrm{div}^M X\,d\mathscr{H}^1\nonumber\\
&=\int_{M}f^\prime_+(\cdot, X) + f\,\mathrm{div}^M X\,d\mathscr{H}^1
+\sigma^\prime(0)\int_{M}f^\prime_+(\cdot, Y)+ f\,\mathrm{div}^M Y\,d\mathscr{H}^1.\label{inequality_for_delta_P_plus}
\end{align}
The identity then follows upon inserting the expression for $\sigma^\prime(0)$ above with $\lambda=-\int_{M}f^\prime_+(\cdot,Y) + f\,\mathrm{div}^M Y\,d\mathscr{H}^1$. The claim follows for $X\in C^1_c(\mathbb{R}^2,\mathbb{R}^2)$ by a density argument.
\qed

\smallskip

\begin{theorem}\label{constant_weighted_mean_curvature}
Let $f$ be as in (\ref{form_of_f}) where $h:[0,+\infty)\rightarrow\mathbb{R}$ is a non-decreasing convex function. Given $v>0$ let $E$ be a minimiser of (\ref{isoperimetric_problem}). Assume that $E$ is a bounded open set in $\mathbb{R}^2$ with $C^{1,1}$ boundary $M$ and suppose that $E\setminus\{0\}=E^{sc}$. Suppose that $M\setminus\overline{\Lambda_1}\neq\emptyset$. Then there exists $\lambda\in\mathbb{R}$ such that 
\begin{itemize}
\item[(i)] $k+\varrho\sin\sigma+\lambda=0$ $\mathscr{H}^1$-a.e. on $M\setminus\overline{\Lambda_1}$;
\item[(ii)] $\varrho_--\lambda\leq k\leq\varrho_+-\lambda$ on $\Lambda_1^+$;
\item[(iii)] $-\varrho_+-\lambda\leq k\leq-\varrho_--\lambda$ on $\Lambda_1^-$.
\end{itemize}
\end{theorem}

\smallskip

\noindent The expression $k+\varrho\sin\sigma$ is called the generalised (mean) curvature of $M$.

\smallskip

\noindent{\em Proof.}
{\em (i)} Let $x\in M$ and $r>0$ such that $M\cap B(x,r)\subset M\setminus\overline{\Lambda_1}$. 
Choose $X\in C^1_c(\mathbb{R}^2,\mathbb{R}^2)$ with $\mathrm{supp}[X]\subset B(x,r)$. We know from Lemma \ref{differentiability_of_f} that $f$ is differentiable $\mathscr{H}^1$-a.e. on $\mathrm{supp}[X]$. Let $\lambda$ be as in Proposition \ref{on_generalised_mean_curvature}. Replacing $X$ by $-X$ we deduce from Proposition \ref{on_generalised_mean_curvature} that
\[
0=\int_M\Big\{\langle\nabla f,X\rangle + f\,\mathrm{div}^M X-\lambda f\langle n, X\rangle\Big\} \,d\mathscr{H}^1.
\]
The divergence theorem on manifolds (cf. \cite{Ambrosio2000} Theorem 7.34) holds also for $C^{1,1}$ manifolds. So
\begin{align}
\int_M\langle\nabla f, X\rangle + f\,\mathrm{div}^MX\,d\mathscr{H}^1
&=\int_M
\partial_nf\,\langle n,\,X\rangle 
+\langle\nabla^M f,\,X\rangle 
+f\,\mathrm{div}^M\,X
\,d\mathscr{H}^1\nonumber\\
&=\int_M
\partial_n f\,\langle n,\,X\rangle
+\mathrm{div}^M(fX)
\,d\mathscr{H}^1\nonumber\\
&=\int_M
\partial_n f\,\langle n,\,X\rangle
-Hf\langle n,\,X\rangle
\,d\mathscr{H}^1
=\int_M
fu
\left\{
\partial_n \log f-H
\right\}\,d\mathscr{H}^1\nonumber
\end{align}
where $u=\langle n,X\rangle$. Combining this with the equality above we see that 
\[
\int_M uf
\left\{
\partial_n\log f-H-\lambda
\right\}
\,d\mathscr{H}^1=0
\]
for all $X\in C^1_c(\mathbb{R}^2,\mathbb{R}^2)$. This leads to the result. 

\smallskip

\noindent{\em (ii)}  Let $x\in M$ and $r>0$ such that $M\cap B(x,r)\subset\Lambda_1^+$. Let $\phi\in C^1(\mathbb{S}^1_r)$ with support in $\mathbb{S}^1_r\cap B(x,r)$. We can construct $X\in C^1_c(\mathbb{R}^2,\mathbb{R}^2)$ with the property that $X=\phi n$ on $M\cap B(x,r)$. By Lemma \ref{differentiability_of_density}, 
\begin{align}
f^\prime_+(\cdot,X)=fh^\prime_+(|x|,\mathrm{sgn}\langle x,X\rangle)|\langle n,X\rangle|
=fh^\prime_+(|x|,\mathrm{sgn}\,\phi\langle x,n\rangle)|\phi|
\nonumber
\end{align}
on $\Lambda_1$. Let us assume that $\phi\geq 0$. As $\langle\cdot,n\rangle>0$ on $\Lambda_1^+$ we have that $f^\prime_+(\cdot,X)=f\phi h^\prime_+(|x|,+1)=f\phi\varrho_+$ so by Proposition \ref{on_generalised_mean_curvature},
\begin{align}
0&\leq\int_M\Big\{f^\prime_+(\cdot,X) + f\,\mathrm{div}^M X-\lambda f\langle n, X\rangle\Big\} \,d\mathscr{H}^1=\int_M f\phi\Big\{\varrho_+-k-\lambda\Big\} \,d\mathscr{H}^1.\nonumber
\end{align}
We conclude that $\varrho_+-k-\lambda\geq 0$ on $M\cap B(x,r)$. Now assume that $\phi\leq 0$. Then $f^\prime_+(\cdot,X)=-f\phi h^\prime_+(|x|,-1)=f\phi\varrho_-$ so
\begin{align}
0&\leq\int_M f\phi\Big\{\varrho_--k-\lambda\Big\} \,d\mathscr{H}^1\nonumber
\end{align}
and hence $\varrho_--k-\lambda\leq 0$ on $M\cap B(x,r)$. This shows {\em (ii)}.

\smallskip

\noindent{\em (iii)} The argument is similar. Assume in the first instance that $\phi\geq 0$. Then $f^\prime_+(\cdot,X)=f\phi h^\prime_+(|x|,-1)=-f\phi\varrho_-$ so
\begin{align}
0&\leq\int_M f\phi\Big\{-\varrho_--k-\lambda\Big\} \,d\mathscr{H}^1.\nonumber
\end{align}
We conclude that $-\varrho_--k-\lambda\geq 0$ on $M\cap B(x,r)$. Next suppose that $\phi\leq 0$. Then $f^\prime_+(\cdot,X)=-f\phi h^\prime_+(|x|,+1)=-f\phi\varrho_+$ so
\begin{align}
0&\leq\int_M f\phi\Big\{-\varrho_+-k-\lambda\Big\} \,d\mathscr{H}^1\nonumber
\end{align}
and $-\varrho_+-k-\lambda\leq 0$ on $M\cap B(x,r)$. 
\qed

\smallskip

\noindent Let $E$ be an open set in $\mathbb{R}^2$ with $C^1$ boundary $M$ and assume that $E\setminus\{0\}=E^{sc}$ and that $\Omega$ is as in (\ref{definition_of_Omega}). Bearing in mind Lemma \ref{properties_of_section_of_boundary_of_E} we may define
\begin{align}
&\theta_2:\Omega\rightarrow(0,\pi);\tau\mapsto L(\tau)/2\tau;\label{definition_of_theta}\\
&\gamma:\Omega\rightarrow M;\tau\mapsto(\tau\cos\theta_2(\tau),\tau\sin\theta_2(\tau)).\label{definition_of_gamma}
\end{align}
The function 
\begin{align}
&u:\Omega\rightarrow[-1,1];\tau\mapsto\sin(\sigma(\gamma(\tau))).\label{definition_of_y}
\end{align}
plays a key role.

\smallskip

\begin{theorem}\label{ode_for_y}
Let $f$ be as in (\ref{form_of_f}) where $h:[0,+\infty)\rightarrow\mathbb{R}$ is a non-decreasing convex function. Given $v>0$ let $E$ be a bounded minimiser of (\ref{isoperimetric_problem}). Assume that $E$ is open with $C^{1,1}$ boundary $M$ and that $E\setminus\{0\}=E^{sc}$. Suppose that $M\setminus\overline{\Lambda_1}\neq\emptyset$ and let $\lambda$ be as in Theorem \ref{constant_weighted_mean_curvature}. Then $u\in C^{0,1}(\Omega)$ and
\[
u^\prime+(1/\tau+\varrho)u+\lambda=0
\]
a.e. on $\Omega$.
\end{theorem}

\smallskip

\noindent{\em Proof.}
Let $\tau\in\Omega$ and $x$ a point in the open upper half-plane such that $x\in M_\tau$. There exists a $C^{1,1}$ parametrisation $\gamma_1:I\rightarrow M$ of $M$ in a neighbourhood of $x$ with $\gamma_1(0)=x$ as above. Put $u_1:=\sin\sigma_1$ on $I$. By shrinking the open interval $I$ if necessary we may assume that $r_1:I\rightarrow r_1(I)$ is a diffeomorphism and that $r_1(I)\subset\subset\Omega$. Note that $\gamma=\gamma_1\circ r_1^{-1}$ and $u=u_1\circ r_1^{-1}$ on $r_1(I)$. It follows that $u\in C^{0,1}(\Omega)$. By (\ref{cos_of_sigma}),
\[
u^\prime=\frac{\dot{u}_1}{\dot{r}_1}\circ r_1^{-1}
=\dot{\sigma}_1\circ r_1^{-1}
\]
a.e. on $r_1(I)$. As $\dot{\alpha}_1=k_1$ a.e. on $I$ and using the identity (\ref{sin_of_sigma}) we see that $\dot{\sigma}_1=\dot{\alpha}_1-\dot{\theta}_1=k_1-(1/r_1)\sin\sigma_1$ a.e on $I$. Thus,
\[
u^\prime=k-(1/\tau)\sin(\sigma\circ\gamma)=k-(1/\tau)u
\]
a.e. on $r_1(I)$. By Theorem \ref{constant_weighted_mean_curvature} there exists $\lambda\in\mathbb{R}$ such that $k+\varrho\sin\sigma+\lambda=0$ $\mathscr{H}^1$-a.e. on $M$. So
\[
u^\prime=-\varrho(\tau)u-\lambda-(1/\tau)u
=-(1/\tau+\varrho(\tau))u-\lambda
\]
a.e. on $r_1(I)$. The result follows.
\qed

\smallskip

\begin{lemma}\label{derivative_of_theta}
Suppose that $E$ is a bounded open set in $\mathbb{R}^2$ with $C^1$ boundary $M$ and that $E\setminus\{0\}=E^{sc}$. Then 
\begin{itemize}
\item[(i)] $\theta_2\in C^1(\Omega)$;
\item[(ii)] $\theta_2^\prime=-\frac{1}{\tau}\frac{u}{\sqrt{1-u^2}}$ on $\Omega$.
\end{itemize}
\end{lemma}

\smallskip

\noindent{\em Proof.}
Let $\tau\in\Omega$ and $x$ a point in the open upper half-plane such that $x\in M_\tau$. There exists a $C^1$ parametrisation $\gamma_1:I\rightarrow M$ of $M$ in a neighbourhood of $x$ with $\gamma_1(0)=x$ as above. By shrinking the open interval $I$ if necessary we may assume that $r_1:I\rightarrow r_1(I)$ is a diffeomorphism and that $r_1(I)\subset\subset\Omega$. It then holds that
\[
\theta_2=\theta_1\circ r_1^{-1}\text{ and }\sigma\circ\gamma=\sigma_1\circ r_1^{-1}
\]
on $r_1(I)$ by choosing an appropriate branch of $\theta_1$. It follows that $\theta_2\in C^1(\Omega)$. By the chain-rule, (\ref{sin_of_sigma}) and (\ref{cos_of_sigma}),
\[
\theta_2^\prime=\frac{\dot{\theta}_1}{\dot{r}_1}\circ r_1^{-1}
=(\frac{1}{r_1}\tan\sigma_1)\circ r_1^{-1}
=(1/\tau)\tan(\sigma\circ\gamma)
\] 
on $r_1(I)$. By Lemma \ref{properties_of_section_of_boundary_of_E}, $\cos(\sigma\circ\gamma)=-\sqrt{1-u^2}$ on $\Omega$. This entails {\em (ii)}.
\qed

\section{Convexity}

\smallskip

\begin{lemma}\label{lemma_on_lb_for_curvature}
Let $E$ be a bounded open set in $\mathbb{R}^2$ with $C^{1,1}$ boundary $M$ and assume that $E\setminus\{0\}=E^{sc}$. Put $d:=\sup\{|x|:x\in M\}>0$ and $b:=(d,0)$. Let $\gamma_1:I\rightarrow M$ be a $C^{1,1}$ parametrisation of $M$ near $b$ with $\gamma_1(0)=b$. Then
\[
\lim_{\delta\downarrow 0}\Big\{\mathrm{ess\,sup}_{[-\delta,\delta]}k_1\Big\}\geq 1/d.
\]
\end{lemma}

\smallskip

\noindent{\em Proof.} For $s\in I$,
\[
\gamma_1(s)=de_1+se_2+\int_0^s\Big\{\dot{\gamma}_1(u)-\dot{\gamma}_1(0)\Big\}\,du
\]
and
\[
\dot{\gamma}_1(u)-\dot{\gamma}_1(0)=\int_0^u k_1 n_1\,dv
\]
by (\ref{difference_of_dot_gamma}). By the Fubini-Tonelli Theorem,
\[
\gamma_1(s)=de_1+se_2+\int_0^s(s-u)k_1(u)n_1(u)\,du
=de_1+se_2+R(s)
\]
for $s\in I$. Assume for a contradiction that
\[
\lim_{\delta\downarrow 0}\Big\{\mathrm{ess\,sup}_{[-\delta,\delta]}k_1\Big\}<l< 1/d
\]
for some $l\in\mathbb{R}$. Then we can find $\delta>0$ such that $k_1<l$ a.e. on $[-\delta,\delta]$. So
\begin{align}
\langle R(s),e_1\rangle&=\int_0^s(s-u)k_1(u)\langle n_1(u),e_1\rangle\,du
>-(1/2)s^2l(1+o(1))\nonumber
\end{align}
as $s\downarrow 0$ and 
\begin{align}
r_1(s)^2-d^2&=2d\langle R(s),e_1\rangle+s^2+o(s^2)
>-dls^2(1+o(1))+s^2+o(s^2)\nonumber
\end{align}
as $s\downarrow 0$. Alternatively,
\begin{align}
\frac{r_1(s)^2-d^2}{s^2}&>1-dl+o(1).\nonumber
\end{align}
As $1-dl>0$ we can find $s\in I$ with $r_1(s)>d$, contradicting the definition of $d$.
\qed

\smallskip

\begin{lemma}\label{upper_bound_for_lambda}
Let $f$ be as in (\ref{form_of_f}) where $h:[0,+\infty)\rightarrow\mathbb{R}$ is a non-decreasing convex function. Given $v>0$ let $E$ be a bounded minimiser of (\ref{isoperimetric_problem}). Assume that $E$ is open with $C^{1,1}$ boundary $M$ and that $E\setminus\{0\}=E^{sc}$. Suppose that $M\setminus\overline{\Lambda_1}\neq\emptyset$. Then $\lambda\leq-1/d-\varrho_-(d)<0$ with $\lambda$ as in Theorem \ref{constant_weighted_mean_curvature}. 
\end{lemma}

\smallskip

\noindent{\em Proof.} We write $M$ as the disjoint union $M=(M\setminus\overline{\Lambda_1})\cup\overline{\Lambda_1}$. Let $b$ be as above. Suppose that $b\in\overline{\Lambda_1}$. Then $b\in\Lambda_1$; in fact, $b\in\Lambda_1^-$. By Theorem \ref{constant_weighted_mean_curvature}, $\lambda\leq-\varrho_--k$ at $b$. By Lemma \ref{lemma_on_lb_for_curvature}, $\lambda\leq-1/d-\varrho_-(d)$ upon considering an appropriate sequence in $M$ converging to $b$. Now suppose that $b$ lies in the open set $M\setminus\overline{\Lambda_1}$ in $M$.  Let $\gamma_1:I\rightarrow M$ be a $C^{1,1}$ parametrisation of $M$ near $b$ with $\gamma_1(I)\subset M\setminus\overline{\Lambda_1}$. By Theorem \ref{constant_weighted_mean_curvature}, $k_1+\varrho(r_1)\sin\sigma_1+\lambda=0$ a.e. on $I$. Now $\sin\sigma_1(s)\rightarrow 1$ as $s\rightarrow 0$. In light of Lemma \ref{lemma_on_lb_for_curvature}, $1/d+\varrho(d-)+\lambda\leq 0$ and $\lambda\leq-1/d-\varrho_-(d)$.
\qed

\smallskip

\begin{theorem}\label{E_is_convex}
Let $f$ be as in (\ref{form_of_f}) where $h:[0,+\infty)\rightarrow\mathbb{R}$ is a non-decreasing convex function. Given $v>0$ let $E$ be a bounded minimiser of (\ref{isoperimetric_problem}). Assume that $E$ is open with $C^{1,1}$ boundary $M$ and that $E\setminus\{0\}=E^{sc}$. Suppose that $M\setminus\overline{\Lambda_1}\neq\emptyset$. Then $E$ is convex. 
\end{theorem}

\smallskip

\noindent{\em Proof.}
The proof runs along similar lines as \cite{MorganPratelli2011} Theorem 6.5. By Theorem \ref{constant_weighted_mean_curvature}, $k+\varrho\sin\sigma+\lambda=0$ $\mathscr{H}^1$-a.e. on $M\setminus\overline{\Lambda_1}$. By Lemma \ref{upper_bound_for_lambda}, 
\[
0\leq k+\varrho_-(d)+\lambda\leq k-1/d
\]
and $k\geq 1/d$ $\mathscr{H}^1$-a.e. on $M\setminus\overline{\Lambda_1}$. On  $\Lambda_1^+$, $k\geq\varrho_--\lambda\geq\varrho_-+\varrho_-(d)+1/d>0$; on the other hand, $k<0$ on $\Lambda_1^+$. So in fact $\Lambda_1^+=\emptyset$. If $b\in\Lambda_1^-$ then $k=1/d$. On $\Lambda_1^-\cap B(0,d)$, $k\geq-\varrho_+-\lambda\geq-\varrho_++\varrho_-(d)+1/d\geq 1/d$. Therefore $k\geq 1/d>0$ $\mathscr{H}^1$-a.e. on $M$. The set $E$ is then convex by a modification of \cite{Spivak1999} Theorem 1.8 and Proposition 1.4. It is sufficient that the function $f$ (here $\alpha_1$) in the proof of the former theorem is non-decreasing.
\qed

\smallskip

\section{A reverse Hermite-Hadamard inequality}

\smallskip

\noindent Let $0\leq a<b<+\infty$ and $\varrho\geq 0$ be a non-decreasing bounded function on $[a,b]$. Let $h$ be a primitive of $\varrho$ on $[a,b]$ so that $h\in C^{0,1}([a,b])$ and introduce the functions
\begin{align}
&\mathtt{f}:[a,b]\rightarrow\mathbb{R};x\mapsto e^{h(x)};\\
&g:[a,b]\rightarrow\mathbb{R};x\mapsto x\mathtt{f}(x).
\end{align}
Then
\begin{align}\label{derivative_of_g}
g^\prime&=(1/x+\varrho)g=\mathtt{f}+g\varrho
\end{align}
a.e. on $(a,b)$.
Define
\begin{equation}\label{definition_of_m}
m=m(\varrho,a,b):=\frac{g(b)-g(a)}{\int_a^b g\,dt}.
\end{equation}
If $\varrho$ takes the constant value $\mathbb{R}\ni\lambda\geq 0$ on $[a,b]$ we use the notation $m(\lambda,a,b)$ and we write $m_0=m(0,a,b)$. A computation gives 
\begin{equation}\label{formula_for_m_0}
m_0=m(0,a,b)=A(a,b)^{-1}
\end{equation}
where $A(a,b):=(a+b)/2$ stands for the arithmetic mean of $a$ and $b$.

\smallskip

\begin{lemma}
Let $0\leq a<b<+\infty$ and $\varrho\geq 0$ be a non-decreasing bounded function on $[a,b]$. Then $m_0\leq m$.
\end{lemma}

\smallskip

\noindent{\em Proof.}
Note that $g$ is convex on $[a,b]$ as can be seen from (\ref{derivative_of_g}). By the Hermite-Hadamard inequality (cf.  \cite{Hermite1883}, \cite{Hadamard1893}),
\begin{align}
\frac{1}{b-a}\int_a^b g\,dt&\leq\frac{g(a)+g(b)}{2}.\label{Hermite_Hadamard_inequality}
\end{align}
The inequality $(b-a)(g(a)+g(b))\leq(a+b)(g(b)-g(a))$ entails
\[
\int_a^b g\,dt\leq\frac{a+b}{2}(g(b)-g(a))
\]
and the result follows on rearrangement.
\qed

\smallskip

\begin{lemma}\label{m_for_constant_varrho}
Let $0\leq a<b<+\infty$ and $\lambda>0$. Then $m(\lambda,a,b)<\lambda+A(a,b)^{-1}$.
\end{lemma}

\smallskip

\noindent{\em Proof.}
First suppose that $\lambda=1$ and take $h:[a,b]\rightarrow\mathbb{R};t\mapsto t$. In this case,
\[
\int_a^b g\,dt=\int_a^b te^t\,dt=(b-1)e^b-(a-1)e^a
\]
and
\[
m(1,a,b)=\frac{be^b-ae^a}{(b-1)e^b-(a-1)e^a}.
\]
The inequality in the statement is equivalent to
\[
(a+b)(be^b-ae^a)<((b-1)e^b-(a-1)e^a)(2+a+b)
\]
which in turn is equivalent to the statement $\tanh[(b-a)/2]<(b-a)/2$ which holds for any $b>a$.

\smallskip

\noindent For $\lambda>0$ take $h:[a,b]\rightarrow\mathbb{R};t\mapsto\lambda t$. Substitution gives
\[
\int_a^b g\,dt=(1/\lambda)^2[(\lambda b-1)e^{\lambda b}-(\lambda a-1)e^{\lambda a}]
\text{ and }
g(b)-g(a)=(1/\lambda)[\lambda be^{\lambda b}-\lambda ae^{\lambda a}]
\]
so from above
\[
m(\lambda,a,b)=\lambda m(1,\lambda a,\lambda b)
<\lambda\Big\{1+A(\lambda a,\lambda b)^{-1}\Big\}
=\lambda+A(a,b)^{-1}.
\]
\qed

\smallskip

\begin{theorem}\label{upper_bound_for_m}
Let $0\leq a<b<+\infty$ and $\varrho\geq 0$ be a non-decreasing bounded function on $[a,b]$. Then
\begin{itemize}
\item[(i)] $m(\varrho,a,b)\leq\varrho(b-)+A(a,b)^{-1}$;
\item[(ii)] equality holds if and only if $\varrho\equiv 0$ on $[a,b)$.
\end{itemize}
\end{theorem}

\smallskip

\noindent{\em Proof.}
{\em (i)} Define $h:=\int_a^\cdot\varrho\,d\tau$ on $[a,b]$ so that $h^\prime=\varrho$ a.e. on $(a,b)$. Define $h_1:[a,b]\rightarrow\mathbb{R};t\mapsto h(b)-\varrho(b-)(b-t)$. Then $h_1(b)=h(b)$, $h_1^\prime=\varrho(b-)\geq\varrho=h^\prime$ a.e. on $(a,b)$ and hence $h\geq h_1$ on $[a,b]$. We derive 
\[
\int_a^b g\,dt=\int_a^b te^{h(t)}\,dt\geq\int_a^b te^{h_1(t)}\,dt=\int_a^b g_1\,dt
\]
and
\[
g(b)-g(a)=be^{h(b)}-ae^{h(a)}=be^{h_1(b)}-ae^{h(a)}
\leq be^{h_1(b)}-ae^{h_1(a)}=g_1(b)-g_1(a)
\]
with obvious notation. This entails that $m(\varrho,a,b)\leq m(\varrho(b-),a,b)$ and the result follows with the help of Lemma \ref{m_for_constant_varrho}.

\smallskip

\noindent{\em (ii)} Suppose that $\varrho\not\equiv 0$ on $[a,b)$. If $\varrho$ is constant on $[a,b]$ the assertion follows from Lemma \ref{m_for_constant_varrho}. Assume then that $\varrho$ is not constant on $[a,b)$. Then $h\not\equiv h_1$ on $[a,b]$ in the above notation and $\int_a^b te^{h(t)}\,dt>\int_a^b te^{h_1(t)}\,dt$ which entails strict inequality in {\em (i)}.
\qed

\smallskip

\noindent With the above notation define
\begin{equation}
\hat{m}=\hat{m}(\varrho,a,b):=\frac{g(a)+g(b)}{\int_a^b g\,dt}.
\end{equation}
A computation gives 
\begin{equation}
\hat{m}_0:=\hat{m}(0,a,b)=\frac{2}{b-a}.
\end{equation}

\smallskip

\begin{lemma}
Let $0\leq a<b<+\infty$ and $\varrho\geq 0$ be a non-decreasing bounded function on $[a,b]$. Then $\hat{m}\geq\hat{m}_0$.
\end{lemma}

\smallskip

\noindent{\em Proof.}
This follows by the Hermite-Hadamard inequality (\ref{Hermite_Hadamard_inequality}). \qed

\smallskip

\noindent We prove a reverse Hermite-Hadamard inequality.

\begin{theorem}\label{reverse_Hermite_Hadamard_inequality}
Let $0\leq a<b<+\infty$ and $\varrho\geq 0$ be a non-decreasing bounded function on $[a,b]$. Then
\begin{itemize}
\item[(i)] $(b-a)\hat{m}(\varrho,a,b)\leq 2+a\varrho(a+)+b\varrho(b-)$;
\item[(ii)] equality holds if and only if $\varrho\equiv 0$ on $[a,b)$.
\end{itemize}
\end{theorem}

\smallskip

\noindent This last inequality can be written in the form
\begin{align}
\frac{g(a)+g(b)}{2+a\varrho(a+)+b\varrho(b-)}&\leq\frac{1}{b-a}\int_a^b g\,dt;\nonumber
\end{align}
comparing with (\ref{Hermite_Hadamard_inequality}) justifies naming this a reverse Hermite-Hadamard inequality. 

\smallskip

\noindent{\em Proof.} {\em (i)} We assume in the first instance that $\varrho\in C^1((a,b))$.  We prove the above result in the form
\begin{align}
\int_a^b g\,dt&\geq(b-a)\frac{g(a)+g(b)}{2+a\varrho(a)+b\varrho(b)}.\label{2nd_form_of_Hermite_Hadamard_inequality}
\end{align}
Put
\[
w:=\frac{(t-a)(g(a)+g)}{2+a\varrho(a)+t\varrho}
\]
for $t\in[a,b]$ so that
\[
\int_a^b w^\prime\,dt=(b-a)\frac{g(a)+g(b)}{2+a\varrho(a)+b\varrho(b)}.
\]
Then using (\ref{derivative_of_g}),
\begin{align}
w^\prime&=\frac{(g(a)+g+(t-a)g^\prime)(2+a\varrho(a)+t\varrho)-(t-a)(g(a)+g)(\varrho+t\varrho^\prime)}{(2+a\varrho(a)+t\varrho)^2}\nonumber\\
&=\frac{(g(a)-ag^\prime+(2+t\varrho)g)(2+a\varrho(a)+t\varrho)-(t-a)(g(a)+g)(\varrho+t\varrho^\prime)}{(2+a\varrho(a)+t\varrho)^2}\nonumber\\
&=\frac{(2+t\varrho)(2+a\varrho(a)+t\varrho)}{(2+a\varrho(a)+t\varrho)^2}g
+\frac{(g(a)-ag^\prime)(2+a\varrho(a)+t\varrho)-(t-a)(g(a)+g)(\varrho+t\varrho^\prime)}{(2+a\varrho(a)+t\varrho)^2}\nonumber\\
&\leq g
-\frac{2g(a)}{(2+a\varrho(a)+b\varrho(b))^2}(t-a)\varrho\label{differential_form_of_Hermite_Hadamard_inequality}
\\
&\leq g\nonumber
\end{align}
on $(a,b)$ as 
\[
g(a)-ag^\prime=a(\mathtt{f}(a)-(1/t+\varrho)g)=a(\mathtt{f}(a)-\mathtt{f}-\varrho g)\leq 0.
\]
An integration over $[a,b]$ gives the result. 

\smallskip

\noindent Let us now assume that $\varrho\geq 0$ is a non-decreasing bounded function on $[a,b]$. Extend $\varrho$ to $\mathbb{R}$ via
\[
\widetilde{\varrho}(t):=
\left\{
\begin{array}{lcl}
\varrho(a+) & \text{ for } & t\in(-\infty,a];\\
\varrho(t) & \text{ for } & t\in(a,b];\\
\varrho(b-) & \text{ for } & t\in(b,+\infty);\\
\end{array}
\right.
\]
for $t\in\mathbb{R}$. Let $(\psi_\varepsilon)_{\varepsilon>0}$ be a family of mollifiers (see e.g. \cite{Ambrosio2000} 2.1) and set $\widetilde{\varrho}_\varepsilon:=\widetilde{\varrho}\star\psi_\varepsilon$ on $\mathbb{R}$ for each $\varepsilon>0$. Then $\widetilde{\varrho}_\varepsilon\in C^\infty(\mathbb{R})$ and is non-decreasing on $\mathbb{R}$ for each $\varepsilon>0$. Put $\varrho_\varepsilon:=\widetilde{\varrho}_\varepsilon\mid_{[a,b]}$ for each $\varepsilon>0$. Then $(\varrho_\varepsilon)_{\varepsilon>0}$ converges to $\varrho$ in $L^1((a,b))$ by \cite{Ambrosio2000} 2.1 for example. Note that $h_{\varepsilon}:=\int_a^\cdot\varrho_\varepsilon\,dt\rightarrow h$ pointwise on $[a,b]$ as $\varepsilon\downarrow 0$ and that $(h_\varepsilon)$ is uniformly bounded on $[a,b]$. Moreover, 
$\varrho_\varepsilon(a)\rightarrow\varrho(a+)$ and $\varrho_\varepsilon(b)\rightarrow\varrho(b-)$ as $\varepsilon\downarrow 0$. By the above result,
\[
(b-a)\hat{m}(\varrho_{\varepsilon},a,b)\leq 2+a\varrho_{\varepsilon}(a)+b\varrho_{\varepsilon}(b)
\]
for each $\varepsilon>0$. The inequality follows on taking the limit $\varepsilon\downarrow 0$ with the help of the dominated convergence theorem.

\smallskip

\noindent{\em (ii)} We now consider the equality case. We claim that
\begin{align}
(b-a)\frac{g(a)+g(b)}{2+a\varrho(a+)+b\varrho(b-)}&\leq\int_a^b g\,dt-\frac{2g(a)}{(2+a\varrho(a+)+b\varrho(b-))^2}\int_a^b(t-a)\varrho\,dt;\label{on_equality_condition}
\end{align}
this entails the equality condition in {\em (ii)}. First suppose that $\varrho\in C^1((a,b))$. In this case the inequality in (\ref{differential_form_of_Hermite_Hadamard_inequality}) implies (\ref{on_equality_condition}) upon integration. Now suppose that $\varrho\geq 0$ is a non-decreasing bounded function on $[a,b]$. Then (\ref{on_equality_condition}) holds with $\varrho_\varepsilon$ in place of $\varrho$ for each $\varepsilon>0$. The inequality for $\varrho$ follows by the dominated convergence theorem.
\qed

\section{Comparison theorems for first-order differential equations}

\smallskip

\noindent Let $\mathscr{L}$ stand for the collection of Lebesgue measurable sets in $[0,+\infty)$. Define a measure $\mu$ on $([0,+\infty),\mathscr{L})$ by $\mu(dx):=(1/x)\,dx$. Let $0\leq a<b<+\infty$. Suppose that $u:[a,b]\rightarrow\mathbb{R}$ is an $\mathscr{L}^1$-measurable function with the property that
\begin{align}
\mu(\{u>t\})&<+\infty\text{ for each }t>0.
\end{align}
 The distribution function  $\mu_u:(0,+\infty)\rightarrow[0,+\infty)$ of $u$ with respect to $\mu$ is given by
\[
\mu_u(t):=\mu(\{u>t\})\text{ for }t>0.
\]
Note that $\mu_u$ is right-continuous and non-increasing on $(0,\infty)$ and $\mu_u(t)\rightarrow 0$ as $t\rightarrow\infty$. 

\smallskip

\noindent Let $u$ be a Lipschitz function on $[a,b]$. Define
\[
Z_1:=\{u\text{ differentiable and }u^\prime=0\},
Z_2:=\{u\text{ not differentiable}\}\text{ and }Z:=Z_1\cup Z_2.
\]
By \cite{Ambrosio2000} Lemma 2.96, $Z\cap\{u=t\}=\emptyset$ for $\mathscr{L}^1$-a.e. $t\in\mathbb{R}$ and hence $N:=u(Z)\subset\mathbb{R}$ is $\mathscr{L}^1$-negligible. We make use of the coarea formula (\cite{Ambrosio2000} Theorem 2.93 and (2.74)),
\begin{equation}\label{coarea_formula}
\int_{[a,b]}\phi|u^\prime|\,dx=\int_{-\infty}^\infty\int_{\{u=t\}}\phi\,d\mathscr{H}^{0}\,dt
\end{equation}
for any $\mathscr{L}^1$-measurable function $\phi:[a,b]\rightarrow[0,\infty]$. 

\smallskip

\begin{lemma}\label{properties_of_D_mu}
Let $0\leq a<b<+\infty$ and $u$ a Lipschitz function on $[a,b]$. Then
\begin{itemize}
\item[(i)] $\mu_u\in\mathrm{BV}_{\mathrm{loc}}((0,+\infty))$;
\item[(ii)] $D\mu_u=-u_\sharp\mu$;
\item[(iii)] $D\mu_u^a=D\mu_u\mres((0,+\infty)\setminus N)$;
\item[(iv)] $D\mu_u^s=D\mu_u\mres N$;
\item[(v)] $A:=\Big\{t\in(0,+\infty):\mathscr{L}^1(Z\cap\{u=t\})>0\Big\}$ is the set of atoms of $D\mu_u$ and $D\mu_u^j=D\mu_u\mres A$;
\item[(vi)] $\mu_u$ is differentiable $\mathscr{L}^1$-a.e. on $(0,+\infty)$ with derivative given by
\[
\mu_u^\prime(t)=-\int_{\{u=t\}\setminus Z}\frac{1}{|u^\prime|}\,\frac{d\mathscr{H}^{0}}{\tau}
\]
for $\mathscr{L}^1$-a.e. $t\in(0,+\infty)$;
\item[(vii)] $\mathrm{Ran}(u)\cap[0,+\infty)=\mathrm{supp}(D\mu_u)$.
\end{itemize}
\end{lemma}

\smallskip

\noindent The notation above $D\mu_u^a$, $D\mu_u^s$, $D\mu_u^j$ stands for the absolutely continuous resp. singular resp. jump part of the measure $D\mu_u$ (see \cite{Ambrosio2000} 3.2 for example).

\smallskip

\noindent{\em Proof.}
For any $\varphi\in C^\infty_c((0,+\infty))$ with $\mathrm{supp}[\varphi]\subset(\tau,+\infty)$ for some $\tau>0$,
\begin{align}
\int_0^\infty\mu_u\varphi^\prime\,dt&=\int_{[a,b]}\varphi\circ u\,d\mu=\int_{[a,b]}\chi_{\{u>\tau\}}\varphi\circ u\,d\mu
\label{integration_by_parts}
\end{align}
by Fubini's theorem; so $\mu_u\in\mathrm{BV}_{\mathrm{loc}}((0,+\infty))$ and $D\mu_u$ is the push-forward of $\mu$ under $u$, $D\mu_u=-u_\sharp\mu$ (cf. \cite{Ambrosio2000} 1.70).
By (\ref{coarea_formula}),
\[
D\mu_u\mres((0,+\infty)\setminus N)(A)=-\mu(\{u\in A\}\setminus Z)=
-\int_A\int_{\{u=t\}\setminus Z}\frac{1}{|u^\prime|}\,\frac{d\mathcal{H}^{0}}{\tau}\,dt
\]
for any $\mathscr{L}^1$-measurable set $A$ in $(0,+\infty)$. In light of the above, we may identify 
$D\mu_u^a=D\mu_u\mres((0,+\infty)\setminus N)$ and  $D\mu_u^s=D\mu_{u}\mres N$. The set of atoms of $D\mu_{u}$ is defined by $A:=\{t\in(0,+\infty):D\mu_u(\{t\})\neq 0\}$. For $t>0$, 
\[
D\mu_u(\{t\})=D\mu_u^s(\{t\})=(D\mu_u\mres N)((\{t\})=-u_\sharp\mu(N\cap\{t\})=-\mu(Z\cap\{u=t\})
\]
and this entails {\em (v)}. The monotone function $\mu_u$ is a good representative within its equivalence class and is differentiable $\mathscr{L}^1$-a.e. on $(0,+\infty)$ with derivative given by the density of $D\mu_{u}$ with respect to $\mathscr{L}^1$ by \cite{Ambrosio2000} Theorem 3.28. Item {\em (vi)} follows from (\ref{coarea_formula}) and {\em (iii)}. Item {\em (vii)} follows from {\em (ii)}.\qed

\smallskip

\noindent Let $0<a<b<+\infty$ and $\varrho\geq 0$ be a non-decreasing bounded function on $[a,b]$. Let $\eta\in\{\pm 1\}^2$. We study solutions to the first-order linear ordinary differential equation
\begin{align}
&u^\prime+(1/x+\varrho)u+\lambda=0\text{ a.e. on }(a,b)\text{ with }u(a)=\eta_1\text{ and }u(b)=\eta_2\label{first_order_ode_with_bc}
\end{align}
where $u\in C^{0,1}([a,b])$ and $\lambda\in\mathbb{R}$. In case $\varrho\equiv 0$ on $[a,b]$ we use the notation $u_0$.

\smallskip

\begin{lemma}\label{first_order_ode}
Let $0<a<b<+\infty$ and $\varrho\geq 0$ be a non-decreasing bounded function on $[a,b]$. Let $\eta\in\{\pm 1\}^2$. Then
\begin{itemize}
\item[(i)] there exists a solution $(u,\lambda)$ of (\ref{first_order_ode_with_bc}) with $u\in C^{0,1}([a,b])$ and $\lambda=\lambda_\eta\in\mathbb{R}$;
\item[(ii)] the pair $(u,\lambda)$ in (i) is unique;
\item[(iii)] $\lambda_\eta$ is given by
\begin{align}
&-\lambda_{(1,1)}=\lambda_{(-1,-1)}=m;\,\lambda_{(1,-1)}=-\lambda_{(-1,1)}=\hat{m};\nonumber
\end{align}
\item[(iv)] if $\eta=(1,1)$ or $\eta=(-1,-1)$ then $u$ is uniformly bounded away from zero on $[a,b]$.
\end{itemize}
\end{lemma}

\smallskip

\noindent{\em Proof.}
{\em (i)} For $\eta=(1,1)$ define $u:[a,b]\rightarrow\mathbb{R}$ by
\begin{equation}\label{formula_for_y}
u(t):=\frac{m\int_a^t g\,ds+g(a)}{g(t)}\text{ for }t\in[a,b]
\end{equation}
with $m$ as in (\ref{definition_of_m}). Then $u\in C^{0,1}([a,b])$ and satisfies (\ref{first_order_ode_with_bc}) with $\lambda=-m$. For $\eta=(1,-1)$ set $u=(-\hat{m}\int_a^\cdot g\,ds+g(a))/g$ with $\lambda=\hat{m}$. The cases  $\eta=(-1,-1)$ and $\eta=(-1,1)$ can be dealt with using linearity.
{\em (ii)} We consider the case $\eta=(1,1)$. Suppose that $(u_1,\lambda_1)$ resp. $(u_2,\lambda_2)$ solve  (\ref{first_order_ode_with_bc}). By linearity $u:=u_1-u_2$ solves
\[
u^\prime+(1/x+\varrho)u+\lambda=0\text{ a.e. on }(a,b)
\text{ with }u(a)=u(b)=0
\]
where $\lambda=\lambda_1-\lambda_2$. An integration gives that $u=(-\lambda\int_a^\cdot g\,ds+c)/g$ for some constant $c\in\mathbb{R}$ and the boundary conditions entail that $\lambda=c=0$. The other cases are similar. {\em (iii)} follows as in {\em (i)}. {\em (iv)} If $\eta=(1,1)$ then $u>0$ on $[a,b]$ from (\ref{formula_for_y}) as $m>0$. 
\qed

\smallskip

\noindent{\em The boundary condition $\eta_1\eta_2=-1$.}

\smallskip

\begin{lemma}\label{first_order_ode_with_eta_having_product_minus_one}
Let $0<a<b<+\infty$ and $\varrho\geq 0$ be a non-decreasing bounded function on $[a,b]$. Let $(u,\lambda)$ solve (\ref{first_order_ode_with_bc}) with $\eta=(1,-1)$. Then
\begin{itemize}
\item[(i)] there exists a unique $c\in(a,b)$ with $u(c)=0$;
\item[(ii)] $u^\prime<0$ a.e. on $[a,c]$ and $u$ is strictly decreasing on $[a,c]$;
\item[(iii)] $D\mu_u^s=0$.
\end{itemize}
\end{lemma}

\smallskip

\noindent{\em Proof.}
{\em (i)} We first observe that $u^\prime\leq-\hat{m}<0$ a.e. on $\{u\geq 0\}$ in view of (\ref{first_order_ode_with_bc}). Suppose $u(c_1)=u(c_2)=0$ for some $c_1,c_2\in(a,b)$ with $c_1<c_2$. We may assume that $u\geq 0$ on $[c_1,c_2]$. This contradicts the above observation. Item {\em (ii)} is plain. For any $\mathscr{L}^1$-measurable set $B$ in $(0,+\infty)$, $D\mu_u^s(B)=\mu(\{u\in B\}\cap Z)=0$ using Lemma \ref{properties_of_D_mu} and {\em (ii)}.
\qed

\smallskip

\begin{lemma}\label{derivative_of_distribution_function_for_y_resp_z}
Let $0<a<b<+\infty$ and $\varrho\geq 0$ be a non-decreasing bounded function on $[a,b]$. Let $(u,\lambda)$ solve (\ref{first_order_ode_with_bc}) with $\eta=(1,-1)$. Assume that
\begin{itemize}
\item[(a)] $u$ is differentiable at both $a$ and $b$ and that (\ref{first_order_ode_with_bc}) holds there;
\item[(b)] $u^\prime(a)<0$ and $u^\prime(b)<0$;
\item[(c)] $\varrho$ is differentiable at $a$ and $b$. 
\end{itemize}
Put $v:=-u$. Then
\begin{itemize}
\item[(i)]
$\int_{\{v=1\}\setminus Z_v}\frac{1}{|v^\prime|}\frac{d\mathscr{H}^0}{\tau}
\geq 
\int_{\{u=1\}\setminus Z_u}\frac{1}{|u^\prime|}\frac{d\mathscr{H}^0}{\tau}$;
\item[(ii)] equality holds if and only if $\varrho\equiv 0$ on $[a,b)$.
\end{itemize}
\end{lemma}

\smallskip

\noindent{\em Proof.}
First, $\{u=1\}=\{a\}$ by Lemma \ref{first_order_ode_with_eta_having_product_minus_one}. Further $0<-au^\prime(a)=1+a[\hat{m}+\varrho(a)]$ from (\ref{first_order_ode_with_bc}). On the other hand $\{v=1\}\supset\{b\}$ and $0<bv^\prime(b)=-1+b[\hat{m}-\varrho(b)]$. Thus
\begin{align}
\int_{\{v=1\}\setminus Z_v}\frac{1}{|v^\prime|}\frac{d\mathscr{H}^0}{\tau}-
\int_{\{u=1\}\setminus Z_u}\frac{1}{|u^\prime|}\frac{d\mathscr{H}^0}{\tau}
&\geq\frac{1}{-1+b[\hat{m}-\varrho(b)]}-\frac{1}{1+a[\hat{m}+\varrho(a)]}.\nonumber
\end{align}
By Theorem \ref{reverse_Hermite_Hadamard_inequality}, $0\leq 2+(a-b)\hat{m}+a\varrho(a)+b\varrho(b)$, noting that $\varrho(a)=\varrho(a+)$ in virtue of (c) and similarly at $b$. A rearrangement leads to the inequality. The equality assertion follows from Theorem \ref{reverse_Hermite_Hadamard_inequality}.
\qed

\smallskip

\begin{theorem}\label{comparision_of_derivatives_of_distribution_functions}
Let $0<a<b<+\infty$ and $\varrho\geq 0$ be a non-decreasing bounded function on $[a,b]$. Suppose that $(u,\lambda)$ solves (\ref{first_order_ode_with_bc}) with $\eta=(1,-1)$ and set $v:=-u$. Assume that $u>-1$ on $[a,b)$. 
Then
\begin{itemize}
\item[(i)] $-\mu_v^\prime\geq-\mu_u^\prime$ for $\mathscr{L}^1$-a.e. $t\in(0,1)$;
\item[(ii)] if $\varrho\not\equiv 0$ on $[a,b)$ then there exists $t_0\in(0,1)$ such that $-\mu_v^\prime>-\mu_u^\prime$ for $\mathscr{L}^1$-a.e. $t\in(t_0,1)$;
\item[(iii)] for $t\in[-1,1]$,
\[
\mu_{u_0}(t)=\log\Big\{\frac{-(b-a)t+\sqrt{(b-a)^2t^2+4ab}}{2a}\Big\}
\]
and $\mu_{v_0}=\mu_{u_0}$ on $[-1,1]$;
\end{itemize}
in obvious notation.
\end{theorem}

\smallskip

\noindent{\em Proof.}
{\em (i)} The set 
\[
Y_u:=Z_u\cup\Big(\{u^\prime+(1/x+\varrho)u+\lambda\neq 0\}\setminus Z_{2,u}\Big)\cup\{\varrho\text{ not differentiable}\}\subset[a,b]
\]
(in obvious notation) is a null set in $[a,b]$ and likewise for $Y_v$. By \cite{Ambrosio2000} Lemma 2.95, $\{u=t\}\cap Y_u=\emptyset$ for a.e. $t\in(0,1)$ and likewise for the function $v$. Let $t\in(0,1)$ and assume that $\{u=t\}\cap Y_u=\emptyset$ and $\{v=t\}\cap Y_v=\emptyset$. Put $c:=\max\{u\geq t\}$. Then $c\in(a,b)$, $\{u>t\}=[a,c)$ by Lemma \ref{first_order_ode_with_eta_having_product_minus_one} and $u$ is differentiable at $c$ with $u^\prime(c)<0$. Put $d:=\max\{v\leq t\}=\max\{u\geq-t\}$. As $u$ is continuous on $[a,b]$ it holds that $a<c<d<b$. Moreover, $u^\prime(d)<0$ as $v(d)=t$ and $d\not\in Z_v$. Put $\widetilde{u}:=u/t$ and $\widetilde{v}:=v/t$ on $[c,d]$. Then
\begin{align}
\widetilde{u}^\prime+(1/\tau+\varrho)\widetilde{u}+\hat{m}/t&=0\text{ a.e. on }(c,d)\text{ and }\widetilde{u}(c)=-\widetilde{u}(d)=1;\nonumber\\
\widetilde{v}^\prime+(1/\tau+\varrho)\widetilde{v}-\hat{m}/t&=0\text{ a.e. on }(c,d)\text{ and }-\widetilde{v}(c)=\widetilde{v}(d)=1.\nonumber
\end{align}
By Lemma \ref{derivative_of_distribution_function_for_y_resp_z},
\begin{align}
\int_{\{v=t\}\setminus Z_v}\frac{1}{|v^\prime|}\frac{d\mathscr{H}^0}{\tau}
&\geq\int_{[c,d]\cap\{v=t\}\setminus Z_v}\frac{1}{|v^\prime|}\frac{d\mathscr{H}^0}{\tau}
=(1/t)\int_{[c,d]\cap\{\widetilde{v}=1\}\setminus Z_v}\frac{1}{|\widetilde{v}^\prime|}\frac{d\mathscr{H}^0}{\tau}\nonumber\\
&\geq(1/t)\int_{[c,d]\cap\{\widetilde{u}=1\}\setminus Z_u}\frac{1}{|\widetilde{u}^\prime|}\frac{d\mathscr{H}^0}{\tau}
=\int_{\{u=t\}\setminus Z_u}\frac{1}{|u^\prime|}\frac{d\mathscr{H}^0}{\tau}.\nonumber
\end{align}
By Lemma \ref{properties_of_D_mu},
\[
-\mu_u^\prime(t)=\int_{\{u=t\}\setminus Z_u}\frac{1}{|u^\prime|}\frac{d\mathscr{H}^0}{\tau}
\]
for $\mathscr{L}^1$-a.e. $t\in(0,1)$ and a similar formula holds for $v$. The assertion in {\em (i)} follows.

\smallskip

\noindent{\em (ii)} Assume that $\varrho\not\equiv 0$ on $[a,b)$. Put $\alpha:=\inf\{\varrho>0\}\in[a,b)$. Note that $\max\{v\leq t\}\rightarrow b$ as $t\uparrow 1$ as $v<1$ on $[a,b)$ by assumption. Choose $t_0\in(0,1)$ such that $\max\{v\leq t_0\}>\alpha$. Then for $t>t_0$,
\[
a<\max\{u\geq t\}<\max\{u\geq -t_0\}=\max\{v\leq t_0\}<\max\{v\leq t\}<d;
\]
that is, the interval $[c,d]$ with $c,d$ as described above intersects $(\alpha,b]$. So for $\mathscr{L}^1$-a.e. $t\in(t_0,1)$,
\begin{align}
\int_{\{v=t\}\setminus Z_v}\frac{1}{|v^\prime|}\frac{d\mathscr{H}^0}{\tau}
&>\int_{\{u=t\}\setminus Z_u}\frac{1}{|u^\prime|}\frac{d\mathscr{H}^0}{\tau}.\nonumber
\end{align}
by the equality condition in Lemma \ref{derivative_of_distribution_function_for_y_resp_z}. The conclusion follows from the representation of $\mu_u$ resp. $\mu_v$ in Lemma \ref{properties_of_D_mu}.

\smallskip

\noindent{\em (iii)} A direct computation gives
\[
u_0(\tau)=\frac{1}{b-a}\Big\{-\tau+\frac{ab}{\tau}\Big\}
\]
for $\tau\in[a,b]$; $u_0$ is strictly decreasing on its domain. This leads to the formula in {\em (iii)}. A similar computation gives
\[
\mu_{v_0}(t)=\log\Big\{\frac{2b}{(b-a)t+\sqrt{(b-a)^2t^2+4ab}}\Big\}
\]
for $t\in[-1,1]$. Rationalising the denominator results in the stated equality.
\qed

\smallskip

\begin{corollary}\label{inequality_for_mu_u_vs_mu_v}
Let $0<a<b<+\infty$ and $\varrho\geq 0$ be a non-decreasing bounded function on $[a,b]$. Suppose that $(u,\lambda)$ solves (\ref{first_order_ode_with_bc}) with $\eta=(1,-1)$ and set $v:=-u$. Assume that $u>-1$ on $[a,b)$. Then 
\begin{itemize}
\item[(i)] $\mu_u(t)\leq\mu_v(t)$ for each $t\in(0,1)$;
\item[(ii)] if $\varrho\not\equiv 0$ on $[a,b)$ then $\mu_u(t)<\mu_v(t)$ for each $t\in(0,1)$.
\end{itemize}
\end{corollary}

\smallskip

\noindent{\em Proof.}
{\em (i)} By \cite{Ambrosio2000} Theorem 3.28 and Lemma \ref{first_order_ode_with_eta_having_product_minus_one}, 
\[
\mu_u(t)=\mu_u(t)-\mu_u(1)=-D\mu_u((t,1])=-D\mu_u^a((t,1])-D\mu_u^s((t,1])
=-\int_{(t,1]}\mu_u^\prime\,ds
\]
for each $t\in(0,1)$ as $\mu_u(1)=0$. On the other hand,
\[
\mu_v(t)=\mu_v(1)+(\mu_v(t)-\mu_v(1))
=\mu_v(1)-D\mu_v((t,1])=\mu_v(1)-\int_{(t,1]}\mu_v^\prime\,ds-D\mu_v^s((t,1])
\]
for each $t\in(0,1)$. The claim follows from Theorem \ref{comparision_of_derivatives_of_distribution_functions} noting that $D\mu_v^s((t,1])\leq 0$ as can be seen from Lemma \ref{properties_of_D_mu}. Item {\em (ii)} follows from Theorem \ref{comparision_of_derivatives_of_distribution_functions} {\em (ii)}.
\qed

\smallskip

\begin{corollary}\label{integral_of_solution_of_first_order_ode}
Let $0<a<b<+\infty$ and $\varrho\geq 0$ be a non-decreasing bounded function on $[a,b]$. Suppose that $(u,\lambda)$ solves (\ref{first_order_ode_with_bc}) with $\eta=(1,-1)$. Assume that $u>-1$ on $[a,b)$. Let $\varphi\in C^1((-1,1))$ be an odd strictly increasing function with $\varphi\in L^1((-1,1))$. Then
\begin{itemize}
\item[(i)] $\int_{\{u>0\}}\varphi(u)\,d\mu<+\infty$;
\item[(ii)] $\int_a^b\varphi(u)\,d\mu\leq 0$;
\item[(iii)] equality holds in (ii) if and only if $\varrho\equiv 0$ on $[a,b)$.
\end{itemize}
In particular, 
\begin{itemize}
\item[(iv)] $\int_a^b\frac{u}{\sqrt{1-u^2}}\,d\mu\leq 0$ with equality if and only if $\varrho\equiv 0$ on $[a,b)$.
\end{itemize}
\end{corollary}

\smallskip

\noindent{\em Proof.}
{\em (i)} Put $I:=\{1>u>0\}$. The function $u:I\rightarrow(0,1)$ is $C^{0,1}$ and $u^\prime\leq-\hat{m}$ a.e. on $I$ by Lemma \ref{first_order_ode_with_eta_having_product_minus_one}. It has $C^{0,1}$ inverse $v:(0,1)\rightarrow I$, $v^\prime=1/(u^\prime\circ v)$ and $|v^\prime|\leq 1/\hat{m}$ a.e. on $(0,1)$. By a change of variables, 
\begin{align}
\int_{\{u>0\}}\varphi(u)\,d\mu&=\int_0^1\varphi(v^\prime/v)\,dt\nonumber
\end{align}
from which the claim is apparent.
{\em (ii)} The integral is well-defined because $\varphi(u)^+=\varphi(u)\chi_{\{u>0\}}\in L^1((a,b),\mu)$ by {\em (i)}. By Lemma \ref{first_order_ode_with_eta_having_product_minus_one} the set $\{u=0\}$ consists of a singleton and has $\mu$-measure zero. So
\begin{align}
\int_a^b\varphi(u)\,d\mu&=\int_{\{u>0\}}\varphi(u)\,d\mu+\int_{\{u<0\}}\varphi(u)\,d\mu
=\int_{\{u>0\}}\varphi(u)\,d\mu-\int_{\{v>0\}}\varphi(v)\,d\mu\nonumber
\end{align}
where $v:=-u$ as $\varphi$ is an odd function. We remark that in a similar way to (\ref{integration_by_parts}),
\[
\int_0^1\varphi^\prime\mu_u\,dt=\int_{\{u>0\}}\Big\{\varphi(u)-\varphi(0)\Big\}\,d\mu
=\int_{\{u>0\}}\varphi(u)\,d\mu
\]
using oddness of $\varphi$ and an analogous formula holds with $v$ in place of $u$. Thus we may write
\begin{align}
\int_a^b\varphi(u)\,d\mu&=\int_0^1\varphi^\prime\mu_u\,dt-\int_0^1\varphi^\prime\mu_v\,dt
=\int_0^1\varphi^\prime\Big\{\mu_u-\mu_v\Big\}\,dt\leq 0\nonumber
\end{align}
by Corollary \ref{inequality_for_mu_u_vs_mu_v} as $\varphi^\prime>0$ on $(0,1)$. {\em (iii)} Suppose that $\varrho\not\equiv 0$ on $[a,b)$. Then strict inequality holds in the above by Corollary \ref{inequality_for_mu_u_vs_mu_v}. If $\varrho\equiv 0$ on $[a,b)$ the equality follows from Theorem \ref{comparision_of_derivatives_of_distribution_functions}. {\em (iv)} follows from {\em (ii)} and {\em (iii)} with the particular choice $\varphi:(-1,1)\rightarrow\mathbb{R};t\mapsto t/\sqrt{1-t^2}$.
\qed

\smallskip

\noindent{\em The boundary condition $\eta_1\eta_2=1$.} Let $0<a<b<+\infty$ and  $\varrho\geq 0$ be a non-decreasing bounded function on $[a,b]$. We study solutions of the auxilliary Riccati equation 
\begin{align}
w^\prime+\lambda w^2&=(1/x+\varrho)w\text{ a.e. on }(a,b)
\text{ with }w(a)=w(b)=1;\label{Riccati_equation}
\end{align}
with $w\in C^{0,1}([a,b])$ and $\lambda\in\mathbb{R}$. If $\varrho\equiv 0$ on $[a,b]$ then we write $w_0$ instead of $w$. Suppose $(u,\lambda)$ solves (\ref{first_order_ode_with_bc}) with $\eta=(1,1)$. Then $u>0$ on $[a,b]$ by Lemma \ref{first_order_ode} and we may set $w:=1/u$. Then $(w,-\lambda)$ satisfies (\ref{Riccati_equation}).

\smallskip

\begin{lemma}\label{solution_of_Riccati_equation}
Let $0<a<b<+\infty$ and $\varrho\geq 0$ be a non-decreasing bounded function on $[a,b]$. Then
\begin{itemize}
\item[(i)] there exists a solution $(w,\lambda)$ of (\ref{Riccati_equation}) with $w\in C^{0,1}([a,b])$ and
$\lambda\in\mathbb{R}$;
\item[(ii)] the pair $(w,\lambda)$ in (i) is unique;
\item[(iii)] $\lambda=m$.
\end{itemize}
\end{lemma}

\smallskip

\noindent{\em Proof.}
{\em (i)} Define $w:[a,b]\rightarrow\mathbb{R}$ by
\[
w(t):=\frac{g(t)}{m\int_a^t g\,ds+g(a)}\text{ for }t\in[a,b].
\]
Then $w\in C^{0,1}([a,b])$ and $(w,m)$ satisfies (\ref{Riccati_equation}). {\em (ii)} We claim that $w>0$ on $[a,b]$ for any solution $(w,\lambda)$ of (\ref{Riccati_equation}). For otherwise, $c:=\min\{w=0\}\in(a,b)$. Then $u:=1/w$ on $[a,c)$ satisfies
\[
u^\prime+(\frac{1}{\tau}+\varrho)u-\lambda=0\text{ a.e. on }(a,c)
\text{ and }u(a)=1, u(c-)=+\infty.
\]
Integrating, we obtain
\[
gu-g(a)-\lambda\int_a^\cdot g\,dt=0\text{ on }[a,c)
\]
and this entails the contradiction that $u(c-)<+\infty$. We may now use the uniqueness statement in Lemma \ref{first_order_ode}. {\em (iii)} follows from {\em (ii)} and the particular solution given in {\em (i)}. 
\qed

\smallskip

\noindent We introduce the mapping
\[
\omega:(0,\infty)\times(0,\infty)\rightarrow\mathbb{R};(t,x)\mapsto-(2/t)\coth(x/2).
\]
For $\xi>0$,
\begin{align}
|\omega(t,x)-\omega(t,y)|&\leq\mathrm{cosech}^2[\xi/2](1/t)|x-y|\label{Lipschitz_estimate_for_omega}
\end{align}
for $(t,x),(t,y)\in(0,\infty)\times(\xi,\infty)$ and $\omega$ is locally Lipschitzian in $x$ on $(0,\infty)\times(0,\infty)$ in the sense of \cite{Hale1969} I.3.
Let $0<a<b<+\infty$ and set $\lambda:=A/G>1$. Here, $A=A(a,b)$ stands for the arithmetic mean of $a,b$ as introduced in the previous Section while $G=G(a,b):=\sqrt{|ab|}$ stands for their geometric mean. We refer to the inital value problem
\begin{equation}\label{auxilliary_ode_1}
z^\prime=\omega(t,z)\text{ on }(0,\lambda)
\text{ and }z(1)=\mu((a,b)).
\end{equation}
Define
\[
z_0:(0,\lambda)\rightarrow\mathbb{R};t\mapsto
2\log\Big\{\frac{\lambda+\sqrt{\lambda^2-t^2}}{t}\Big\}.
\]

\smallskip

\begin{lemma}\label{properties_of_w_0}
Let $0<a<b<+\infty$.  Then
\begin{itemize}
\item[(i)] $w_0(\tau)=\frac{2A\tau}{G^2+\tau^2}$ for $\tau\in[a,b]$;
\item[(ii)] $\|w_0\|_\infty=\lambda$;
\item[(iii)] $\mu_{w_0}=z_0$ on $[1,\lambda)$;
\item[(iv)] $z_0$ satisfies (\ref{auxilliary_ode_1}) and this solution is unique;
\item[(v)] $\int_{\{w_0=1\}}\frac{1}{|w_0^\prime|}\,\frac{d\mathscr{H}^0}{\tau}=2\coth(\mu((a,b))/2)$;
\item[(vi)] $\int_a^b\frac{1}{\sqrt{w_0^2-1}}\,\frac{dx}{x}=\pi$. 
\end{itemize}
\end{lemma}

\smallskip

\noindent{\em Proof.} {\em (i)} follows as in the proof of Lemma \ref{solution_of_Riccati_equation} with $g(t)=t$ while {\em (ii)} follows by calculus. {\em (iii)} follows by solving the quadratic equation $t\tau^2-2A\tau+G^2t=0$ for   $\tau$ with $t\in(0,\lambda)$. Uniqueness in {\em (iv)} follows from \cite{Hale1969} Theorem 3.1 as $\omega$ is locally Lipschitzian with respect to $x$ in $(0,\infty)\times(0,\infty)$. For {\em (v)} note that $|aw_0^\prime(a)|=1-a/A$ and $|bw_0^\prime(b)|=b/A-1$ and 
\[
2\coth(\mu((a,b))/2)=2(a+b)/(b-a).
\]
{\em (vi)} We may write
\begin{align}
\int_a^b\frac{1}{\sqrt{w_0^2-1}}\,\frac{d\tau}{\tau}
&=\int_a^b\frac{ab+\tau^2}{\sqrt{(a+b)^2\tau^2-(ab+\tau^2)^2}}
\frac{d\tau}{\tau}=\int_a^b\frac{ab+\tau^2}{\sqrt{(\tau^2-a^2)(b^2-\tau^2)}}
\frac{d\tau}{\tau}.\nonumber
\end{align}
The substitution $s=\tau^2$ followed by the Euler substitution (cf. \cite{Gradshteynetal1965} 2.251)
$\sqrt{(s-a^2)(b^2-s)}=t(s-a^2)$ gives
\begin{align}
\int_a^b\frac{1}{\sqrt{w_0^2-1}}\,\frac{d\tau}{\tau}
&=\int_0^\infty\frac{1}{1+t^2}+\frac{ab}{b^2+a^2t^2}\,dt=\pi.\nonumber
\end{align}
\qed

\smallskip

\begin{lemma}\label{monotonicity_of_rational_functions}
Let $0<a<b<+\infty$. Then
\begin{itemize}
\item[(i)] for $y>a$ the function $x\mapsto\frac{by-ax}{(y-a)(b-x)}$ is strictly increasing on $(-\infty,b]$;
\item[(ii)] the function $y\mapsto\frac{(b-a)y}{(y-a)(b-y)}$ is strictly increasing on $[G,b]$;
\item[(iii)] for $x<b$ the function $y\mapsto\frac{by-ax}{(y-a)(b-x)}$ is strictly decreasing on $[a,+\infty)$
\end{itemize}
\end{lemma}

\smallskip

\noindent{\em Proof.} The proof is an exercise in calculus.
\qed

\smallskip

\begin{lemma}\label{lemma_on_derivative_of_mu_v}
Let $0<a<b<+\infty$ and $\varrho\geq 0$ be a non-decreasing bounded function on $[a,b]$. Let $(w,\lambda)$ solve (\ref{Riccati_equation}). Assume
\begin{itemize}
\item[(i)] $w$ is differentiable at both $a$ and $b$ and that (\ref{Riccati_equation}) holds there;
\item[(ii)] $w^\prime(a)>0$ and $w^\prime(b)<0$;
\item[(iii)] $w>1$ on $(a,b)$;
\item[(iv)] $\varrho$ is differentiable at $a$ and $b$.
\end{itemize}
Then
\[
\int_{\{w=1\}\setminus Z_w}\frac{1}{|w^\prime|}\frac{d\mathscr{H}^0}{\tau}\geq 2\coth(\mu((a,b))/2)
\]
with equality if and only if $\varrho\equiv 0$ on $[a,b)$.
\end{lemma}

\smallskip

\noindent{\em Proof.}
At the end-points $x=a,b$ the condition {\em (i)} entails that $w^\prime+m-\varrho=1/x=w_0^\prime+m_0$ so that 
\begin{align}
w^\prime-w_0^\prime&=m_0-m+\varrho\text{ at }x=a,b.\label{equation_for_derivatives}
\end{align}
We consider the four cases
\begin{itemize}
\item[(a)] $w^\prime(a)\geq w_0^\prime(a)\text{ and }w^\prime(b)\geq w_0^\prime(b)$;
\item[(b)] $w^\prime(a)\geq w_0^\prime(a)\text{ and }w^\prime(b)\leq w_0^\prime(b)$;
\item[(c)] $w^\prime(a)\leq w_0^\prime(a)\text{ and }w^\prime(b)\geq w_0^\prime(b)$;
\item[(d)] $w^\prime(a)\leq w_0^\prime(a)\text{ and }w^\prime(b)\leq w_0^\prime(b)$;
\end{itemize}
in turn. 

\smallskip

\noindent (a)
Condition (a) together with (\ref{equation_for_derivatives}) means that $m_0-m+\varrho(a)\geq 0$;
that is, $m-\varrho(a)\leq m_0$. By {\em (i)} and {\em (ii)}, $bm-b\varrho(b)-1=-bw^\prime(b)> 0$; or $m-\varrho(b)>1/b$. Therefore,
\[
0<1/b<m-\varrho(b)\leq m-\varrho(a)\leq 1/A
\]
by (\ref{formula_for_m_0}). Put $x:=1/(m-\varrho(b))$ and $y:=1/(m-\varrho(a))$. Then
\[
a<A\leq y\leq x<b.
\]
We write
\begin{align}
aw^\prime(a)&=-(m-\varrho(a))a+1=-(1/y)a+1>0;\nonumber\\
bw^\prime(b)&=-(m-\varrho(b))b+1=-(1/x)b+1<0.\nonumber
\end{align}
Making use of assumption {\em (iii)},
\begin{align}
\int_{\{w=1\}\setminus Z_w}\frac{1}{|w^\prime|}\frac{d\mathscr{H}^0}{x}
&=\frac{1}{-(1/y)a+1}-\frac{1}{-(1/x)b+1}
=\frac{by-ax}{(y-a)(b-x)}.\nonumber
\end{align}
By Lemma \ref{monotonicity_of_rational_functions} {\em (i)} then {\em (ii)},
\begin{align}
\int_{\{w=1\}}\frac{1}{|w^\prime|}\frac{d\mathscr{H}^0}{x}
&\geq\frac{(b-a)y}{(y-a)(b-y)}
\geq\frac{(b-a)A}{(A-a)(b-A)}
=2\frac{a+b}{b-a}
=2\coth(\mu((a,b))/2).\nonumber
\end{align}
If equality holds then $\varrho(a)=\varrho(b)$ and $\varrho$ is constant on $[a,b]$. By Theorem \ref{upper_bound_for_m} we conclude that $\varrho\equiv 0$ on $[a,b)$.

\smallskip

\noindent (b) Condition (b) together with (\ref{equation_for_derivatives}) entails that $0\leq m_0-m+\varrho(a)$ and 
$0\leq -m_0+m-\varrho(b)$ whence $0\leq\varrho(a)-\varrho(b)$ upon adding; so $\varrho$ is constant on the interval $[a,b]$ by monotonicity. Define $x$ and $y$ as above. Then $x=y$ and $y\geq A$. The result now follows in a similar way to case (a). 

\smallskip

\noindent (c) In this case,
\[
\frac{1}{aw^\prime(a)}-\frac{1}{bw^\prime(b)}\geq\frac{1}{aw_0^\prime(a)}-\frac{1}{bw_0^\prime(b)}
=2\coth(\mu((a,b))/2)
\]
by Lemma \ref{properties_of_w_0}.  If equality holds then $w^\prime(b)=w_0^\prime(b)$ so that $m_0-m+\varrho(b)=0$ and $\varrho$ vanishes on $[a,b]$ by Theorem \ref{upper_bound_for_m}.

\smallskip

\noindent (d) Condition (d) together with (\ref{equation_for_derivatives}) means that $m_0-m+\varrho(b)\leq 0$;
that is, $m\geq \varrho(b)+m_0$. On the other hand, by Theorem \ref{upper_bound_for_m}, $m\leq\varrho(b)+m_0$. In consequence, $m=\varrho(b)+m_0$. It then follows that $\varrho\equiv 0$ on $[a,b]$ by Theorem \ref{upper_bound_for_m}. Now use Lemma \ref{properties_of_w_0}.\qed

\smallskip

\begin{lemma}\label{On_convex_decreasing_function}
Let $\phi:(0,+\infty)\rightarrow(0,+\infty)$ be a convex non-increasing function with $\inf_{(0,+\infty)}\phi>0$. Let $\Lambda$ be an at most countably infinite index set and $(x_h)_{h\in\Lambda}$ a sequence of points in $(0,+\infty)$ with $\sum_{h\in\Lambda}x_h<+\infty$. Then
\[
\sum_{h\in\Lambda}\phi(x_h)\geq \phi(\sum_{h\in\Lambda} x_h)
\]
and the left-hand side takes the value $+\infty$ in case $\Lambda$ is countably infinite and is otherwise finite.
\end{lemma}

\smallskip

\noindent{\em Proof.}
Suppose $0<x_1<x_2<+\infty$. By convexity $\phi(x_1)+\phi(x_2)\geq 2\phi(\frac{x_1+x_2}{2})\geq \phi(x_1+x_2)$ as $\phi$ is non-increasing. The result for finite $\Lambda$ follows by induction. \qed

\smallskip

\begin{theorem}\label{derivative_of_Ricatti_distribution_function}
Let $0<a<b<+\infty$ and $\varrho\geq 0$ be a non-decreasing bounded function on $[a,b]$. Let $(w,\lambda)$ solve (\ref{Riccati_equation}). Assume that $w>1$ on $(a,b)$. Then
\begin{itemize}
\item[(i)] for $\mathscr{L}^1$-a.e. $t\in(1,\|w\|_\infty)$,
\begin{equation}\label{inequality_for_mu_v}
-\mu_w^\prime\geq(2/t)\coth((1/2)\mu_w);
\end{equation}
\item[(ii)] if $\varrho\not\equiv 0$ on $[a,b)$ then there exists $t_0\in(1,\|w\|_\infty)$ such that strict inequality holds in (\ref{inequality_for_mu_v}) for $\mathscr{L}^1$-a.e. $t\in(1,t_0)$.
\end{itemize}
\end{theorem}

\smallskip

\noindent{\em Proof.}
{\em (i)} The set 
\[
Y_w:=Z_w\cup\Big(\{w^\prime+mw^2\neq(1/x+\varrho)w\}\setminus Z_{2,w}\Big)\cup\{\varrho\text{ not differentiable}\}\subset[a,b]
\]
is a null set in $[a,b]$. By \cite{Ambrosio2000} Lemma 2.95, $\{w=t\}\cap Y_w=\emptyset$ for a.e. $t>1$. Let $t\in(1,\|w\|_\infty)$ and assume that $\{w=t\}\cap Y_w=\emptyset$. We write $\{w>t\}=\bigcup_{h\in\Lambda}I_h$
where $\Lambda$ is an at most countably infinite index set and $(I_h)_{h\in\Lambda}$ are disjoint non-empty  well-separated open intervals in $(a,b)$. The term well-separated means that for each $h\in\Lambda$, $\inf_{k\in\Lambda\setminus\{h\}}d(I_h,I_k)>0$. This follows from the fact that $w^\prime\neq 0$ on $\partial I_h$ for each $h\in\Lambda$. Put $\widetilde{w}:=w/t$ on $\overline{\{w>t\}}$ so
\begin{align}
\widetilde{w}^\prime+(mt)\widetilde{w}^2&=(1/x+\varrho)\widetilde{w}\text{ a.e. on }\{w>t\}\text{ and }\widetilde{w}=1\text{ on }\{w=t\}.\nonumber
\end{align}
We use the fact that the mapping  $\phi:(0,+\infty)\rightarrow(0,+\infty);t\mapsto\coth t$ satisfies the hypotheses of Lemma \ref{On_convex_decreasing_function}. By Lemmas \ref{lemma_on_derivative_of_mu_v} and \ref{On_convex_decreasing_function},
\begin{align}
(0,+\infty]\ni\int_{\{w=t\}\setminus Z_w}\frac{1}{|w^\prime|}\frac{d\mathscr{H}^0}{x}
&=(1/t)\int_{\{\widetilde{w}=1\}}\frac{1}{|\widetilde{w}^\prime|}\frac{d\mathscr{H}^0}{\tau}\nonumber\\
&=(1/t)\sum_{h\in\Lambda}\int_{\partial I_h}\frac{1}{|\widetilde{w}^\prime|}\frac{d\mathscr{H}^0}{\tau}\nonumber\\
&\geq(2/t)\sum_{h\in\Lambda}\coth((1/2)\mu(I_h))\nonumber\\
&\geq(2/t)\coth((1/2)\sum_{h\in\Lambda}\mu(I_h))\nonumber\\
&=(2/t)\coth((1/2)\mu(\{w>t\})))=(2/t)\coth((1/2)\mu_w(t)).\nonumber
\end{align}
The statement now follows from Lemma \ref{properties_of_D_mu}.

\smallskip

\noindent{\em (ii)} Suppose that $\varrho\not\equiv 0$ on $[a,b)$. Put $\alpha:=\min\{\varrho>0\}\in[a,b)$. Now  that $\{w>t\}\uparrow(a,b)$ as $t\downarrow 1$ as $w>1$ on $(a,b)$. Choose $t_0\in(1,\|w\|_\infty)$ such that $\{w>t_0\}\cap(\alpha,b)\neq\emptyset$. Then for each $t\in(1,t_0)$ there exists $h\in\Lambda$ such that $\varrho\not\equiv 0$ on $I_h$. The statement then follows by Lemma \ref{lemma_on_derivative_of_mu_v}.
\qed

\smallskip

\begin{lemma}\label{measurability_lemma}
Let $\emptyset\neq S\subset\mathbb{R}$ be bounded and suppose $S$ has the property that for each $s\in S$ there exists $\delta>0$ such that $[s,s+\delta)\subset S$. Then $S$ is $\mathscr{L}^1$-measurable and $|S|>0$.
\end{lemma}

\smallskip

\noindent{\em Proof.} For each $s\in S$ put $t_s:=\inf\{t>s:t\not\in S\}$. Then $s<t_s<+\infty$, $[s,t_s)\subset S$ and $t_s\not\in S$. Define
\[
\mathscr{C}:=\Big\{[s,t]:s\in S\text{ and }t\in(s,t_s)\Big\}.
\]
Then $\mathscr{C}$ is a Vitali cover of $S$ (see \cite{Carothers2000} Chapter 16 for example). By Vitali's Covering Theorem (cf. \cite{Carothers2000} Theorem 16.27) there exists an at most countably infinite subset $\Lambda\subset\mathscr{C}$ consisting of pairwise disjoint intervals such that
\[
|S\setminus\bigcup_{I\in\Lambda}I|=0.
\]
Note that $I\subset S$ for each $I\in\Lambda$. Consequently, $S=\bigcup_{I\in\Lambda}I\cup N$ where $N$ is an $\mathscr{L}^1$-null set and hence $S$ is $\mathscr{L}^1$-measurable. The positivity assertion is clear.
\qed

\smallskip

\begin{theorem}\label{distribution_function_inequality_for_Ricatti_equation}
Let $0<a<b<+\infty$ and $\varrho\geq 0$ be a non-decreasing bounded function on $[a,b]$. Let $(w,\lambda)$ solve (\ref{Riccati_equation}). Assume that $w>1$ on $(a,b)$. Put $T:=\min\{\|w_0\|_\infty,\|w\|_\infty\}>1$. Then 
\begin{itemize}
\item[(i)] $\mu_w(t)\leq\mu_{w_0}(t)$ for each $t\in[1,T)$;
\item[(ii)] $\|w\|_\infty\leq\|w_0\|_\infty$;
\item[(iii)] if $\varrho\not\equiv 0$ on $[a,b)$ then there exists $t_0\in(1,\|w\|_\infty)$ such that $\mu_w(t)<\mu_{w_0}(t)$ for each $t\in(1,t_0)$.
\end{itemize}
\end{theorem}

\smallskip

\noindent{\em Proof.} {\em (i)} We adapt the proof of \cite{Hale1969} Theorem I.6.1. The assumption entails that $\mu_w(1)=\mu_{w_0}(1)=\mu((a,b))$. Suppose for a contradiction that $\mu_w(t)>\mu_{w_0}(t)$ for some $t\in(1,T)$. 

\smallskip

\noindent For $\varepsilon>0$ consider the initial value problem
\begin{equation}\label{auxilliary_ode_2}
z^\prime=\omega(t,z)+\varepsilon
\text{ and }z(1)=\mu((a,b))+\varepsilon
\end{equation}
on $(0,T)$. Choose $\upsilon\in(0,1)$ and $\tau\in(t,T)$. By \cite{Hale1969} Lemma I.3.1 there exists $\varepsilon_0>0$ such that for each $0\leq\varepsilon<\varepsilon_0$ (\ref{auxilliary_ode_2}) has a continuously differentiable solution $z_\varepsilon$ defined on $[\upsilon,\tau]$ and this solution is unique by \cite{Hale1969} Theorem I.3.1. Moreover, the sequence $(z_\varepsilon)_{0<\varepsilon<\varepsilon_0}$ converges uniformly to $z_0$ on $[\upsilon,\tau]$. 

\smallskip

\noindent Given $0<\varepsilon<\eta<\varepsilon_0$ it holds that $z_0\leq z_\varepsilon\leq z_\eta$ on $[1,\tau]$  by \cite{Hale1969} Theorem I.6.1. Note for example that $z_0^\prime\leq\omega(\cdot,z_0)+\varepsilon$ on $(1,\tau)$. In fact, $(z_\varepsilon)_{0<\varepsilon<\varepsilon_0}$ decreases strictly to $z_0$ on $(1,\tau)$. For if, say, $z_0(s)=z_\varepsilon(s)$ for some $s\in(1,\tau)$ then $z_\varepsilon^\prime(s)=\omega(s,z_\varepsilon(s))+\varepsilon>\omega(s,z_0(s))=z_0^\prime(s)$ by (\ref{auxilliary_ode_2}); while on the other hand $z_\varepsilon^\prime(s)\leq z_0^\prime(s)$ by considering the left-derivative at $s$ and using the fact that $z_\varepsilon\geq z_0$ on $[1,\tau]$. This contradicts the strict inequality.
 
 \smallskip

\noindent Choose $\varepsilon_1\in(0,\varepsilon_0)$ such that $z_\varepsilon(t)<\mu_w(t)$ for each $0<\varepsilon<\varepsilon_1$. Now $\mu_w$ is right-continuous and strictly decreasing as $\mu_w(t)-\mu_w(s)=-\mu(\{s<w\leq t\})<0$ for $1\leq s<t<\|w\|_\infty$ by continuity of $w$. So the set $\{z_\varepsilon<\mu_w\}\cap(1,t)$ is open and non-empty in $(0,+\infty)$ for each $\varepsilon\in(0,\varepsilon_1)$. Thus there exists a unique $s_\varepsilon\in[1,t)$ such that
\[
\mu_w>z_\varepsilon\text{ on }(s_\varepsilon,t]
\text{ and }\mu_w(s_\varepsilon)=z_\varepsilon(s_\varepsilon)
\]
for each $\varepsilon\in(0,\varepsilon_1)$. As $z_\varepsilon(1)>\mu((a,b))$ it holds that each $s_\varepsilon>1$. Note that $1<s_\varepsilon<s_\eta$ whenever $0<\varepsilon<\eta$ as $(z_\varepsilon)_{_{0<\varepsilon<\varepsilon_0}}$ decreases strictly to $z_0$ as $\varepsilon\downarrow 0$. 

\smallskip

\noindent Define
\[
S:=\Big\{s_\varepsilon:0<\varepsilon<\varepsilon_1\Big\}\subset(1,t).
\]
We claim that for each $s\in S$ there exists $\delta>0$ such that $[s,s+\delta)\subset S$. This entails that $S$ is $\mathscr{L}^1$-measurable with positive $\mathscr{L}^1$-measure by Lemma \ref{measurability_lemma}.
 
\smallskip

\noindent Suppose $s=s_\varepsilon\in S$ for some $\varepsilon\in(0,\varepsilon_1)$ and put $z:=z_\varepsilon(s)=\mu_w(s)$. Put $k:=\mathrm{cosech}^2(z_0(t)/2)$. For $0\leq\zeta<\eta<\varepsilon_1$ define
\[
\Omega_{\zeta,\eta}:=\Big\{(u,y)\in\mathbb{R}^2:u\in(0,t)\text{ and }z_\zeta(u)<y<z_\eta(u)\Big\}
\]
and note that this is an open set in $\mathbb{R}^2$. We remark that for each $(u,y)\in\Omega_{\zeta,\eta}$ there exists a unique $\nu\in(\zeta,\eta)$ such that $y=z_\nu(u)$. Given $r>0$ with $s+r<t$ set
\[
Q=Q_r:=\Big\{(u,y)\in\mathbb{R}^2:s\leq u<s+r\text{ and }|y-z|<\|z_\varepsilon-z\|_{C([s,s+r])}\Big\}.
\]
Choose $r\in(0,t-s)$ and $\varepsilon_2\in(\varepsilon,\varepsilon_1)$ such that
\begin{itemize}
\item[(a)]  $Q_r\subset\Omega_{0,\varepsilon_1}$;
\item[(b)] $\|z_\varepsilon-z\|_{C([s,s+r])}<s\varepsilon/(2k)$;
\item[(c)] $\sup_{\eta\in(\varepsilon,\varepsilon_2)}\|z_\eta-z\|_{C([s,s+r])}\leq\|z_\varepsilon-z\|_{C([s,s+r])}$;
\item[(d)] $z_\eta<\mu_w$ on $[s+r,t]$ for each $\eta\in(\varepsilon,\varepsilon_2)$.
\end{itemize}
We can find $\delta\in(0,r)$ such that $z_\varepsilon<\mu_w<z_{\varepsilon_2}$ on $(s,s+\delta)$ as $z_{\varepsilon_2}(s)>z$; in other words, the graph of $\mu_w$ restricted to $(s,s+\delta)$ is contained in $\Omega_{\varepsilon,\varepsilon_2}$.

\smallskip

\noindent Let $u\in(s,s+\delta)$. Then $\mu_w(u)=z_\eta(u)$ for some $\eta\in(\varepsilon,\varepsilon_2)$ as above. We claim that $u=s_\eta$ so that $u\in S$. This implies in turn that $[s,s+\delta)\subset S$. Suppose for a contradiction that $z_\eta\not<\mu_w$ on $(u,t]$. Then there exists $v\in(u,t]$ such that $\mu_w(v)=z_\eta(v)$. In view of condition (d), $v\in(u,s+r)$. By \cite{Ambrosio2000} Theorem 3.28 and Theorem \ref{derivative_of_Ricatti_distribution_function},
\begin{align}
\mu_w(v)-\mu_w(u)&=D\mu_w((u,v])
=D\mu_w^a((u,v])+D\mu_w^s((u,v])
\nonumber\\
&\leq D\mu_w^a((u,v])
=\int_u^v\mu_w^\prime\,d\tau
\leq\int_u^v\omega(\cdot,\mu_w)\,d\tau.\nonumber
\end{align}
On the other hand,
\begin{align}
z_\eta(v)-z_\eta(u)&=\int_u^v z_\eta^\prime\,d\tau
=\int_u^v\omega(\cdot,z_\eta)\,d\tau+\eta(v-u).\nonumber
\end{align}
We derive that
\begin{align}
\varepsilon(v-u)&\leq\eta(v-u)
\leq\int_u^v\Big\{\omega(\cdot,\mu_w)-\omega(\cdot,z_\eta)\Big\}\,d\tau
\leq k\int_u^v|\mu_w-z_\eta|\,d\mu\nonumber
\end{align}
using the estimate (\ref{Lipschitz_estimate_for_omega}). Thus
\begin{align}
\varepsilon&\leq  k\frac{1}{v-u}\int_u^v|\mu_w-z_\eta|\,d\mu\nonumber\\
&\leq(k/s)\|\mu_w-z_\eta\|_{C([u,v])}\nonumber\\
&\leq(k/s)\Big\{\|\mu_w-z\|_{C([s,s+r])}+\|z_\eta-z\|_{C([s,s+r])}\Big\}\nonumber\\
&\leq(2k/s)\|z_\varepsilon-z\|_{C([s,s+r])}
<\varepsilon\nonumber
\end{align}
by (b) and (c) giving rise to the desired contradiction.

\smallskip

\noindent By Theorem \ref{derivative_of_Ricatti_distribution_function}, $\mu_w^\prime\leq\omega(\cdot,\mu_w)$ for $\mathscr{L}^1$-a.e. $t\in S$. Choose $s\in S$ such that $\mu_w$ is differentiable at $s$ and the latter inequality holds at $s$. Let $\varepsilon\in(0,\varepsilon_1)$ such that $s=s_\varepsilon$. For any $u\in(s,t)$,
\[
\mu_w(u)-\mu_w(s)>z_\varepsilon(u)-z_\varepsilon(s).
\]
We deduce that $\mu_w^\prime(s)\geq z_\varepsilon^\prime(s)$. But then
\[
\mu_w^\prime(s)\geq z_\varepsilon^\prime(s)=\omega(s,z_\epsilon(s))+\varepsilon>\omega(s,\mu_w(s)).
\]
This strict inequality holds on a set of full measure in $S$. This contradicts Theorem \ref{derivative_of_Ricatti_distribution_function}.

\smallskip

\noindent{\em (ii)} Use the fact that $\|w\|_\infty=\sup\{t>0:\mu_w(t)>0\}$.

\smallskip

\noindent{\em (iii)} Assume that $\varrho\not\equiv 0$ on $[a,b)$. Let $t_0\in(1,\|w\|_\infty)$ be as in Lemma \ref{derivative_of_Ricatti_distribution_function}. Then for $t\in(1,t_0)$,
\begin{align}
\mu_w(t)-\mu_w(1)&=D\mu_w((1,t])=D\mu_w^a((1,t])+D\mu_w^s((1,t])\leq D\mu_w^a((1,t])\nonumber\\
&=\int_{(1,t]}\mu_w^\prime\,ds<\int_{(1,t]}\omega(s,\mu_w)\,ds\leq\int_{(1,t]}\omega(s,\mu_{w_0})\,ds=\mu_{w_0}(t)-\mu_{w_0}(1)\nonumber
\end{align}
by Theorem \ref{derivative_of_Ricatti_distribution_function}, Lemma \ref{properties_of_w_0} and the inequality in {\em (i)}.
\qed

\smallskip

\begin{corollary}\label{integral_of_w}
Let $0<a<b<+\infty$ and $\varrho\geq 0$ be a non-decreasing bounded function on $[a,b]$. Suppose that $(w,\lambda)$ solves (\ref{Riccati_equation}). Assume that $w>1$ on $(a,b)$. Let $0\leq\varphi\in C^1((1,+\infty))$ be strictly decreasing with $\int_a^b\varphi(w_0)\,d\mu<+\infty$. Then
\begin{itemize}
\item[(i)] $\int_a^b\varphi(w)\,d\mu\geq\int_a^b\varphi(w_0)\,d\mu$;
\item[(ii)] equality holds in (i) if and only if $\varrho\equiv 0$ on $[a,b)$.
\end{itemize}
In particular, 
\begin{itemize}
\item[(iii)] $\int_a^b\frac{1}{\sqrt{w^2-1}}\,d\mu\geq\pi$ with equality if and only if $\varrho\equiv 0$ on $[a,b)$.
\end{itemize}
\end{corollary}

\smallskip

\noindent{\em Proof.}
{\em (i)} Let $\varphi\geq 0$ be a decreasing function on $(1,+\infty)$ which is piecewise $C^1$. Suppose that $\varphi(1+)<+\infty$. By Tonelli's Theorem,
\begin{align}
\int_{[1,+\infty)}\varphi^\prime\mu_w\,ds
&=\int_{[1,+\infty)}\varphi^\prime\Big\{\int_{(a,b)}\chi_{\{w>s\}}\,d\mu\Big\}\,ds\nonumber\\
&=\int_{(a,b)}\Big\{\int_{[1,+\infty)}\varphi^\prime\chi_{\{w>s\}}\,ds\Big\}\,d\mu\nonumber\\
&=\int_{(a,b)}\Big\{\varphi(w)-\varphi(1)\Big\}\,d\mu
=\int_{(a,b)}\varphi(w)\,d\mu-\varphi(1)\mu((a,b))\nonumber
\end{align}
and a similar identity holds for $\mu_{w_0}$. By Theorem \ref{distribution_function_inequality_for_Ricatti_equation}, $\int_a^b\varphi(w)\,d\mu\geq\int_a^b\varphi(w_0)\,d\mu$. Now suppose that $0\leq\varphi\in C^1((1,+\infty))$ is strictly decreasing with $\int_a^b\varphi(w_0)\,d\mu<+\infty$. The inequality holds for the truncated function $\varphi\wedge n$ for each $n\in\mathbb{N}$. An application of the monotone convergence theorem establishes the result for $\varphi$.

\smallskip

\noindent{\em (ii)} Suppose that equality holds in {\em (i)}. For $c\in(1,+\infty)$ put $\varphi_1:=\varphi\vee\varphi(c)-\varphi(c)$ and $\varphi_2:=\varphi\wedge\varphi(c)$. By {\em (i)} we deduce
\[
\int_a^b\varphi_2(w)\,d\mu=\int_a^b\varphi_2(w_0)\,d\mu;
\]
and hence by the above that
\[
\int_{[c,+\infty)}\varphi^\prime\Big\{\mu_w-\mu_{w_0}\Big\}\,ds=0.
\]
This means that $\mu_w=\mu_{w_0}$ on $(c,+\infty)$ and hence on $(1,+\infty)$. By Theorem \ref{distribution_function_inequality_for_Ricatti_equation} we conclude that $\varrho\equiv 0$ on $[a,b)$. {\em (iii)} flows from {\em (i)} and {\em (ii)} noting that the function $\varphi:(1,+\infty)\rightarrow\mathbb{R};t\mapsto 1/\sqrt{t^2-1}$ satisfies the integral condition by Lemma \ref{properties_of_w_0}.
\qed

\smallskip

\noindent{\em The case $a=0$.} Let $0<b<+\infty$ and $\varrho\geq 0$ be a non-decreasing bounded function on $[0,b]$. We study solutions to the first-order linear ordinary differential equation
\begin{align}
&u^\prime+(1/x+\varrho)u+\lambda=0\text{ a.e. on }(0,b)\text{ with }u(0)=0\text{ and }u(b)=1\label{first_order_ode_with_bc_and_a_zero}
\end{align}
where $u\in C^{0,1}([0,b])$ and $\lambda\in\mathbb{R}$. If $\varrho\equiv 0$ on $[0,b]$ then we write $u_0$ instead of $u$.

\smallskip

\begin{lemma}\label{first_order_ode_with_a_zero}
Let $0<b<+\infty$ and $\varrho\geq 0$ be a non-decreasing bounded function on $[0,b]$. Then
\begin{itemize}
\item[(i)] there exists a solution $(u,\lambda)$ of (\ref{first_order_ode_with_bc_and_a_zero}) with $u\in C^{0,1}([0,b])$ and $\lambda\in\mathbb{R}$;
\item[(ii)] $\lambda$ is given by $\lambda=-g(b)/G(b)$ where $G:=\int_0^\cdot g\,ds$;
\item[(iii)] the pair $(u,\lambda)$ in (i) is unique;
\item[(iv)] $u>0$ on $(0,b]$.
\end{itemize}
\end{lemma}

\smallskip

\noindent{\em Proof.}
{\em (i)} The function $u:[a,b]\rightarrow\mathbb{R}$ given by
\begin{equation}\label{representation_of_u_in_case_a_is_0}
u=\frac{g(b)}{G(b)}\frac{G}{g}
\end{equation}
on $[0,b]$ solves (\ref{first_order_ode_with_bc_and_a_zero}) with $\lambda$ as in {\em (ii)}.
{\em (iii)} Suppose that $(u_1,\lambda_1)$ resp. $(u_2,\lambda_2)$ solve  (\ref{first_order_ode_with_bc_and_a_zero}). By linearity $u:=u_1-u_2$ solves
\[
u^\prime+(1/x+\varrho)u+\lambda=0\text{ a.e. on }(0,b)
\text{ with }u(0)=u(b)=0
\]
where $\lambda=\lambda_1-\lambda_2$. An integration gives that $u=(-\lambda G+c)/g$ for some constant $c\in\mathbb{R}$ and the boundary conditions entail that $\lambda=c=0$. {\em (iv)} follows from the formula (\ref{representation_of_u_in_case_a_is_0}) and unicity. \qed

\smallskip

\begin{lemma}\label{convexity_lemma}
Suppose $-\infty<a<b<+\infty$ and that $\phi:[a,b]\rightarrow\mathbb{R}$ is convex. Suppose that there exists $\xi\in(a,b)$ such that
\[
\phi(\xi)=\frac{b-\xi}{b-a}\phi(a)+\frac{\xi-a}{b-a}\phi(b).
\]
Then 
\[
\phi(c)=\frac{b-c}{b-a}\phi(a)+\frac{c-a}{b-a}\phi(b)
\]
for each $c\in[a,b]$.
\end{lemma}

\smallskip

\noindent{\em Proof.}
Let $c\in(\xi,b)$. By monotonicity of chords,
\[
\frac{\phi(\xi)-\phi(a)}{\xi-a}\leq\frac{\phi(c)-\phi(\xi)}{c-\xi}
\]
so
\begin{align}
\phi(c)&\geq\frac{c-a}{\xi-a}\phi(\xi)-\frac{c-\xi}{\xi-a}\phi(a)\nonumber\\
&=\frac{c-a}{\xi-a}\Big\{\frac{b-\xi}{b-a}\phi(a)+\frac{\xi-a}{b-a}\phi(b)\Big\}
-\frac{c-\xi}{\xi-a}\phi(a)\nonumber\\
&=\frac{b-c}{b-a}\phi(a)+\frac{c-a}{b-a}\phi(b)\nonumber
\end{align}
and equality follows. The case $c\in(a,\xi)$ is similar.
\qed

\smallskip

\begin{lemma}\label{comparision_result_for_u_in_case_a_is_zero}
Let $0<b<+\infty$ and $\varrho\geq 0$ be a non-decreasing bounded function on $[0,b]$. Let $(u,\lambda)$ satisfy (\ref{first_order_ode_with_bc_and_a_zero}). Then 
\begin{itemize}
\item[(i)] $u\geq u_0$ on $[0,b]$;
\item[(ii)] if $\varrho\not\equiv 0$ on $[0,b)$ then $u>u_0$ on $(0,b)$.
\end{itemize}
\end{lemma}

\smallskip

\noindent{\em Proof.}
{\em (i)} The mapping $G:[0,b]\rightarrow[0,G(b)]$ is a bijection with inverse $G^{-1}$. Define $\eta:[0,G(b)]\rightarrow\mathbb{R}$ via $\eta:=(tg)\circ G^{-1}$. Then
\[
\eta^\prime=\frac{(tg)^\prime}{g}\circ G^{-1}
=(2+t\varrho)\circ G^{-1}
\]
a.e. on $(0,G(b))$ so $\eta^\prime$ is non-decreasing there. This means that $\eta$ is convex on $[0,G(b)]$. In particular, $\eta(s)\leq[\eta(G(b))/G(b)]s$ for each $s\in[0,G(b)]$. For $t\in[0,b]$ put $s:=G(t)$ to obtain
$tg(t)\leq(bg(b)/G(b))G(t)$. A rearrangement gives $u\geq u_0$ on $[0,b]$ noting that $u_0:[0,b]\rightarrow\mathbb{R};t\mapsto t/b$. {\em (ii)} Assume $\varrho\not\equiv 0$ on $[0,b)$. Suppose that $u(c)=u_0(c)$ for some $c\in(0,b)$. Then $\eta(G(c))=[\eta(G(b))/G(b)]G(c)$. By Lemma \ref{convexity_lemma}, $\eta^\prime=0$ on $(0,G(b))$. This implies that $\varrho\equiv 0$ on $[0,b)$. \qed

\smallskip

\begin{lemma}\label{integral_of_solution_of_first_order_ode_with_a_equal_to_0}
Let $0<b<+\infty$. Then $\int_0^b\frac{u_0}{\sqrt{1-u_0^2}}\,d\mu=\pi/2$. 
\end{lemma}

\smallskip

\noindent{\em Proof.} The integral is elementary as $u_0(t)=t/b$ for $t\in[0,b]$.
\qed

\section{Proof of Main Results}

\smallskip

\begin{lemma}\label{lemma_on_cones}
Let $x\in H$ and $v$ be a unit vector in $\mathbb{R}^2$ such that the pair $\{x,v\}$ forms a positively oriented orthogonal basis for $\mathbb{R}^2$. Put $b:=(\tau,0)$ where $|x|=\tau$ and $\gamma:=\theta(x)\in(0,\pi)$. Let $\alpha\in(0,\pi/2)$ such that
\[
\frac{\langle v,x-b\rangle}{|x-b|}=\cos\alpha.
\]
Then
\begin{itemize}
\item[(i)] $C(x,v,\alpha)\cap H\cap\overline{C}(0,e_1,\gamma)=\emptyset$;
\item[(ii)] for any $y\in C(x,v,\alpha)\cap H\setminus\overline{B}(0,\tau)$ the line segment $[b,y]$ intersects $\mathbb{S}^1_\tau$ outside the closed cone $\overline{C}(0,e_1,\gamma)$.
\end{itemize}
\end{lemma}

\smallskip

\noindent  We point out that $C(0,e_1,\gamma)$ is the open cone with vertex $0$ and axis $e_1$ which contains the point $x$ on its boundary. We note that $\cos\alpha\in(0,1)$ because
\begin{equation}\label{v_dot_b}
\langle v,x-b\rangle=-\langle v,b\rangle=-\langle(1/\tau)Ox,b\rangle=-\langle Op,e_1\rangle=\langle x,O^\star e_1\rangle=\langle x,e_2\rangle>0
\end{equation}
and if $|x-b|=\langle v,x-b\rangle$ then $b=x-\lambda v$ for some $\lambda\in\mathbb{R}$ and hence $x_1=\langle e_1,x\rangle=\tau$ and $x_2=0$.

\smallskip

\noindent{\em Proof.}
{\em (i)} For $\omega\in\mathbb{S}^1$ define the open half-space
\[
H_{\omega}:=\{y\in\mathbb{R}^2:\langle y,\omega\rangle>0\}.
\]
We claim that $C(x,v,\alpha)\subset H_v$. For given $y\in C(x,v,\alpha)$,
\[
\langle y,v\rangle=\langle y-x,v\rangle>|y-x|\cos\alpha>0.
\]
On the other hand, it holds that $\overline{C}(0,e_1,\gamma)\cap H\subset\overline{H}_{-v}$. This establishes {\em (i)}. 

\smallskip

\noindent{\em (ii)}  By some trigonometry $\gamma=2\alpha$. Suppose that $\omega$ is a unit vector in $C(b,-e_1,\pi/2-\alpha)$. Then $\lambda:=\langle\omega,e_1\rangle<\cos\alpha$ since upon rewriting the membership condition for $C(b,-e_1,\pi/2-\alpha)$ we obtain the quadratic inequality
\[
\lambda^2-2\cos^2\alpha\lambda+\cos\gamma>0.
\]
For $\omega$ a unit vector in $\overline{C}(0,e_1,\gamma)$ the opposite inequality $\langle\omega,e_1\rangle\geq\cos\alpha$ holds. This shows that
\[
C(b,-e_1,\pi/2-\alpha)\cap\overline{C}(0,e_1,\gamma)\cap\mathbb{S}^1_\tau=\emptyset.
\]

\smallskip

\noindent The set $C(x,v,\alpha)\cap H$ is contained in the open convex cone $C(b,-e_1,\pi/2-\alpha)$. Suppose $y\in C(x,v,\alpha)\cap H\setminus\overline{B}(0,\tau)$. Then the line segment $[b,y]$ is contained in $C(b,-e_1,\pi/2-\alpha)\cup\{b\}$. Now the set $C(b,-e_1,\pi/2-\alpha)\cap\mathbb{S}^1_\tau$ disconnects $C(b,-e_1,\pi/2-\alpha)\cup\{b\}$. This entails that $(b,y]\cap C(b,-e_1,\pi/2-\alpha)\cap\mathbb{S}^1_\tau\neq\emptyset$.  The foregoing paragraph entails that $(b,y]\cap\overline{C}(0,e_1,\gamma)\cap\mathbb{S}^1_\tau=\emptyset$. This establishes the result.
\qed

\smallskip

\begin{lemma}\label{non_convexity_criterion}
Let $E$ be an open set in $\mathbb{R}^2$ such that $M:=\partial E$ is a $C^{1,1}$ hypersurface in $\mathbb{R}^2$. Assume that $E\setminus\{0\}=E^{sc}$. Suppose 
\begin{itemize}
\item[(i)] $x\in(M\setminus\{0\})\cap H$;
\item[(ii)] $\sin(\sigma(x))=-1$.
\end{itemize}
Then $E$ is not convex. 
\end{lemma}

\smallskip

\noindent{\em Proof.}
Let $\gamma_1:I\rightarrow M$ be a $C^{1,1}$ parametrisation of $M$ in a neighbourhood of $x$  with $\gamma_1(0)=x$ as above. As $\sin(\sigma(x))=-1$, $n(x)$ and hence $n_1(0)$ point in the direction of $x$. Put $v:=-t_1(0)=-t(x)$. We may write
\[
\gamma_1(s)=\gamma_1(0)+st_1(0)+R_1(s)=x-sv+R_1(s)
\]
for $s\in I$ where $R_1(s)=s\int_0^1\dot{\gamma}_1(ts)-\dot{\gamma}_1(0)\,dt$ and we can find a finite positive constant $K$ such that $|R_1(s)|\leq Ks^2$ on a symmetric open interval $I_0$ about $0$ with $I_0\subset\subset I$. Then 
\begin{align}
\frac{\langle\gamma_1(s)-x,v\rangle}{|\gamma_1(s)-x|}&=\frac{\langle-sv+R_1,v\rangle}{|-sv+R_1|}
=\frac{1-\langle(R_1/s),v\rangle}{|v-R_1/s|}\rightarrow 1\nonumber
\end{align}
as $s\uparrow 0$. Let $\alpha$ be as in Lemma \ref{lemma_on_cones} with $x$ and $v$ as just mentioned. The above estimate entails that $\gamma_1(s)\in C(x,v,\alpha)$ for small $s<0$. By (\ref{cos_of_sigma}) and Lemma \ref{properties_of_section_of_boundary_of_E} the function $r_1$ is non-increasing on $I$. In particular, $r_1(s)\geq r_1(0)=|x|=:\tau$ for $I\ni s<0$ and $\gamma_1(s)\not\in B(0,\tau)$. 

\smallskip

\noindent Choose $\delta_1>0$ such that $\gamma_1(s)\in C(x,v,\alpha)\cap H$ for each $s\in[-\delta_1,0)$. Put $\beta:=\inf\{s\in[-\delta_1,0]:r_1(s)=\tau\}$. Suppose first that $\beta\in[-\delta_1,0)$. Then $E$ is not convex (see Lemma \ref{property_of_E_sc}). Now suppose that $\beta=0$. Let $\gamma$ be as in Lemma \ref{lemma_on_cones}. Then the open circular arc $\mathbb{S}^1_\tau\setminus \overline{C}(0,e_1,\gamma)$ does not intersect $\overline{E}$: for otherwise, $M$ intersects $\mathbb{S}^1_\tau\setminus \overline{C}(0,e_1,\gamma)$ and $\beta<0$ bearing in mind Lemma \ref{property_of_E_sc}. Choose $s\in[-\delta_1,0)$. Then  the points $b$ and $\gamma_1(s)$ lie in $\overline{E}$. But by Lemma \ref{lemma_on_cones} the line segment $[b,\gamma_1(s)]$ intersects $\mathbb{S}^1_\tau$ in $\mathbb{S}^1_\tau\setminus \overline{C}(0,e_1,\gamma)$. Let $c\in[b,\gamma_1(s)]\cap\mathbb{S}^1_\tau$. Then $c\not\in\overline{E}$. This shows that $\overline{E}$ is not convex. But if $E$ is convex then $\overline{E}$ is convex. Therefore $E$ is not convex.
\qed

\smallskip

\begin{theorem}\label{Omega_empty}
Let $f$ be as in (\ref{form_of_f}) where $h:[0,+\infty)\rightarrow\mathbb{R}$ is a non-decreasing convex function. Given $v>0$ let $E$ be a bounded minimiser of (\ref{isoperimetric_problem}). Assume that $E$ is open, $M:=\partial E$ is a $C^{1,1}$ hypersurface in $\mathbb{R}^2$ and $E\setminus\{0\}=E^{sc}$. Put
\begin{equation}\label{definition_of_R}
R:=\inf\{\varrho>0\}\in[0,+\infty).
\end{equation}
Then $\Omega\cap(R,+\infty)=\emptyset$ with $\Omega$ as in (\ref{definition_of_Omega}).
\end{theorem}

\smallskip

\noindent{\em Proof.}
Suppose that $\Omega\cap(R,+\infty)\neq\emptyset$. As $\Omega$ is open in $(0,+\infty)$ by Lemma \ref{Omega_is_open} we may write $\Omega$ as a countable union of disjoint open intervals in $(0,+\infty)$. By a suitable choice of one of these intervals we may assume that $\Omega=(a,b)$ for some $0\leq a<b<+\infty$ and that $\Omega\cap(R,+\infty)\neq\emptyset$. Let us assume for the time being that $a>0$. Note that $[a,b]\subset\pi(M)$ and $\cos\sigma$ vanishes on $M_a\cup M_b$. 

\smallskip

\noindent Let $u:\Omega\rightarrow[-1,1]$ be as in (\ref{definition_of_y}). Then $u$ has a continuous extension to $[a,b]$ and $u=\pm 1$ at $\tau=a,b$. This may be seen as follows. For $\tau\in(a,b)$ the set $M_\tau\cap\overline{H}$ consists of a singleton by Lemma \ref{properties_of_section_of_boundary_of_E}. The limit $x:=\lim_{\tau\downarrow a}M_\tau\cap\overline{H}\in\mathbb{S}^1_a\cap\overline{H}$ exists as $M$ is $C^1$. There exists a $C^{1,1}$ parametrisation $\gamma_1:I\rightarrow M$ with $\gamma_1(0)=x$ as above. By (\ref{cos_of_sigma}) and Lemma \ref{properties_of_section_of_boundary_of_E}, $r_1$ is decreasing on $I$. So $r_1>a$ on $I\cap\{s<0\}$ for otherwise the $C^1$ property fails at $x$. It follows that $\gamma_1=\gamma\circ r_1$ and $\sigma_1=\sigma\circ\gamma\circ r_1$ on $I\cap\{s<0\}$. Thus $\sin(\sigma\circ\gamma)\circ r_1=\sin\sigma_1$ on $I\cap\{s<0\}$. Now the function $\sin\sigma_1$ is continuous on $I$. So $u\rightarrow\sin\sigma_1(0)\in\{\pm 1\}$ as $\tau\downarrow a$. Put $\eta_1:=u(a)$ and $\eta_2:=u(b)$. 

\smallskip

\noindent Let us consider the case $\eta=(\eta_1,\eta_2)=(1,1)$. According to Theorem \ref{constant_weighted_mean_curvature} the generalised (mean) curvature is constant $\mathscr{H}^1$-a.e. on $M$ with value $-\lambda$, say. Note that $u<1$ on $(a,b)$ for otherwise $\cos(\sigma\circ\gamma)$ vanishes at some point in $(a,b)$ bearing in mind Lemma \ref{properties_of_section_of_boundary_of_E}. By Theorem \ref{ode_for_y} the pair $(u,\lambda)$ satisfies (\ref{first_order_ode_with_bc}) with $\eta=(1,1)$. By Lemma \ref{first_order_ode}, $u>0$ on $[a,b]$. Put $w:=1/u$. Then $(w,-\lambda)$ satisfies (\ref{Riccati_equation}) and $w>1$ on $(a,b)$. By Lemma \ref{derivative_of_theta},
\[
\theta_2(b)-\theta_2(a)=\int_a^b\theta_2^\prime\,d\tau
=-\int_a^b\frac{u}{\sqrt{1-u^2}}\,\frac{d\tau}{\tau}
=-\int_a^b\frac{1}{\sqrt{w^2-1}}\,\frac{d\tau}{\tau}.
\]
By Corollary \ref{integral_of_w}, $|\theta_2(b)-\theta_2(a)|>\pi$. But this contradicts the definition of $\theta_2$ in (\ref{definition_of_theta}) as $\theta_2$ takes values in $(0,\pi)$ on $(a,b)$. If $\eta=(-1,-1)$ then $\lambda>0$ by Lemma \ref{first_order_ode}; this contradicts Lemma \ref{upper_bound_for_lambda}.

\smallskip

\noindent Now let us consider the case $\eta=(-1,1)$. Using the same formula as above, $\theta_2(b)-\theta_2(a)<0$ by Corollary \ref{integral_of_solution_of_first_order_ode}. This means that $\theta_2(a)\in(0,\pi]$. As before the limit $x:=\lim_{\tau\downarrow a}M_\tau\cap\overline{H}\in\mathbb{S}^1_a\cap\overline{H}$ exists as $M$ is $C^1$. Using a local parametrisation it can be seen that $\theta_2(a)=\theta(x)$ and $\sin(\sigma(x))=-1$. If $\theta_2(a)\in(0,\pi)$ then $E$ is not convex by Lemma \ref{non_convexity_criterion}. This contradicts Theorem \ref{E_is_convex}. Note that we may assume that $\theta_2(a)\in(0,\pi)$. For otherwise, $\langle\gamma,e_2\rangle<0$ for $\tau>a$ near $a$, contradicting the definition of $\gamma$ (\ref{definition_of_gamma}). If $\eta=(1,-1)$ then $\lambda>0$ by Lemma \ref{first_order_ode} and this contradicts Lemma \ref{upper_bound_for_lambda} as before.

\smallskip

\noindent Suppose finally that $a=0$. By Lemma \ref{behaviour_at_0}, $u(0)=0$ and $u(b)=\pm 1$. Suppose $u(b)=1$. Again employing the formula above, $\theta_2(b)-\theta_2(0)<-\pi/2$ by Lemma \ref{comparision_result_for_u_in_case_a_is_zero}, the fact that the function $\phi:(0,1)\rightarrow\mathbb{R};t\mapsto t/\sqrt{1-t^2}$ is strictly increasing and Lemma \ref{integral_of_solution_of_first_order_ode_with_a_equal_to_0}. This means that $\theta_2(0)>\pi/2$. This contradicts the $C^1$ property at $0\in M$. If $u(b)=-1$ then then $\lambda>0$ by Lemma \ref{first_order_ode} giving a contradiction.
\qed

\smallskip

\begin{lemma}\label{boundary_of_E_as_union_of_circles}
Let $f$ be as in (\ref{form_of_f}) where $h:[0,+\infty)\rightarrow\mathbb{R}$ is a non-decreasing convex function. Let $v>0$.
\begin{itemize}
\item[(i)] Let $E$ be a bounded minimiser of (\ref{isoperimetric_problem}). Assume that $E$ is open, $M:=\partial E$ is a $C^{1,1}$ hypersurface in $\mathbb{R}^2$ and $E\setminus\{0\}=E^{sc}$. Then for any $r>0$ with $r\geq R$, $M\setminus\overline{B}(0,r)$ consists of a finite union of disjoint centred circles.
\item[(ii)] There exists a minimiser $E$ of (\ref{isoperimetric_problem}) such that $\partial E$ consists of a countable union of disjoint centred circles whose radii accumulate at $0$ if at all.
\end{itemize}
\end{lemma}

\smallskip

\noindent{\em Proof.} {\em (i)} First observe that
\begin{align}
\emptyset\neq\pi(M)&=\Big[\pi(M)\cap[0,r]\Big]\cup
\Big[\pi(M)\cap(r,+\infty)\Big]\setminus\Omega\nonumber
\end{align}
by Lemma \ref{Omega_empty}. We assume that the latter member is non-empty. By definition of $\Omega$, $\cos\sigma=0$ on $M\cap A((r,+\infty))$. Let $\tau\in\pi(M)\cap(r,+\infty)$. We claim that $M_\tau=\mathbb{S}^1_\tau$. Suppose for a contradiction that $M_\tau\neq\mathbb{S}^1_\tau$. By Lemma \ref{property_of_E_sc}, $M_\tau$ is the union of two closed  spherical arcs in $\mathbb{S}^1_\tau$. Let $x$ be a point on the boundary of one of these spherical arcs relative to $\mathbb{S}^1_\tau$. There exists a $C^{1,1}$ parametrisation $\gamma_1:I\rightarrow M$ of $M$ in a neighbourhood of $x$ with $\gamma_1(0)=x$ as before. By shrinking $I$ if necessary we may assume that $\gamma_1(I)\subset  A((r,+\infty))$ as $\tau>r$. By (\ref{cos_of_sigma}), $\dot{r}_1=0$ on $I$ as $\cos\sigma_1=0$ on $I$ because $\cos\sigma=0$ on $M\cap A((r,+\infty))$; that is, $r_1$ is constant on $I$. This means that $\gamma_1(I)\subset\mathbb{S}^1_\tau$. As the function $\sin\sigma_1$ is continuous on $I$ it takes the value $\pm 1$ there. By (\ref{sin_of_sigma}), $r_1\dot{\theta}_1=\sin\sigma_1=\pm 1$ on $I$. This means that $\theta_1$ is either strictly decreasing or strictly increasing on $I$. This entails that the point $x$ is not a boundary point of $M_\tau$ in $\mathbb{S}^1_\tau$ and this proves the claim.

\smallskip

\noindent It follows from these considerations that $M\setminus\overline{B}(0,r)$ consists of a finite union of disjoint centred circles. Note that $f\geq e^{h(0)}=:c>0$ on $\mathbb{R}^2$. As a result, $+\infty>P_f(E)\geq cP(E)$ and in particular the relative perimeter $P(E,\mathbb{R}^2\setminus\overline{B}(0,r))<+\infty$. This explains why $M\setminus\overline{B}(0,r)$ comprises only finitely many circles.

\smallskip

\noindent{\em (ii)} Let $E$ be a bounded minimiser of (\ref{isoperimetric_problem}) such that $E$ is open, $M:=\partial E$ is a $C^{1,1}$ hypersurface in $\mathbb{R}^2$ and $E\setminus\{0\}=E^{sc}$ as in Theorem \ref{M_is_spherical_cap_symmetric}. Assume that $R>0$. By {\em (i)}, $M\setminus\overline{B}(0,R)$ consists of a finite union of disjoint centred circles. We claim that only one of the possibilities
\begin{align}
M_R&=\emptyset,\,M_R=\mathbb{S}^1_R,\,M_R=\{Re_1\}\text{ or }M_R=\{-Re_1\}\label{options_for_M_R}
\end{align}
holds. To prove this suppose that $M_R\neq\emptyset$ and $M_R\neq\mathbb{S}^1_R$. Bearing in mind Lemma \ref{property_of_E_sc} we may choose $x\in M_R$ such that $x$ lies on the boundary of $M_R$ relative to $\mathbb{S}^1_R$. Assume that $x\in H$. Let $\gamma_1:I\rightarrow M$ be a local parametrisation of $M$ with $\gamma_1(0)=x$ with the usual conventions. We first notice that $\cos(\sigma(x))=0$ for otherwise we obtain a contradiction to Theorem \ref{Omega_empty}. As $r_1$ is decreasing on $I$ and $x$ is a relative boundary point it holds that $r_1<R$ on $I^+:=I\cap\{s>0\}$. As $M\setminus\overline{\Lambda_1}$ is open in $M$ we may suppose that $\gamma_1(I^+)\subset M\setminus\overline{\Lambda_1}$. According to Theorem \ref{constant_weighted_mean_curvature} the curvature $k$ of $\gamma_1(I^+)\cap B(0,R)$ is a.e. constant as $\varrho$ vanishes on $(0,R)$. Hence $\gamma_1(I^+)\cap B(0,R)$ consists of a line or circular arc. The fact that $\cos(\sigma(x))=0$ means that $\gamma_1(I^+)\cap B(0,R)$ cannot be a line. So $\gamma_1(I^+)\cap B(0,R)$ is an open arc of a circle $C$ containing $x$ in its closure with centre on the line-segment $[0,x]$ and radius $r\in(0,R)$. By considering a local parametrisation, it can be seen that $C\cap B(0,R)\subset M$. But this contradicts the fact that $E\setminus\{0\}=E^{sc}$. In summary, $M_R\subset\{\pm Re_1\}$. Finally note that if $M_R=\{\pm Re_1\}$ then $M_R=\mathbb{S}^1_R$ by Lemma \ref{property_of_E_sc}. This establishes (\ref{options_for_M_R}).

\smallskip

\noindent Suppose that $M_R=\emptyset$. As both sets $M$ and $\mathbb{S}^1_R$ are compact, $d(M,\mathbb{S}^1_R)>0$. Assume first that $\mathbb{S}^1_R\subset E$. Put $F:=B(0,R)\setminus E$ and suppose $F\neq\emptyset$. Then $F$ is a set of finite perimeter, $F\subset\subset B(0,R)$ and $P(F)=P(E,B(0,R))$. Let $B$ be a centred ball with $|B|=|F|$. By the classical isoperimetric inequality, $P(B)\leq P(F)$. Define $E_1:=(\mathbb{R}^2\setminus B)\cap(B(0,R)\cup E)$. Then $V_f(E_1)=V_f(E)$ and $P_f(E_1)\leq P_f(E)$. That is, $E_1$ is a minimiser of (\ref{isoperimetric_problem}) such that $\partial E_1$ consists of a finite union of disjoint centred circles. Now suppose that $\mathbb{S}^1_R\subset\mathbb{R}^2\setminus\overline{E}$. In like fashion we may redefine $E$ via $E_1:=B\cup(E\setminus\overline{B}(0,R))$ with $B$ a centred ball in $B(0,R)$. The remaining cases in (\ref{options_for_M_R}) can be dealt with in a similar way. The upshot of this argument is that there exists a minimiser of (\ref{isoperimetric_problem}) whose boundary $M$ consists of a finite union of disjoint centred circles in case $R>0$. 

\smallskip

\noindent The assertion for $R=0$ follows from {\em (i)}.
\qed

\smallskip

\begin{lemma}\label{on_the_J_function}
Suppose that the function $J:[0,+\infty)\rightarrow[0,+\infty)$ is continuous non-decreasing and $J(0)=0$. Let $N\in\mathbb{N}\cup\{+\infty\}$ and $\{t_h:h=0,\ldots,2N+1\}$ a sequence of points in $[0,+\infty)$ with
\[
t_0>t_1>\cdots>t_{2h}>t_{2h+1}>\cdots\geq 0.
\]
Then 
\[
+\infty\geq\sum_{h=0}^{2N+1}J(t_h)\geq J(\sum_{h=0}^{2N+1}(-1)^{h}t_h).
\]
\end{lemma}

\smallskip

\noindent{\em Proof.}
We suppose that $N=+\infty$. The series $\sum_{h=0}^\infty(-1)^{h}t_h$ converges by the alternating series test. For each $n\in\mathbb{N}$,
\[
\sum_{h=0}^{2n+1}(-1)^{h}t_h\leq t_0
\]
and the same inequality holds for the infinite sum. As in Step 2 in \cite{Bettaetal2008} Theorem 2.1,
\[
+\infty\geq\sum_{h=0}^\infty J(t_h)
\geq J(t_0)
\geq J(\sum_{h=0}^\infty(-1)^{h}t_h)
\]
as $J$ is non-decreasing.
\qed

\smallskip

\noindent{\em Proof of Theorem \ref{main_theorem}.}
There exists a minimiser $E$ of (\ref{isoperimetric_problem}) with the property that $\partial E$ consists of a countable union of disjoint centred circles whose radii accumulate at $0$ if at all according to Lemma \ref{boundary_of_E_as_union_of_circles}. As such we may write
\[
E=\bigcup_{h=0}^NA((a_{2h+1},a_{2h}))
\]
where $N\in\mathbb{N}\cup\{+\infty\}$ and $+\infty>a_0>a_1>\cdots>a_{2h}>a_{2h+1}>\cdots\geq 0$. Define
\begin{align}
&\mathtt{f}:[0,+\infty)\rightarrow\mathbb{R};t\mapsto e^{h(t)};\nonumber\\
&g:[0,+\infty)\rightarrow\mathbb{R};t\mapsto t\mathtt{f}(t);\nonumber\\
&G:[0,+\infty)\rightarrow\mathbb{R};t\mapsto\int_0^t g\,d\tau.\nonumber
\end{align}
Then $G:[0,+\infty)\rightarrow[0,+\infty)$ is a bijection with inverse $G^{-1}$.  Define the strictly increasing function
\begin{align}
&J:[0,+\infty)\rightarrow\mathbb{R};t\mapsto g\circ G^{-1}.\nonumber
\end{align}
With $\{a_0,a_1,\ldots\}$ as above put $t_h:=G(a_h)$ for $h=0,\ldots,2N+1$. Then $+\infty>t_0>t_1>\cdots>t_{2h}>t_{2h+1}>\cdots\geq 0$. Put $B:=B(0,r)$ where $r:=G^{-1}(v/2\pi)$ so that $V_f(B)=v$. Note that
\[
v=V_f(E)=2\pi\sum_{h=0}^N\Big\{G(a_{2h})-G(a_{2h+1})\Big\}=2\pi\sum_{h=0}^{2N+1}(-1)^{h}t_h.
\]
By Lemma \ref{on_the_J_function},
\[
P_f(E)=2\pi\sum_{h=0}^{2N+1}g(a_h)=2\pi\sum_{h=0}^{2N+1}J(t_h)
\geq 2\pi J(\sum_{h=0}^{2N+1}(-1)^h t_h)=2\pi J(v/2\pi)=P_f(B).
\]
\qed

\smallskip

\noindent{\em Proof of Theorem \ref{uniqueness_theorem}.}
Let $v>0$ and $E$ be a minimiser for (\ref{isoperimetric_problem}). Then $E$ is essentially bounded by Theorem \ref{first_properties_of_isoperimetric_set}. By Theorem \ref{M_is_spherical_cap_symmetric}
there exists an $\mathscr{L}^2$-measurable set $\widetilde{E}$ with the properties
\begin{itemize}
\item[(a)] $\widetilde{E}$ is a minimiser of (\ref{isoperimetric_problem});
\item[(b)] $L_{\widetilde{E}}=L_E$ a.e. on $(0,+\infty)$;
\item[(c)] $\widetilde{E}$ is open, bounded and has $C^{1,1}$ boundary;
\item[(d)] $\widetilde{E}\setminus\{0\}=\widetilde{E}^{sc}$.
\end{itemize}
{\em (i)} Suppose that $0<v\leq v_0$ so that $R>0$. Choose $r\in(0,R]$ such that $V(B(0,r))=V(E)=v$. Suppose that $\widetilde{E}\setminus\overline{B}(0,R)\neq\emptyset$. By Lemma \ref{boundary_of_E_as_union_of_circles} there exists $t>R$ such that $\mathbb{S}^1_t\subset M$. As $g$ is strictly increasing, $g(t)>g(r)$. So $P_f(E)=P_f(\widetilde{E})\geq\pi g(t)>\pi g(r)=P_f(B(0,r))$. This contradicts the fact that $E$ is a minimiser for (\ref{isoperimetric_problem}). So $\widetilde{E}\subset\overline{B}(0,R)$ and $L_{\widetilde{E}}=0$ on $(R,+\infty)$. By property (b), $|E\setminus\overline{B}(0,R)|=0$. By the uniqueness property in the classical isoperimetric theorem (see for example \cite{Fusco2004} Theorem 4.11) the set $E$ is equivalent to a ball $B$ in $\overline{B}(0,R)$.

\smallskip

\noindent{\em (ii)} With $r>0$ as before, $V(B(0,r))=V(E)=v>v_0=V(B(0,R))$ so $r>R$. If $\widetilde{E}\setminus\overline{B}(0,r)\neq\emptyset$ we derive a contradiction in the same way as above. Consequently, $\widetilde{E}=B:=B(0,r)$. Thus, $L_E=L_B$ a.e. on $(0,+\infty)$; in particular, $|E\setminus B|=0$. This entails that $E$ is equivalent to $B$.
\qed

\end{document}